\documentclass[twoside,a4paper,reqno,nobysame,initials]{amsart}

\usepackage{xr-hyper,xparse,ifthen,amssymb,graphicx}
\PassOptionsToPackage{pdfauthor={Vladimir V. Kisil},%
    pdftitle={Cross-Toeplitz Operators on the Fock--Segal--Bargmann Spaces
  and Two-Sided   Convolutions  on the Heisenberg Group},%
    pdfsubject={mathematics},%
    backref=page,%
  pdfkeywords={symmetries}
  colorlinks=true,%
  bookmarks=true,%
  backref=page,%
  pagebackref=true
}{hyperref}
  \AtBeginDocument{\hypersetup{breaklinks=true,%
      colorlinks=true,%
      bookmarks=true,%
      backref=page,%
      pagebackref=true
}}
\IfFileExists{eulervm.sty}{%
\usepackage{eulervm}}{}
\usepackage{hyperref}
   \PassOptionsToPackage{british}{babel}
\newtheorem{thm}{Theorem}[section]
    \newtheorem{prop}[thm]{Proposition}
    \newtheorem{lem}[thm]{Lemma}
    \newtheorem{cor}[thm]{Corollary}
   \theoremstyle{definition}

    \newtheorem{defn}[thm]{Definition}

    \newtheorem{example}[thm]{Example}

   \theoremstyle{remark}
    \newtheorem{rem}[thm]{Remark}

\makeatother
\usepackage[backrefs,msc-links
]{amsrefs}
\BibSpecAlias{incollection}{inproceedings}
\BibSpec{article}{%
    +{}  {\PrintAuthors}                {author}
    +{.} { }                            {title}
    +{.} { }                            {part}
    +{:} { }                            {subtitle}
    +{.} { \PrintContributions}         {contribution}
    +{.} { \PrintPartials}              {partial}
    +{.} { \emph}                       {journal}
    +{,} { \textbf}                            {volume}
    +{}  { \parenthesize}               {number}
    +{:} {}                             {pages}
    +{,} { \PrintDateB}                 {date}
    +{,} { }                            {status}
    +{.} { \PrintTranslation}           {translation}
    +{.} { Reprinted in \PrintReprint}  {reprint}
    +{.} { }                            {note}
    +{.} {}                             {transition}
}

\BibSpec{partial}{%
    +{}  {}                             {part}
    +{:} { }                            {subtitle}
    +{.} { \PrintContributions}         {contribution}
    +{.} { \emph}                       {journal}
    +{,} { \textbf}                            {volume}
    +{}  { \parenthesize}               {number}
    +{:} {}                             {pages}
    +{,} { \PrintDateB}                 {date}
}

\BibSpec{book}{%
    +{}  {\PrintPrimary}                {transition}
    +{.} { \emph}                       {title}
    +{.} { }                            {part}
    +{:} { \emph}                       {subtitle}
    +{.} { }                            {series}
    +{,} { \voltext}                    {volume}
    +{.} { Edited by \PrintNameList}    {editor}
    +{.} { Translated by \PrintNameList}{translator}
    +{.} { \PrintContributions}         {contribution}
    +{.} { }                            {publisher}
    +{.} { }                            {organization}
    +{,} { }                            {address}
    +{,} { \PrintEdition}               {edition}
    +{,} { \PrintDateB}                 {date}
    +{.} { }                            {note}
    +{.} {}                             {transition}
    +{.} { \PrintTranslation}           {translation}
    +{.} { Reprinted in \PrintReprint}  {reprint}
    +{.} {}                             {transition}
}

\BibSpec{collection.article}{%
    +{}  {\PrintAuthors}                {author}
    +{.} { }                            {title}
    +{.} { }                            {part}
    +{:} { }                            {subtitle}
    +{.} { \PrintContributions}         {contribution}
    +{.} { \PrintConference}            {conference}
    +{.} { \PrintBook}                  {book}
    +{.} { In \PrintEditorsA}              {editor}
    +{.} { \emph}                         {booktitle}
    +{,} { pages~}                      {pages}
    +{,} { }                            {publisher}
    +{,} { }                            {organization}
    +{,} { }                            {address}
    +{,} { \PrintDateB}                 {date}
    +{.} { \PrintTranslation}           {translation}
    +{.} { Reprinted in \PrintReprint}  {reprint}
    +{.} { }                            {note}
    +{.} {}                             {transition}
}

 \BibSpec{inproceedings}{%
     +{}  {\PrintAuthors}                {author}
     +{.} { }                            {title}
     +{.} { }                            {part}
     +{:} { }                            {subtitle}
     +{.} { \PrintContributions}         {contribution}
     +{.} { \PrintConference}            {conference}
     +{.} { \PrintBook}                  {book}
     +{.} { In \PrintEditorsA}              {editor}
     +{.} {  \emph}                         {booktitle}
     +{,} { pages~}                      {pages}
     +{,} { }                            {publisher}
     +{,} { }                            {organization}
     +{,} { }                            {address}
     +{,} { \PrintDateB}                 {date}
     +{.} { \PrintTranslation}           {translation}
     +{.} { Reprinted in \PrintReprint}  {reprint}
     +{.} { }                            {note}
     +{.} {}                             {transition}
 }

\BibSpec{conference}{%
    +{}  {}                        {title}
    +{}  {\PrintConferenceDetails} {transition}
}

\BibSpec{innerbook}{%
    +{.} { \emph}                       {title}
    +{.} { }                            {part}
    +{:} { \emph}                       {subtitle}
    +{.} { }                            {series}
    +{,} { \voltext}                    {volume}
    +{.} { Edited by \PrintNameList}    {editor}
    +{.} { Translated by \PrintNameList}{translator}
    +{.} { \PrintContributions}         {contribution}
    +{.} { }                            {publisher}
    +{.} { }                            {organization}
    +{,} { }                            {address}
    +{,} { \PrintEdition}               {edition}
    +{,} { \PrintDateB}                 {date}
    +{.} { }                            {note}
    +{.} {}                             {transition}
}

\BibSpec{report}{%
    +{}  {\PrintPrimary}                {transition}
    +{.} { \emph}                       {title}
    +{.} { }                            {part}
    +{:} { \emph}                       {subtitle}
    +{.} { \PrintContributions}         {contribution}
    +{.} { Technical Report }           {number}
    +{,} { }                            {series}
    +{.} { }                            {organization}
    +{,} { }                            {address}
    +{,} { \PrintDateB}                 {date}
    +{.} { \PrintTranslation}           {translation}
    +{.} { Reprinted in \PrintReprint}  {reprint}
    +{.} { }                            {note}
    +{.} {}                             {transition}
}

\BibSpec{thesis}{%
    +{}  {\PrintAuthors}                {author}
    +{,} { \emph}                       {title}
    +{:} { \emph}                       {subtitle}
    +{.} { \PrintThesisType}            {type}
    +{.} { }                            {organization}
    +{,} { }                            {address}
    +{,} { \PrintDateB}                 {date}
    +{.} { \PrintTranslation}           {translation}
    +{.} { Reprinted in \PrintReprint}  {reprint}
    +{.} { }                            {note}
    +{.} {}                             {transition}
}

\RequirePackage{yhmath}
\DeclareSymbolFont{extrasym}  {U}{zeuex}{m}{n}
\DeclareMathSymbol{\infty}{\mathord}{extrasym}{153}

\usepackage[all]{xy}
\numberwithin{equation}{section}

\providecommand{\half}{{\textstyle\frac{1}{2}}}

\providecommand{\comment}[1]{}
\providecommand{\Space}[3][]{\ensuremath{\mathbb{#2}^{#3}_{#1}{}}}
\providecommand{\FSpace}[3][]{\ensuremath{\ifx#2l \ell_{#3}^{#1}{}\else
    \mathcal{#2}_{#3}^{#1}{}\fi}} 
\providecommand{\rmi}{\mathrm{i}}
\providecommand{\rme}{\mathrm{e}}
\providecommand{\rmd}{\mathrm{d}}
\providecommand{\oper}[1]{\mathsf{#1}}
\providecommand{\uir}[3][0]{\ifcase #1{\rho^{#2}_{#3}}%
\or {\breve{\rho}^{#2}_{#3}}%
\or {\tilde{\rho}^{#2}_{#3}}%
\or {\hat{\rho}^{#2}_{#3}}\fi}
\providecommand{\covar}[2][\relax]{\oper{W}^{#1}_{#2}}
\providecommand{\contravar}[2][\relax]{\oper{M}^{#1}_{#2}}
\providecommand{\algebra}[1]{\ensuremath{\mathfrak{#1}}}

\providecommand{\norm}[2][\relax]{\left\|#2\right\|\ifx#1\relax\else_{#1}\fi}
\providecommand{\modulus}[2][\relax]{\left| #2 \right|\ifx#1\relax\else_{#1}\fi}
\providecommand{\map}[1]{\mathsf{#1}}
\providecommand{\such}{\,\mid\,}
\providecommand{\myhbar}{\hbar}
\providecommand{\myh}{h}
\providecommand{\MR}[1]{MR\href{http://www.ams.org/mathscinet-getitem?mr=#1}{#1}}
\providecommand{\Zbl}[1]{Zbl\href{http://www.emis.de:80/cgi-bin/zmen/ZMATH/en/zmathf.html?first=1&maxdocs=3&type=html&an=#1&format=complete}{#1}}
\providecommand{\scalar}[3][\relax]{\left\langle #2,#3 
        \right\rangle\ifx#1\relax\else_{#1}\fi}
\providecommand{\myeprint}[2]{E-print: \href{#1}{\texttt{#2}}}
\providecommand{\twist}{\mathbin{\mathop{\natural}}}
\providecommand{\Btwist}{\mathbin{\mathop{\tilde{\natural}}}}

\DeclareDocumentCommand\ladder{o m o}{%
  \IfNoValueTF{#3}%
  {L_{#1}^{\!#2}}%
  { \overset{#3}{L}{}_{#1}^{\!#2}}%
}
\DeclareMathOperator{\sign}{sign}

\newcounter{order}

\usepackage{mathtools}

\newcommand\increase[1]{%
  \stepcounter{order}
  \count0=1
  \count1=#1
  \advance\count0 by \count1
  \the\count0
}

\let\pyginac=\iffalse
\let\endpyginac=\fi

\NewDocumentCommand \D { > { \SplitList { , } } o }
{\setcounter{order}{0} \partial_{\ProcessList {#1} { \increase }}
\ifthenelse{\equal{1}{\arabic{order}}}{}{^{\arabic{order}}}}

\DeclareFontEncoding{LS1}{}{}
\DeclareFontSubstitution{LS1}{stix2}{m}{n}
\DeclareSymbolFont{symbols2}{LS1}{stix2frak} {m} {n}
\DeclareMathSymbol{\operp}{\mathbin}{symbols2}{"A8}

\title[Cross-Toeplitz Operators and Two-Sided Convolutions]
{Cross-Toeplitz Operators\\ on the Fock--Segal--Bargmann Spaces\\
  and Two-Sided   Convolutions
  on the Heisenberg Group}

\author[Vladimir V. Kisil]%
{\href{http://www.maths.leeds.ac.uk/~kisilv/}{Vladimir V. Kisil}}
\thanks{On  leave from Odessa University.}

\address{%
School of Mathematics\\
University of Leeds\\
Leeds LS2\,9JT\\
UK
}

\email{\href{mailto:kisilv@maths.leeds.ac.uk}{kisilv@maths.leeds.ac.uk}}
\email{\href{mailto:V.Kisiv@leeds.ac.uk}{V.Kisil@leeds.ac.uk}}

\urladdr{\url{http://www.maths.leeds.ac.uk/~kisilv/}}

\date{\today}
\begin{document}
\begin{abstract}
  We introduce an extended class of cross-Toeplitz operators which act between Fock--Segal--Bargmann spaces with different weights. It is natural to consider these operators in the framework of representation theory of the Heisenberg group. Our main technique is  representation of cross-Toeplitz by two-sided relative convolutions from the Heisenberg group. In turn, two-sided convolutions are reduced to usual (one-sided) convolutions on the Heisenberg group of the doubled dimensionality. This allows us to utilise the powerful group-representation technique of coherent states, co- and contra-variant transforms, twisted convolutions, symplectic Fourier transform, etc. We discuss connections of (cross-)Toeplitz operators with pseudo-differential operators, localisation operators in time-frequency analysis, and characterisation of kernels in terms of ladder operators. The paper is written in a detailed and reasonably self-contained manner to be suitable as an introduction into group-theoretical methods in phase space and time-frequency operator theory.
\end{abstract}
\maketitle

\tableofcontents

\section{Introduction}
\label{sec:introduction}

Motivated by applications~\cite{Kisil21c,Kisil22a} this paper starts a systematic treatment of mixed coherent state transforms (aka time-frequency analysis) with various Gaussian windows. In particular, we introduce cross-Toeplitz operators acting between two Fock--Segal--Bargmann spaces with different Gaussian weights, which are studied through two-sided relative convolutions on the phase space. The paper contains a comprehensive description of the theory out of the main objects of the representation theory. The next subsection provides a summary for readers familiar with the field. A sloping introduction in a wider context starts from \S~\ref{sec:uncert-relat-sque}.

\subsection{Brief summary for specialists}
\label{sec:summ-spec-fsb}

Let \(\FSpace{F}{2}\) be the Fock--Segal--Bargmann (FSB) space~\cites{Bargmann61,Segal60,Zhu12a,Folland89} on the phase space, which can be defined as a certain irreducible component of the unitary representation of the Heisenberg group in \(\FSpace{L}{2}\) space~\cites{Coburn99,Coburn19a}. Let \(\oper{P}: \FSpace{L}{2} \rightarrow \FSpace{F}{2}\) be  the respective FSB orthoprojection. For a function \(\psi\) on the phase space one defines the Toeplitz operator \(\oper{T}_\psi: \FSpace{F}{2} \rightarrow \FSpace{F}{2}\) by the identity
\begin{equation}
  \label{eq:Toeplitz-initial-defn}
  \oper{T}_\psi f=\oper{P} (\psi f), \qquad  \text{ for } f\in \FSpace{F}{2}.
\end{equation}
A study of Toeplitz operators~\cites{Zhu12a,Coburn19a,CoburnIsralowitzLi11a,BauerCoburnIsralowitz10a}, their connections to PDOs~\cite{Guillemin84} and a parallel theory of localisation operators on the phase space~\cites{Grochenig01a,BoggiattoCorderoGrochenig04a,Coburn01a,AbreuFaustino15a} is a vivid research area.

Coburn discovered~\cite{Coburn01b} some fundamental limitations for calculus of Toeplitz operators. For example, a composition operator \([C_\theta  f](z) = f(e^{\rmi \theta } z)\) with \(\pi/4< \theta < \pi/2\) can be represented as a Toeplitz operator \(\oper{T}_\phi\) for a Gaussian type function \(\phi\). However, its square \(C_\theta^2\) can\emph{not} be given in the form~\eqref{eq:Toeplitz-initial-defn} for any \(\psi\).

This paper includes Toeplitz operators into a larger class of \emph{two-sided} relative convolutions from the Heisenberg group~\cites{Kisil94e,Kisil13b,Kisil13a}. Convolutions with certain families of kernels form algebras closed under composition. They are also naturally linked to pseudo-differential operators (PDO) as will be shown in details. Yet,  the Coburn counterexample \(C_\theta\) from the previous paragraph is not a PDO with a symbol in a classical H\"ormander class.  Thus, a proper treatment of \(C_\theta\) requires more general Fourier integral operators~\cite{Hormander83IV}, which is beyond the scope of the present paper.

There are clear advantages of mixed coherent states decompositions over bounded domains~\citelist{\cite{Kisil21c}\cite{Kisil22a}\cite{MazyaSchmidt07a}*{Ch.~9}}.  To this end we introduce and study in this work a new class of cross-Toeplitz operators. They act between FSB spaces with different Gaussian weights~\cite{Zhu12a} or different Gaussian windows in terms time-frequency analysis~\cite{Grochenig01a}. Justification of cross-Toeplitz operators from a wider context is discussed in subsequent subsections of this Introduction.  To work with cross-Toeplitz operator we need numerous adjustments to the framework, starting from the definition of pre-FSB spaces. We systematically present the theory here from the group-representation perspective, which is not yet dominating in the field (with the some notable exceptions, e.g.~\cite{BergerCoburn86b,Coburn99,Coburn19a,Folland89}).

While Toeplitz operators admit representation as usual (one-sided) convolutions on the Heisenberg group, cross-Toeplitz operators need two-sided convolutions because they mix different irreducible components under the Heisenberg group action. Yet, our main tool is a transformation of two-sided convolutions to one-sided convolutions on the Heisenberg group of the doubled dimensionality. This is an example of the Heisenberg group universality for relative convolutions~\cite{Kisil94e}*{Thm~3.8}.  A more costly alternative is a consideration of the nilpotent \emph{step 3} Dynin~\cite{Kisil13a} (aka Dynin--Folland~\cite{RottensteinerRuzhansky20a}) group build on top of the simplest nilpotent step 2 Heisenberg group. Another alternative is a solvable group, which is a semi-direct product of the Heisenberg group and its one-parameter group of automorphisms by inhomogeneous dilations.

\subsubsection{The paper outline}
\label{sec:paper-outline}
 In Section~\ref{sec:heis-group-essent} we recall the minimal background on the Heisenberg group and its representation theory. The standard construction of the induced representations puts all FSB spaces with different Gaussian weights (or rather their unitary equivalent counterparts) as irreducible subspaces of \(\FSpace{L}{2}(\Space{R}{2n})\). In the traditional approach from holomorphic function perspective all these FSB spaces are detained in differently weighted \(\FSpace{L}{2}\)-spaces, which prevents a natural consideration of cross-Toeplitz operators so far. The key ingredients for two-sided convolutions---the left and right pulled action---are introduced here as well.

In Section~\ref{sec:covar-amd-contr} we collect fundamental facts on the co- and contra-variant transforms also known under numerous other names in various fields, e.g. coherent states, voice, Berezin, Wigner, FSB, etc. transforms. Together with twisted convolution and symplectic Fourier transform these will be our main tools. We review complexification of the theory and appearance of (poly-)analytic functions from the representation theory in Section~\ref{sec:compl-vari-anal}. Ladder operators for the left and right actions play different but equally important r\^oles here.

We turn to our study of (cross-)Toeplitz operators in Section~\ref{sec:repr-oper-pdo}. To larger extend results of this section go back to original works of Howe~\cite{Howe80a,Howe80b} and Guillemin~\cite{Guillemin84}. Yet we are able to add some more group-theoretic perspectives to this material as well. We discuss connections to the theory of PDO, localisation operators in time-frequency analysis and a significance of the heat flow on symbols.

The final Section~\ref{sec:two-sided-conv-toeplltz} contains most of the original material. Here we represent  cross-Toeplitz operators as two-sided relative convolutions from the Heisenberg group. Two-sided convolutions are reduced to one-sided convolutions from the Heisenberg group of double dimensionality. Thus, we can re-cycle all the theory presented in earlier sections for new needs. Cross-Toeplitz operators are treated through the symbolic calculus of PDO and their symbols  are characterised through certain identities with ladder operators. A point of interest is that a PDO symbol of \(\oper{T}_\psi\) is expressed through the invertible FSB-like transform of \(\psi\) instead of smoothing heat flow considered in the literature so far.

Now we turn to a discussion of cross-Toeplitz operators within a wider context of harmonic analysis. 

\subsection{Uncertainity relation and squeeze transformation}
\label{sec:uncert-relat-sque}

The fundamental idea of the Fourier analysis is that a function \(f(t)\) of a real variable can be represented as a superposition of simple harmonics \(\rme^{2\pi \rmi k t}\) with all real wave numbers \(k\). The original application to differential equations---due to the fact that harmonics are eigenfunctions of the derivative---was joined later by numerous others.

There is a fundamental restriction to combine the Fourier analysis with another common technique---localisation: a smaller support of a function \(f(t)\) implies a wider bandwidth of representing harmonics. The phenomenon was manifested as the Heisenberg uncertainty relation in quantum mechanics, which assigns certain de Broglie wave to every elementary object. It is equally relevant for time-frequency analysis of voice/signal functions. Since these  theories are mathematically equivalent we will interchangeably use terms from both areas. 

For quantitative description of the uncertainty we need a coefficient to make coordinates and wave numbers to be co-measurable. In quantum mechanics it is known as the Planck constant \(\myhbar\). 
Then, for the coordinate \({Q}\) and momentum \(P\) operators there is the
{Heisenberg--Kennard uncertainty relation}~\cite{Folland89}*{\S~1.3}:
\begin{equation}
  \label{eq:heisenberg-uncertainty}
  \Delta_\phi(Q) \cdot \Delta_\phi(P) \geq \frac{\myhbar}{2},
\end{equation}
where \(\Delta_\phi(A)\) is the dispersion of an operator \(A\) on the state \(\phi\).
The equality holds if and only if \(\phi\) is the Schr\"odinger coherent state:
\begin{equation}
  \label{eq:schrodinger-coherent-state}
  \phi_{qp}(t)=\sqrt[4]{\tau}\rme^{2\pi\rmi \myhbar p t} \rme^{-\pi \myhbar (t-q)^2/(2 \tau)}.
\end{equation}
which are \((q,p)\)-translations in the phase space of the squeezed Gaussian \(\phi_\tau(t)=\sqrt[4]{\tau}\rme^{-\pi  \myhbar t^2/(2 \tau)}\), see~\eqref{eq:schroedinger-rep} below. Fig.~\ref{fig:two-sque-gauss} shows two differently squeezed Gaussians on the configuration space (on the left graph) and the respective regions of their essential support in the phase space (on the right graph). Green-solid and blue-dashed pencils are used for the respective drawings on both graphs.
\begin{figure}
  \centering
  \includegraphics[scale=.8]{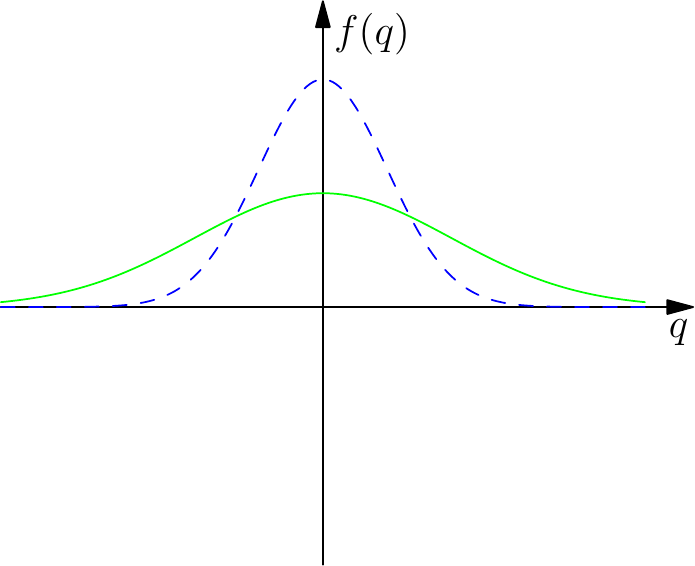}\hfill
  \includegraphics[scale=.8]{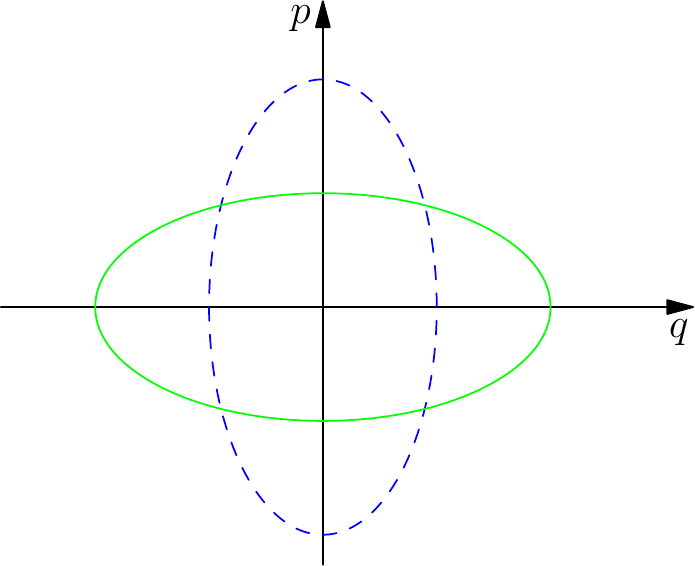}
  \caption[Two squeezed Gaussians]{Two differently squeezed Gaussians on the configurational space (left) and their phase space footprints (right). Green-solid and blue-dashed pencils are used for the respective drawings on both graphs.}
  \label{fig:two-sque-gauss}
\end{figure}

\begin{figure}
  \centering
  \includegraphics{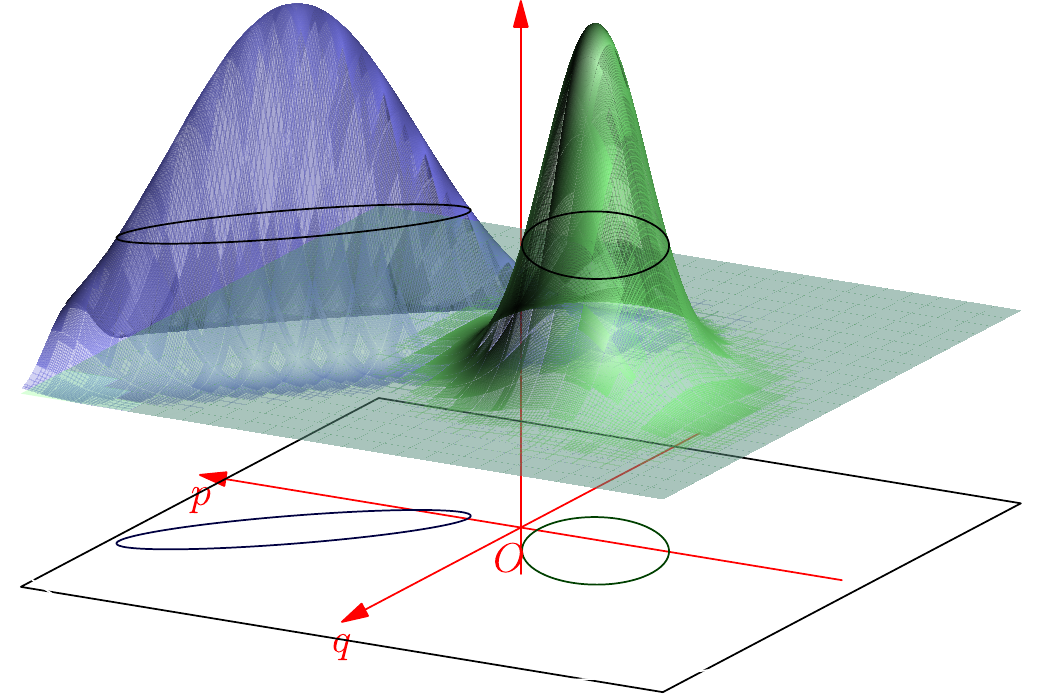}
  \caption{Squeezed states on the phasespace and their blobs}
  \label{fig:sque-stat-phas}
\end{figure}

Clearly, squeezed Gaussians are fully determined by these types of ellipses---their footprints on the phase space, see Fig.~\ref{fig:sque-stat-phas}. This concept (branded as quantum blobs~\cites{deGosson13a,deGosson11a}) nicely visualises the phase space uncertainty.
The squeeze parameter \(\tau\) allows us to trade a better localisation in the  configuration space for a more defused bandwidth. Thus a possibility to interchangeably operate with differently squeezed Gaussians is an advantage.

\subsection{Fock--Segal--Bargmann spaces and Toeplitz operators}
\label{sec:fock-segal-bargmann}

The Fourier transform through simple harmonics can be seen as the limiting case of a decomposition into a superposition of the Schr\"odinger coherent states~\eqref{eq:schrodinger-coherent-state} for \(\tau \rightarrow 0\). It achieves a Dirac delta-like localisation in the frequency scale on the expense of no localisation whatsoever in coordinates. For a fixed \(\tau \neq 0\) such a decomposition  is known as the Fock--Segal--Bargmann (FSB) transform~\citelist{\cite{Fock04a}*{\S~28-3} \cite{Segal60} \cite{Bargmann61} } with  a quantum mechanical origin. Alternatively it appears in time-frequency analysis as the Gabor transform of a signal into simple harmonics \(\rme^{2\pi\rmi \myhbar p t}\) modulated by the Gaussian window \(\sqrt{ \tau}\rme^{-\pi  \tau(t-q)^2/2}\)~\cite{Grochenig01a}. FSB transform is an archetypal source of numerous developments: coherent states~\cites{Perelomov86,AliAntGaz14a}, atomic decompositions~\cites{FeichGroech89b,FeichGroech89a}, etc.

It is common to adjust the construction in such a way that the image space of the FSB transform consists of holomorphic functions. Specifically, define Gaussian-weighted measure on \(\Space{C}{n}\):
\begin{displaymath}
  \rmd \mu_{\tau}(z) = (4 \tau)^{-n} \exp(-\pi \modulus{z}^2/(4\tau))\,\rmd z \wedge \rmd \overline{z}, \quad \text{ where } \tau>0.
\end{displaymath}
The FSB space \(\FSpace[\tau]{F}{}\) is commonly defined as the closed subspace of \(\FSpace{L}{2}(\Space{C}{n}, \rmd \mu_{\tau})\) consisting of holomorphic functions. There is the associated orthoprojection 
\(\oper{P}_\tau: \FSpace{L}{2}(\Space{C}{n}, \rmd \mu_{\tau}) \rightarrow \FSpace[\tau]{F}{}\). Then any bounded function \(\psi\) on \(\Space{C}{n}\) defines the Toeplitz operator \(\oper{T}_\psi: \FSpace[\tau]{F}{} \rightarrow \FSpace[\tau]{F}{}\) by the rule \(\oper{T}_\psi f= \oper{P}_\tau(\psi f)\), for \(f\in \FSpace[\tau]{F}{}\).

The correspondence of functions \(\psi\) on the phase space (classical observables) to Toeplitz operators \(\oper{T}_\psi \) (quantum observables) is known as  Berezin--Toeplitz quantization---a physically significant~\citelist{ \cite{Berezin86} \cite{Folland89}} and mathematically rich~\citelist{ \cite{Zhu12a} \cite{Coburn19a}} concept. 
Another mainstream quantisation employs the Fourier transform in a different way: a pseudo-differential operator (PDO) \(a(X,D)\) with a Weyl symbol \(a(x,\lambda)\)  is defined by:
\begin{displaymath}
  [a(X,D)f](x) =\int\int a(\half(x+y),\lambda )\, f(y)\, \rme^{2\pi \rmi \lambda (x-y) } \,\rmd \lambda \,\rmd y.
\end{displaymath}
Relations between Berezin--Toeplitz and Weyl quantisation was already addressed by Guillemin~\cite{Guillemin84} and continues to be in the focus of current research~\cite{Coburn19a}.  Effectively, the Toeplitz operator \(\oper{T}_\psi\) is unitary equivalent to PDO \(\psi_\tau(X,D)\) for
\begin{equation}
  \label{eq:smoothin-PDO-symbol}
  \psi_\tau=\psi*\Phi_\tau, 
\end{equation}
which is a smoothing of \(\psi\) by the convolution with the Gaussian \(\Phi_\tau(q,p) = \rme^{-\pi \myhbar( q^2/\tau+ \tau p^2)/2}\)~\cite{Folland89}*{(3.5)}. It can be also interpreted as a heat/diffusion flow transformation~\cite{BergerCoburn94a}. Of course, the degree of smoothing depends on \(\tau\) which again put this parameter in a focus of our attention. 

\subsection{Variable squeezing}
\label{sec:variable-squeezing}

We have seen that the Planck constant \(\myhbar\) limits the joint localisability in the phase space~\eqref{eq:heisenberg-uncertainty} while the squeeze parameter \(\tau\) controls shares of uncertainty between coordinates and momentum (or time and frequency). While the fundamental physical constants, e.g. the Planck constant, cannot be affected, we can and sometime want to change parameters like squeezing. Here are few illustrations:
\begin{enumerate}
\item An upgraded version of FSB transform which treats \(\tau\) as a parameter on a par with \(p\) and \(q\) is known as FBI  (Fourier--Bros--Iagolnitzer) transform. It is used, for example, to analyse wave fronts~\cite{Folland89}*{\S~3.3}.
\item An adaptation of FSB transform with \(\tau\) being a power function of \(p\) was employed~\citelist{\cite{CordobaFefferman78a} \cite{Folland89}*{(3.6)}} for a better approximation of a PDO by a Toeplitz operator. Indeed the smoothing~\eqref{eq:smoothin-PDO-symbol} harder affects high frequencies. To reduce this effect we need a more narrow Gaussian window with its width to be much smaller than targeted wavelengths.
\item There is a possibility to geometrise a quantum dynamics by extending the classical phase space with additional coordinates and the squeeze parameter \(\tau\) is a suitable option. Examples of coherent states with oscillating squeeze were described in~\cites{AlmalkiKisil18a,AlmalkiKisil19a}.
\item\label{item:archimedes} Adopting coherent states decomposition (aka time-frequency analysis) on bounded domains we can use Gaussian windows of narrowing widths through the Archimedes' method of exhaustion~\citelist{\cite{Kisil22a}\cite{MazyaSchmidt07a}*{Ch.~9}}.
\end{enumerate}

The present paper widens the study of Toeplitz operators between FSB spaces with different squeeze parameters. Recall, that for each value \(\tau\) the FSB space \(\FSpace[\tau]{F}{}\) is unitary equivalent (through FSB transform) to the space \(\FSpace{L}{2}(\Space{R}{n})\) of functions on the configurational space. Thus, for two squeeze parameters \(\tau\) and \(\varsigma\) there is a natural unitary equivalence of \(\FSpace[\tau]{F}{}\) and \(\FSpace[\varsigma]{F}{}\). That is, this map recalculates the FSB transform with a squeeze parameter \(\tau\) into the FSB  transform of the same function in \(\FSpace{L}{2}(\Space{R}{n})\) for a different squeeze parameter \(\varsigma\).
Alternatively, the equivalence  of \(\FSpace[\tau]{F}{}\) and \(\FSpace[\varsigma]{F}{}\) corresponds to the unitary dilation \(f(t) \mapsto r^{n/2}f(r t)\) with \(r=\sqrt{\tau/\varsigma}\) on \(\FSpace{L}{2}(\Space{R}{n})\). 

Furthermore, for a function \(\psi\) on the phase space we consider cross-Toeplitz operator \(\oper{T}_\psi: \FSpace[\tau]{F}{} \rightarrow \FSpace[\varsigma]{F}{}\) by the identity:
\begin{displaymath}
  \oper{T}_\psi f= \oper{P}_\varsigma (\psi f), \qquad \text{ where } f \in \FSpace[\tau]{F}{}.
\end{displaymath}
We develop some basic result on the calculus of cross-Toeplitz operators using group representation technique notably two-sided relative convolutions from the Heisenberg group~\citelist{\cite{VasTru88} \cite{VasTru94}
  \cite{Vasilevski94a} \cite{Kisil92} \cite{Kisil93e} \cite{Kisil93b}
  \cite{Kisil94a} \cite{Kisil94f} \cite{Kisil96e}}. Interestingly, some obtained results are new even for the case \(\tau=\varsigma\) of the traditional Toeplitz operators. One of applications of cross-Toeplitz operators  is related to the above item~\ref{item:archimedes}.  Zones with different Gaussian widths are overlapping  and cross-Toeplitz operators are a right tool to account the imbrication~\cite{Kisil22a}. 

\subsection{Coherent states and group representations}
\label{sec:coher-stat-group}

Our choice of group representations technique is well-motivated and historically based. 
Connections of the Heisenberg group to quantum mechanics are rooted in the Stone-von Neumann theorem.  It established the coincidence of formalisms based on the Heisenberg matrix mechanics and the Schr\"odinger equation through the classification of the unitary irreducible representations of the Heisenberg group. Various values of the Planck constants parameterise non-equivalent classes of representations.  Furthermore, the Schr\"odinger coherent states~\eqref{eq:schrodinger-coherent-state} are obtained from the vacuum vector---the squeezed Gaussian \(\phi_\tau(t)=\sqrt[4]{\tau} \rme^{-\pi  \myhbar t^2/(2 \tau )}\)---by the irreducible Schr\"odinger representation of the Heisenberg group on \(\FSpace{L}{2}(\Space{R}{})\). Correspondingly, the  FSB transform is the prime example of the coherent states transform generated by a square-integrable unitary group representation~\cites{AliAntGaz14a,Perelomov86,Folland89}.

On the other hand, a link of  the calculus of PDO to the Heisenberg group was indicated
as early as~\cite{GrossmannLoupiasStein68}. It was spectacularly
developed in~\cites{Howe80b,Howe80a} and spelt in details
in~\cite{Folland89}. Toeplitz operators on the Fock space from the Heisenberg group appeared
in~\cites{Howe80b,Howe80a} implicitly (see
Sect.~\ref{sec:toepl-oper-cald} below), their connection to PDO (and
therefore to the Heisenberg group) was explicitly revealed by
Guillemin~\cite{Guillemin84}. Folland elaborated this in his
step-stone monograph~\cite{Folland89}. Further usage of the Heisenberg group for 
Toeplitz operators can be found in~\cites{Coburn01b,Coburn99,Coburn01a,Coburn12a,Coburn19a}.

To enable various values of the squeeze parameter we extend the Heisenberg group using the symplectic transformation \((q,p)\mapsto(rq,p/r)\) of the phase space, see the right graph on Fig.~\ref{fig:two-sque-gauss} for a visualisation.  The Heisenberg group acts on the phase-space by shifts in a reducible way, in particular every FSB space \(\FSpace[\tau]{F}{}\) is invariant (and irreducible). Thereafter, the von Neumann algebra of operators generated by the phase space shifts contains Toeplitz operators (with suitable class of symbols). However, it can not contain any cross-Toeplitz operators \(\FSpace[\tau]{F}{} \rightarrow \FSpace[\varsigma]{F}{}\) for \(\tau\neq \varsigma\).  On the other hand, the extension of the phase shift with squeeze transformations acts on \(\FSpace{L}{2}(\Space{C}{n})\) irreducibly and cross-Toeplitz operators become covered.

From general principles, we can use  other extensions of phase space shifts to achieve irreducibility on \(\FSpace{L}{2}(\Space{C}{n})\), see the mentioned in \S~\ref{sec:summ-spec-fsb} Dynin group as an example.  Our choice allows to achieve results virtually at no extra cost (e.g. without increments of nilpotency of the group), we combine the initial (left) phase space shifts by their right counterparts. This allows us to study cross-Toeplitz operators  as two-sided relative convolutions from the Heisenberg group. The latter can be reduced to usual (one-sided) convolutions on a Heisenberg group of the doubled dimensionality~\citelist{\cite{Kisil92} \cite{Kisil93e} \cite{Kisil93b}  \cite{Kisil94a} \cite{Kisil94f} \cite{Kisil96e}}. This technique successfully treats cross-Toeplitz operators and also gives some new perspective on the traditional Toeplitz operators.




 \section{The Heisenberg group essentials}
\label{sec:heis-group-essent}
This section contains a minimal presentation of the Heisenberg
group required for this paper. We refer to other
works~\cites{Folland89,Kisil17a,Kirillov04a,Alamer19a} for further details.

An element of the \(n\)-dimensional \emph{Heisenberg group}
\(\Space{H}{n}\)%
\index{Heisenberg!group|indef}%
\index{group!Heisenberg|indef}%
\index{$\Space{H}{n}$ (Heisenberg group)} is
\((s,x,y)\in\Space{R}{2n+1}\), where \(s\in\Space{R}{}\) and \(x\),
\(y\in \Space{R}{n}\).  The group law on
\(\Space{H}{n}\) is given as follows:
\begin{equation}
  \label{eq:H-n-group-law}
    (s,x,y)\cdot(s',x',y')= \textstyle (s+s'+\frac{1}{2}\omega(x,y;x',y'),x+x',y+y'),\\
\end{equation} 
where \(\omega\) is the symplectic form:
\begin{align}
  \label{eq:symplectic-form}
  \omega(x,y;x',y')&=xy'-x'y
.
\end{align}
It is often convenient to identify
\((x,y)\in\Space{R}{2n}\) with \(z=x+\rmi y\in\Space{C}{n}\) with the
group law presented by:
\begin{equation}
  \label{eq:H-n-group-law-complex-early}
    (s,z)\cdot(s',z')= \textstyle (s+s'+\frac{1}{2}\Im (z\overline{z}'),z+z'),
\end{equation} 
with the symplectic form \( \omega(z;z') =\Im (z\overline{z}')\) in the
complex coordinates \(z=x+\rmi y\) and
\(z'=x'+\rmi y'\in\Space{C}{n}\). However, some additional
considerations are required for this, see
\S~\ref{sec:ladd-oper-compl}.

The two-sided Haar (the left and the right invariant) measure on the
Heisenberg group coincides with the Lebesgue measure on
\(\Space{R}{2n+1}\).

\subsection{Induced representations}
\label{sec:induc-repr}

The bare minimum of the induced representation construction for the Heisenberg group is as follows, see numerous sources for detailed treatments of the induced representation technique~\cites{Kirillov04a,Kirillov76,Folland16a}.

For a real \(\myhbar\neq 0\) we denote by \(\chi_\myhbar\) the
non-trivial unitary character of the group \((\Space{R}{},+)\):
\begin{equation}
  \label{eq:character-defn}
  \chi_\myhbar(s) = \rme^{2\pi \rmi \myhbar s}\,.
\end{equation}
Consider functions on \(\Space{H}{n}\)
having the following covariant property with the respect to the
central action on \(\Space{H}{n}\):
\begin{equation}
  \label{eq:Z-covariance}
  F(s'+s,x,y)=\overline{\chi}_\myhbar(s)F(s',x,y)=\rme^{-2\pi \rmi \myhbar s}F(s',x,y) \qquad \text{ for all
  }\ s\in\Space{R}{}.
\end{equation}
This characteristic is preserved under the left \(\Lambda(g)\) and right
shifts \(R(g)\):
\begin{equation}
  \label{eq:left-right-shift-group}
  \Lambda(g): F(g') \mapsto F(g^{-1}g'), \qquad
  R(g): F(g') \mapsto F(g'g),\quad g,g'\in \Space{H}{n}.
\end{equation}
Clearly, a function \(F(s,x,y)\) satisfying~\eqref{eq:Z-covariance} is
completely defined by its values \(F(0,x,y)\) on \(\Space{R}{2n}\subset
\Space{H}{n}\). Thus, for functions \(F\) and \(f\) on
\(\Space{H}{n}\) and \(\Space{R}{2n}\) respectively we define the
\emph{lifting} \(\oper{L}_\myhbar\) and \emph{pulling} \(\oper{P}\) as
follows:
\begin{equation}
  \label{eq:lifting-pulling}
  [\oper{L}_{\myhbar} f] (s,x,y)=\rme^{-2\pi \rmi \myhbar s}f(x,y),\qquad
  [\oper{P} F](x,y)=F(0,x,y)
  .
\end{equation}
Obviously, the pulling is a left inverse of the lifting:
\([(\oper{P}\circ \oper{L}_\myhbar)  f](x,y)=f(x,y)\). Also
\([(\oper{L}_\myhbar\circ \oper{P}) F](s,x,y)=F(s,x,y)\) if \(F\)
satisfies~\eqref{eq:Z-covariance}. Thus, we can pull the regular
actions~\eqref{eq:left-right-shift-group} on \(\Space{R}{2n}\):
\begin{equation}
  \label{eq:left-right-action-pulled-def}
  \Lambda_\myhbar(g)\coloneqq \oper{P}\circ \Lambda(g) \circ
   \oper{L}_\myhbar, \qquad
  R_\myhbar(g)\coloneqq \oper{P}\circ R(g) \circ
  \oper{L}_\myhbar,
\end{equation}
or, explicitly:
\begin{align}
  \label{eq:left-action-pulled}
  {}[\Lambda_\myhbar(s,x,y)f] (x',y') &= \rme^{\pi \rmi \myhbar
                                      (2s+xy'-yx')} f(x'-x,y'-y), \\
  \label{eq:right-action-pulled}
  [R_\myhbar(s,x,y)f] (x',y') &= \rme^{\pi \rmi \myhbar
                                (-2s+xy'-yx')} f(x'+x,y'+y)
                                       .
\end{align}
These are infinite-dimensional representations of \(\Space{H}{n}\)
which are induced in the sense of Mackey~\cite{Kirillov04a}*{App.~V.2}
from the character
\(\chi_\myhbar(s,0,0) = \rme^{2\pi \rmi \myhbar
  s}\)~\eqref{eq:character-defn} of the centre of \(\Space{H}{n}\). To
distinguish \(\Lambda_\myhbar\) and \(R_\myhbar\)
\eqref{eq:left-action-pulled}--\eqref{eq:right-action-pulled} from
\(\Lambda\) and \(R\)~\eqref{eq:left-right-shift-group} we call formers
as \emph{pulled actions}.


\begin{rem}[General scheme of induced representation]
  \label{re:induced-representation}
  For a later use in \S~\ref{sec:pdo-representing-an} we need
  a bit more general setup of induced
  representations~\citelist{\cite{Kirillov76}*{\S~13.2}
    \cite{MTaylor86}*{Ch.~5} \cite{Folland16a}*{Ch.~6}
    \cite{Kirillov04a}*{\S~V.2} \cite{Mackey70a}}.  Let \(G\) be a
  locally compact group and \(*\) denote its multiplication, \(H\) be
  its closed subgroup and \(X=G/H\)--the respective homogeneous
  space. Let \(\map{s}: X \rightarrow G\) be a section of the bundle defined by the
  natural projection \(\map{p}: G \rightarrow X\). Then, any
  \(g\in G\) has a unique decomposition of the form
  \(g= \map{s}(x)* h\) where \(x=\map{p}(g)\in X\) and \(h\in H\). We
  also define the map \(\map{r}: G\rightarrow H\):%
  \index{map!$r$}\index{$r$ (map)}
\begin{equation}
  \label{eq:r-map1}
  \map{r}(g)={\map{s}(x)}^{-1} * g, \qquad \text{ where } x=\map{p}(g).
\end{equation}
This map is our substitution~\citelist{\cite{Kisil97c}*{\S~3.1}
  \cite{Kisil09e} \cite{Kisil10a} \cite{Kisil13a}} for the so-called
master equation~\cite{Kirillov04a}*{\S~V.2}. Then, the left action of
\(G: X\rightarrow X\) is given by
\begin{equation}
  \label{eq:group-left-action-hs}
  g: x \mapsto g\cdot x=\map{p}(g*\map{s}(x)), \qquad \text{ where }
  g\in G,\ x\in X.
\end{equation}
A character \(\chi\) of \(H\) induces a representation
\(\uir{}{\chi}\) of \(G\) in a
space of function on \(X\) by the formula:
\begin{equation}
  \label{eq:induced-representation-def}
  \uir{}{\chi}(g): f(x) \rightarrow \overline{\chi}(\map{r}(g^{-1}*\map{s}(x)))\,f(g^{-1}\cdot x),
\end{equation}
The complex conjugation of the character \(\chi\) is used here to make
our formula in line the general induced representation construction.
If \(X\) posses a \(G\)-invariant measure \(\rmd x\) then this
representation acts by isometries of \(\FSpace{L}{p}(X,\rmd
x)\). There is a suitable adaptation for quasi-invariant measures as well. 

\end{rem}



The
representations~\eqref{eq:left-action-pulled}--\eqref{eq:right-action-pulled}
of \(\Space{H}{n}\) are reducible. An irreducible representation of
\(\Space{H}{n}\) can be induced from the character
\(\chi_\myhbar(s,0,y)=\rme^{2\pi \rmi \myhbar s}\) of the continuous
two-dimensional maximal commutative subgroup
\begin{displaymath}
  H_x=\{(s,0,y)\in \Space{H}{n}:\ s\in \Space{R}{}, y\in
\Space{R}{n}\}
\end{displaymath}
 Using the above steps with the respective lifting
and pulling (or maps \(\map{p}\), \(\map{s}\), \(\map{r}\) from
Rem.~\ref{re:induced-representation}), the \emph{Schr\"odinger
  representation} \(\uir{}{\myhbar}\) of \(\Space{H}{n}\) on
\(\FSpace{L}{2}(\Space{R}{n})\) can be written as:
\begin{equation}
  \label{eq:schroedinger-rep}
  \textstyle
  [\uir{}{\myhbar}(s,x,y)f](t)= \rme^{\pi \rmi \myhbar (2
  s-2yt+xy)}\, f(t-x).
\end{equation}

By the Stone-von Neumann theorem, any irreducible infinite-dimensional
representation of \(\Space{H}{n}\)---in particular, an irreducible
subrepresentation
of~\eqref{eq:left-action-pulled}--\eqref{eq:right-action-pulled}---is
equivalent to~\eqref{eq:schroedinger-rep} with the same value of
\(\myhbar\). The intertwining operator between
representations~\eqref{eq:schroedinger-rep}
and~\eqref{eq:left-action-pulled} will be constructed as coherent
states transform in Ex.~\ref{ex:gaussian-fsb}.

In the context of representations of Lie groups the corresponding
derived representations of their Lie algebras are important. Take a basis:
\begin{equation}
  \label{eq:Lie-algebra-basis}
  S=(1,0,0), \quad X=(0,1,0), \quad Y=(0,0,1) 
\end{equation}
of the Lie algebra \(\algebra{h}_1\) of \(\Space{H}{1}\). For the
representations~\eqref{eq:left-action-pulled},
\eqref{eq:right-action-pulled} and \eqref{eq:schroedinger-rep} the derived representation is expressed through the following differential operators:
\begin{align}
  \label{eq:left-action-derived}
  \rmd \Lambda^S_\myhbar&= 2\pi \rmi \myhbar I ,&
  \rmd \Lambda^X_\myhbar&= \pi \rmi \myhbar y -\partial_x,&
  \rmd \Lambda^Y_\myhbar&= -\pi \rmi \myhbar x -\partial_y;\\
  \label{eq:right-action-derived}
  \rmd R^S_\myhbar&=  -2\pi \rmi \myhbar I,&
  \rmd R^X_\myhbar&= \pi \rmi \myhbar y +\partial_x,&
  \rmd R^Y_\myhbar&= -\pi \rmi \myhbar x +\partial_y;\\
  \label{eq:schrodinger-derived}
  \rmd \uir{S}{\myhbar}&= 2\pi \rmi \myhbar I ,&
  \rmd \uir{X}{\myhbar}&= -\partial_t ,&
  \rmd \uir{Y}{\myhbar}&=  -2\pi \rmi \myhbar t.
\end{align}
Each of these derived representations represents the canonical (aka Heisenberg) \emph{commutator relation}
\begin{displaymath}
  [X,Y]=S.
\end{displaymath}
It is the source of the Heisenberg--Kennard uncertainty relation~\eqref{eq:heisenberg-uncertainty}, see~\cite{Folland89}*{\S~1.3} for further details and~\cite{Kisil13c} for some recent developments.

\subsection{Pulled actions and the phase space}
\label{sec:heis-group-phase}

There are several important links be\-twe\-en the left and right actions:
\begin{enumerate}
\item In general, the left and right regular representations for any
  group \(G\) are equivalent and are intertwined by the reflection
  \(g \mapsto g^{-1}\), \(g\in G\). For the Heisenberg group this
  implies that the reflection of the function \emph{domain}
  \(\oper{R}: f(x,y) \mapsto f(-x,-y)\) intertwines the left and the right pulled
  actions with \emph{opposite} Planck constants:
  \begin{equation}
    \label{eq:right-through-left-reflection}
    \Lambda_\myhbar(s,x,y) \circ \oper{R} = \oper{R} \circ R_{-\myhbar}(s,x,y)\,.
  \end{equation}
\item The complex conjugation (that is a reflection in the function
  \emph{range}) intertwines the left and right actions with
  near-inverse elements:
  \begin{equation}
    \label{eq:conjugation-left-right}
    \Lambda_\myhbar(s,-x,-y)\overline{f}= \overline{R_\myhbar(s,x,y)f}\,.
  \end{equation}
\item Finally, a specific feature of nilpotent step two Lie groups and
  \(\Space{H}{n}\) particularly implies:
  \begin{equation}
    \label{eq:right-through-left}
    \Lambda_\myhbar(-s,-x,-y){f} (x',y')= {R_{\myhbar}(s,x',y')f (x,y)}\,.
  \end{equation}
  Note the \emph{swap} of primed and unprimed variables in the last expression
  in comparison to~\eqref{eq:right-through-left-reflection}.  This will
  be used later for representing the FSB projection as a right
  integrated representation in~\eqref{eq:reproducing-projection-double}. Also combination of~\eqref{eq:conjugation-left-right} and~\eqref{eq:right-through-left} yields:
  \begin{equation}
    \label{eq:swap-variables-conjugated}
    \overline{\Lambda_\myhbar(s,x,y)f (x',y')}=     \Lambda_\myhbar(-s,x',y')\overline{\oper{R}f (x,y)}\,.
  \end{equation}
\end{enumerate}
The pulled actions \(\Lambda_\myhbar\) and \(R_\myhbar\)
are also connected to main operators of the classical abelian harmonic analysis---the
Euclidean shift \(\map{S}\) and multiplication \(\map{E}_\myhbar\) by a
unimodular plane wave:
\begin{align}
   \label{eq:shift-defn}
  \map{S} (x,y):\quad f(x',y') &\mapsto f(x'-x,y'-y)\,,
  \\
   \label{eq:exp-mult-defn}
  \map{E}_\myhbar (x,y):\quad f(x',y') &\mapsto \rme^{\pi \rmi \myhbar
                                 (xy'-yx')}f(x',y')\,.
\end{align}
Both \(\map{S}\) and \(\map{E}_\myhbar\) are reducible representations
of the abelian group \((\Space{R}{2n},+)\) by invertible isometries in
\(\FSpace{L}{p}(\Space{R}{2n})\). More accurately, \(\map{E}_\myhbar\)
is a representation of the \emph{dual group} \(\hat{G}\) of
\(G=(\Space{R}{2n},+)\) by (symplectic) \emph{phase shifts}. Of course, \(\hat{G}\)
is isomorphic to \(G\) but there is no a canonical
isomorphism between two groups and the Planck constant \(\myhbar\) labels
different identifications.  Operators \(\map{S}\) and \(\map{E}_\myhbar\)  have the maximal localisation in the configuration and frequency spaces respectively, cf. \S~\ref{sec:uncert-relat-sque}.

The direct product of groups
\(G\times \hat{G}\) is called \emph{phase space}. For \(n=1\) this
space is also known as \emph{time-frequency} space with
\(\map{E}_\myhbar\) called \emph{frequency shift} or
\emph{modulation}~\cite{Grochenig01a}. In this context the Heisenberg
group is the central extension of the phase space which contains different
automorphisms of the phase space interchanging \(G\) and
\(\hat{G}\):
\begin{displaymath}
  (x,y) \mapsto (\myhbar y, -\myhbar x ),\qquad (x,y)\in G,
\end{displaymath}
 parametrised by the non-zero Planck constants.

Since spatial and phase shifts do not commute with each other they are
able to represent the non-commutative pulled actions:
\begin{align}
  \label{eq:left-action-from-shift-mul}
    \Lambda_\myhbar(s,x,y) &= \chi_\myhbar(s)\cdot \map{E}_\myhbar (x,y) \circ \map{S} (x,y)\,,\\
  \label{eq:right-action-from-shift-mul}
  R_\myhbar(s,x,y) &= \overline{\chi}_\myhbar(s)\cdot \map{E}_\myhbar (x,y) \circ \map{S} (-x,-y)\,.
\end{align}
Expressions in
the opposite direction are:
\begin{align}
  \label{eq:shift-by-left-right}
  \map{S} (x,y)&=\Lambda_\myhbar (s,\half x,\half y)\circ R_\myhbar
                 (s,-\half x,-\half y) \\
  \label{eq:shift-by-left-right-exp}
                 & =\exp\left(\half x(\rmd\Lambda^X -\rmd R^X)+\half y(\rmd\Lambda^Y -\rmd R^Y) \right),\\
  \label{eq:mult-by-left-right}
  \map{E}_\myhbar (x,y)&=\Lambda_\myhbar (s,\half x,\half y)\circ R_\myhbar
                         (s,\half x,\half y)\\
  \label{eq:mult-by-left-right-exp}
                         & = \exp\left(\half x(\rmd\Lambda^X +\rmd R^X) + \half  y(\rmd\Lambda^Y +\rmd R^Y) \right)
\end{align}
for an arbitrary \(s\in\Space{R}{}\).

Of course, \(\FSpace{L}{p}\)-norms,  \(1\leq p \leq \infty\), are  invariant
under spatial and phase shift. There are other invariant norms, which are yet less studied, for
example~\cite{Johansson08a,Miheisi10a}:
\begin{equation}
  \label{eq:johansson-norm}
  \norm[J]{f}= \sup_{(x,y)\in\Space{R}{2n}} \left(\int_Q
    \modulus{f(x+x',y+y')}^p\,\rmd x'\,\rmd y'\right)^{\frac{1}{p}},
\end{equation}
where \(Q\) is the unit cube (or, equivalently, any compact set) in
\(\Space{R}{2n}\). The
expressions~\eqref{eq:left-action-from-shift-mul}--\eqref{eq:right-action-from-shift-mul}
of pulled actions through the spatial and phase shifts imply:
\begin{lem}
  \label{le:shift-invariant-norm}
  Let \(V\) be a space of functions on \(\Space{R}{2n}\) with a norm
  which is invariant under the spatial and phase shifts. Then, the pulled
  actions~\eqref{eq:left-action-pulled}--\eqref{eq:right-action-pulled}
  are invertible isometries \(V\).
\end{lem}

\begin{rem}
  Since in both
  cases~\eqref{eq:left-action-pulled}--\eqref{eq:right-action-pulled}
  the centre's action reduces to multiplication by a constant it is
  tempting to use only the ``essential'' actions
  \(\Lambda_\myhbar(0,x,y)\) and \(R_\myhbar(0,x,y)\). However, this
  removes the advantages of the group structure, e.g. a composition
  \(\Lambda_\myhbar(0,x,y)\Lambda_\myhbar(0,x',y')\) cannot be
  presented as \(\Lambda_\myhbar(0,x'',y'')\) without an additional
  factor.
\end{rem}

\begin{rem}[Spaces of analytical functions]
  \label{re:peeling-map}
  The raw representation~\eqref{eq:left-action-pulled} of the
  Heisenberg group are rarely used in harmonic analysis. Researchers
  often prefer to deal with some spaces of analytic functions. This
  is easily achieved by intertwining \(\Lambda_\myhbar\)
  with an operator of multiplication:
  \begin{equation}
    \label{eq:peeling-map}
    \oper{E}_d: f(x,y) \mapsto \rme^{d(x,y)} f(x,y)\,
  \end{equation}
  by the exponent of a function \(d(x,y):\, \Space{R}{2n}\rightarrow \Space{C}{}\).
  In an obvious way, \(\oper{E}_d\) is unitary as a map:
  \begin{displaymath}
    \oper{E}_d: \FSpace{L}{2}(\Space{R}{2n}) \rightarrow
  \FSpace{L}{2}(\Space{R}{2n}, \rme^{-2\Re d(x,y)} \,\rmd x\,\rmd y)
  \end{displaymath}
  for the appropriate weighted measures. We call \(\oper{E}_d\)
  \emph{peeling}~\cite{Kisil13a} and will discuss the choice of
  \(d(x,y)\) in Example~\ref{ex:gaussina-peeling}.
\end{rem}

\section{Covariant and contravariant transforms}
\label{sec:covar-amd-contr}

 \subsection{The Fourier--Wigner and the covariant transforms}
 \label{sec:four-wign-transf}

Let \(\uir{}{}\) be a strongly continuous unitary representation of a
group \(G\) in a Hilbert space  \(\FSpace{H}{}\), a pair of vectors \(u\),
\(v\in \FSpace{H}{}\) defines the \emph{matrix coefficients} of \(\uir{}{}\):
\begin{equation}
  \label{eq:matrix-coefficients}
  \oper{W}(u,v)(g)=\scalar{u}{\uir{}{}(g)v},\qquad g\in G,
\end{equation}
which is a continuous bounded function on \(G\). Since the
construction occurs in many seemingly disconnected contexts there is an
extensive list of its names: wavelet transform~\cite{AliAntGaz14a}, voice
transform~\cite{FeichtingerPap14a}, coherent state
transform~\cite{Perelomov86}, covariant transform~\cite{Berezin72},
etc. To add a diversity we recall that all sorts of special
functions are matrix coefficients for various group
representations~\cite{Vilenkin68,VilenkinKlimyk95a,%
BindenharnGustafsonLoheLouckMilne84a}. Furthermore, many classical function spaces can be defined as coorbits~\cite{FeichtingerGrochenig88a} with natural norms transported through  covariant transforms to new settings~\cites{Kisil12d,Kisil13a}. 

For the case of \(G=\Space{H}{n}\),
the elementary action of the centre \(Z\) can be omitted. The
corresponding matrix coefficients for the Schr\"odinger representation~\eqref{eq:schroedinger-rep} are:
\begin{align}
  \label{eq:fourier-wigner}
  \oper{W} (f,\phi)(x,y) &=\scalar{f}{\uir{}{\myhbar}(0,x,y)\phi}\\
  &=\int_{\Space{R}{n}}\textstyle
     \rme^{2\pi \rmi \myhbar
                            y(t-\frac{1}{2}x)}\, f(t)\,
    \overline {\phi}(t-x)\,\rmd t\nonumber \\
  \nonumber 
  &=\scalar{\uir{}{\myhbar}(0,-\half x,-\half y)f}{\uir{}{\myhbar}(0,\half x, \half y)\phi}\\
  &=\int_{\Space{R}{n}}\textstyle
     \rme^{2\pi \rmi \myhbar
     yt'}\, f(t'+\frac{1}{2}x)\, \overline {\phi}(t'-\frac{1}{2}x)\, \rmd
    t'.\nonumber 
\end{align}
This was called the \emph{Fourier--Wigner transform}
in~\cite{Folland89}*{\S~I.4}, since the Fourier transform
\((x,y)\mapsto (q,p)\) of \(\oper{W} (\phi,\phi)(x,y)\) produces the
\emph{Wigner transform} of a quantum state \(\phi\), see
also~Prop.~\ref{pr:kernal-as-PDO-Wigner}.

The Fourier--Wigner transform is initially defined on
\(\FSpace{L}{2}(\Space{R}{n})\times \FSpace{L}{2}(\Space{R}{n}) \) and
can be naturally extended by linearity to a unitary map on
\(\FSpace{L}{2}(\Space{R}{n})\otimes \FSpace{L}{2}(\Space{R}{n}) \sim
\FSpace{L}{2}(\Space{R}{2n})\)%
.
In particular, \(\oper{W}\) posses the crucial property of
``sesqui-unitary''~\cite{Folland89}*{(1.42)}:
\begin{equation}
  \label{eq:wigner-sesqui-linear}
  \scalar[\Space{R}{2n}]{\oper{W} (f_1,\phi_1)}{\oper{W}
    (f_2,\phi_2)}=\scalar[\Space{R}{n}]{f_1}{f_2}
  \overline{\scalar[\Space{R}{n}]{\phi_1}{\phi_2}}\,.
\end{equation}
This identity is also known as the \emph{orthogonality relation} and valid for a more general coherent states transform~\cite{AliAntGaz14a}*{\S~8.2}. Also, the identity~\eqref{eq:wigner-sesqui-linear} does not contain anything specific for the Schr\"odinger representation~\eqref{eq:schroedinger-rep} on \( \FSpace{L}{2}(\Space{R}{n}) \) and holds for a unitary irreducible representation of \(\Space{H}{n}\) on a Hilbert space \(\FSpace{H}{}\).
For a fixed non-zero
\(\phi\in\FSpace{L}{2}(\Space{R}{n})\) the linear map:
\begin{equation}
  \label{eq:covariant-heisenberg-defn}
  \covar{\phi }:  \FSpace{L}{2}(\Space{R}{n}) \rightarrow
  \FSpace{L}{2}(\Space{R}{2n}): \ f \mapsto \oper{W}(f,\phi )
\end{equation}
is a particular case of the general construction of \emph{covariant
  transform}.  Some other names (e.g. coherent states transform, wave packet
expansion, voice transform, etc.) are also used in the literature. In physical language the overcomplete system of vectors
\(\phi_{(x,y)}=\uir{}{\myhbar}(x,y) \phi\) are known as \emph{Schr\"odinger coherent
  states} generated by the \emph{ground} (or \emph{vacuum}, or \emph{fiducial}, etc.) state \(\phi\).

Let \(\uir{}{\myhbar}\) be an irreducible unitary representation of \(\Space{H}{n}\) in a Hilbert space \(\FSpace{H}{}\) (for example, the Schr\"odinger representation~\eqref{eq:schroedinger-rep} on \( \FSpace{L}{2}(\Space{R}{n}) \)) and \(\theta \in \FSpace{H}{}\) be a unit vector. The main properties of the covariant transform~\eqref{eq:covariant-heisenberg-defn} follow
from~\eqref{eq:wigner-sesqui-linear} and can be summarised as
follows~\citelist{ \cite{AliAntGaz14a}*{Thm.~7.3.1}
  \cite{Perelomov86}*{\S~1.2}}.
\begin{enumerate}
\item The map \(\covar{\theta} \)~\eqref{eq:covariant-heisenberg-defn} is a unitary operator
  into the space
  \begin{equation}
    \label{eq:image-space-coav-tr}
    \FSpace[\theta]{F}{}\coloneqq \{\covar{\theta} (f): \ f\in
    \FSpace{H}{}\},
  \end{equation}
  which is a closed subspace of \(\FSpace{L}{2}(\Space{R}{2n})\).
\item \(\FSpace[\theta]{F}{}\) is invariant and irreducible under the pulled action
  \(\Lambda_\myhbar\)~\eqref{eq:left-action-pulled} of
  \(\Space{H}{n}\).
\item The map \(\covar{\theta} \) intertwines \(\uir{}{\myhbar}\) and
  \(\Lambda_\myhbar\)~\eqref{eq:left-action-pulled} restricted to
  \(\FSpace[\theta]{F}{}\):
  \begin{equation}
    \label{eq:left-shift-itertwine}
    \covar{\theta} \circ \uir{}{\myhbar} = \Lambda_{\myhbar}\circ \covar{\theta} \,.
  \end{equation}
\item The orthogonal \emph{FSB projection} \(\oper{P}_\theta: \FSpace{L}{2}(\Space{R}{2n})
  \rightarrow \FSpace[\theta]{F}{}\) is an integral operator:
  \begin{align}
    \label{eq:reproducing-projection}
    [\oper{P}_\theta f](x,y)&=\scalar{f}{\Lambda_{\myhbar}(0,x,y)\Theta}\\
    \nonumber 
              &=\int_{\Space{R}{2n}} f(x',y')\,
    \rme^{-\pi \rmi \myhbar(xy'-yx')} \overline{\Theta}(x'-x,y'-y)\,\rmd x'\,\rmd y',\quad 
  \end{align}
  where \(\Theta= \covar{\theta} (\theta)\in \FSpace[\theta]{F}{}\) and
  \(f\in \FSpace{L}{2}(\Space{R}{2n})\).  An immediate consequence is
  the following commutativity:
  \begin{equation}
    \label{eq:Bargmann-proj-commutes-left-action}
    \oper{P}_\theta \circ \Lambda_{\myhbar}(0,x,y)
    = \Lambda_{\myhbar}(0,x,y) \circ \oper{P}_\theta\,.
  \end{equation}
\item \label{it:cyclicity}
  The vector \(\Theta\) is \emph{cyclic} in \(\FSpace[\theta]{F}{}\) in the sense that \(\FSpace[\theta]{F}{}\) is the closed linear span of all \(\Lambda_{\myhbar}(0,x,y)\Theta\), for \((x,y)\in\Space{R}{2n}\).
\item \label{it:reproducing-property}
  In particular,  \(\FSpace[\theta]{F}{}\) is a \emph{reproducing kernel} Hilbert space
  with the reproducing kernel
  \begin{align}
    \label{eq:reproducing-kernel}
    K_{(x,y)}(x',y')&= \overline{\Lambda_{\myhbar}(0,x,y)  \Theta
    (x',y')}
                      =\scalar{\uir{}{\myhbar}(x',y')
                       \theta}{\uir{}{\myhbar}(x,y) \theta}\\
    \nonumber 
                    &
                      =\Lambda_{\myhbar}(0,x',y')  \Theta(x,y)\,.
  \end{align}
  For reasons elaborated in Example~\ref{ex:gaussian-fsb}, cf. also~\cite{Folland89}*{\S~4.5}, the important   case is the ground state \(\theta\) defined by the \emph{squeezed Gaussian}
  \(\varphi_{\tau}(t)=2^{n/4}\rme^{-\pi \myhbar t^2/\tau}\) for some constant
  \(\tau>0\). Then:
  \begin{align}
    \label{eq:FSB-Gaussian}
    \Phi_{\tau}(x,y)& \coloneqq \oper{W}_{\phi_\tau}(\phi_\tau)=\rme^{-\pi \myhbar /(2\tau )(x^2+\tau ^2 y^2)}\,, \\
    \label{eq:FSB-repro-kernel}
    K_{(x,y)}(x',y')&= \rme^{-\pi \rmi
      \myhbar(xy'-yx')} \rme^{-\pi\myhbar/(2\tau )((x'-x)^2+\tau ^2 (y'-y)^2)}\,.
  \end{align}
\item \label{it:contravariant-adjoint}
  The contravariant transform
  \(\contravar[\uir{}{}]{\psi}: \FSpace[\theta]{F}{} \rightarrow \FSpace{H}{}\) is provided the reconstruction formula:
  \begin{equation}
    \label{eq:reconstruction-formula}
    \contravar[\uir{}{}]{\psi}: f \mapsto \int_{\Space{R}{2n}} f(x,y)\,
      \uir{}{}(x,y)\,\psi\,\rmd x\,\rmd y  \qquad \text{ where }
      f\in \FSpace[\theta]{F}{}\,,
    \end{equation}
    see Defn.~\ref{de:contravariant-transform} below for details.  By the orthogonality relation~\eqref{eq:wigner-sesqui-linear} the identity on \(\FSpace{H}{}\):
     \begin{equation}
      \label{eq:wave-trans-inverse}
      \contravar{\psi} \circ \covar{\theta} = \scalar{\psi}{\theta} I\,.
    \end{equation}
\end{enumerate}

\begin{rem}
    Since the Planck constant \(\myhbar\) is fixed within this paper (and, possibly, in the physical reality as well) we are not indicating it in our notation, e.g. \(\Phi_\tau\). On the other hand, we will consider several different values of the squeeze parameter \(\tau\) simultaneously, in most cases we denote two such values as \(\tau\) and \(\varsigma\). To make our notation simpler the dependence on \(\phi_\tau\) will be indicated by \(\tau\) alone, e.g. \(\FSpace[\tau]{F}{} \coloneqq \FSpace[\phi_\tau]{F}{}\), \( \contravar{\tau} \coloneqq \contravar{\phi_\tau} \), etc.
\end{rem}

\subsection{Irreducible components for the pulled action}
\label{sec:covar-transf-irred}

Because the representation \(\Lambda_\myhbar\)~\eqref{eq:left-action-pulled} is not too far from the left regular representation it  is highly reducible on \(\FSpace{L}{2}(\Space{R}{2n})\). For example. this follows from the observation that \(\Lambda_\myhbar\) commutes with \(R_\myhbar(0,x,y)\)~\eqref{eq:right-action-pulled}, for all \((x,y)\in\Space{R}{2n}\) and the latter operators are not scalar multiplies of the identity on~\(\FSpace{L}{2}(\Space{R}{2n})\).

The intertwining property \(\covar[]{\theta } \circ \uir{}{\myhbar} = \Lambda_{\myhbar}\circ \covar[]{\theta } \)~\eqref{eq:left-shift-itertwine} implies that \(\FSpace[\theta ]{F}{}\)~\eqref{eq:image-space-coav-tr} is invariant. Furthermore, the irreducibility of the Schr\"odinger representation implies irreducibility of the restriction of \(\Lambda_\myhbar\) to \(\FSpace[\theta ]{F}{}\). Together with the standard properties of unitary representations we conclude. 
\begin{lem}
  For any \(\theta \in\FSpace{L}{2}(\Space{R}{n})\), the both spaces \(\FSpace[\theta ]{F}{}\) and its orthogonal complement \(\FSpace[\perp]{F}{\theta }\) are \(\Lambda_\myhbar\) invariant subspaces of \(\FSpace{L}{2}(\Space{R}{2n})\). The space \(\FSpace[\theta ]{F}{}\) is irreducible and its complement \(\FSpace[\perp]{F}{\theta}\) is not.
\end{lem}
In the obvious way, two orbits of a group actions either are disjoint or coincide.
We can obtain a similar conclusion for two \(\Lambda_\myhbar\)-invariant subspaces.
\begin{lem}
  Two irreducible components \(\FSpace[\theta]{F}{}\) and  \(\FSpace[\psi]{F}{}\) have non-trivial intersection if and only if \(\theta\) and \(\psi\) are linearly dependent. In the latter case \(\FSpace[\theta]{F}{} = \FSpace[\psi]{F}{}\). 
\end{lem}
\begin{proof}
  Sufficiency is trivial. For necessity, without loss of generality we assume that \(\norm{\theta}=\norm{\psi}=1\) and a non-zero common element \(F \in  \FSpace[\theta]{F}{} \cap \FSpace[\psi]{F}{}\) has unit norm. Then, by the unitarity of the covariant transform  \(F=\covar{\theta} f_1\) and \(F=\covar{\psi} f_2\) for some \(f_{1,2}\in\FSpace{L}{2}(\Space{R}{n})\) with unit norm. Then from sesqui-unitarity~\eqref{eq:wigner-sesqui-linear} and the Cauchy--Schwartz inequality:
  \begin{displaymath}
    1=\scalar[\Space{R}{2n}]{\covar{\theta} f_1}{\covar{\psi} f_2} = \scalar[\Space{R}{n}]{f_1}{f_2}
  \overline{\scalar[\Space{R}{n}]{\theta }{\psi}} \leq \norm{f_1}\cdot \norm{f_2} \cdot \norm{\theta } \cdot \norm{\psi} =1\,.
\end{displaymath}
The identity is achieved only if \(\theta =\lambda \psi\) and \(f_1=\lambda f_2\) for an unimodular \(\lambda\in\Space{C}{}\).
\end{proof}
A systematic usage of representation theory requires the following concept of the induced covariant transform~\cite{Kisil11c,Kisil13a}, which is based on Perelomov's coherent states~\cite{Perelomov86}.
\begin{defn}
  \label{de:covariant-transform}
  Let  \(\uir{}{}\) be a unitary representation of a group \(G\)  in a space \(\FSpace{H}{}\). Let there exist a common eigenvector  \(\phi \in\FSpace{H}{}\) for all operators \(\uir{}{}(h)\) for all \(h\in H\) for some subgroup \(H\subset G\). For a fixed section \(\map{s}: G/H \rightarrow G\) the \emph{induced covariant transform} \(\covar[\rho]{\phi}: \FSpace{H}{} \rightarrow \FSpace{L}{}(G/H)\) is
  \begin{equation}
    \label{eq:covariant-transform-defn}
    \covar[\rho]{\phi} \,:   f \mapsto [\covar[\rho]{\phi} f](x)=\scalar{f}{\uir{}{}(\map{s}(x)) \phi}\,,\qquad \text{where } f\in \FSpace{H}{} \text{ and }  x\in G/H.
  \end{equation}
\end{defn}
Here \(\FSpace{L}{}(G/H)\) is a certain linear space of functions on \(G/H\) depending on properties of \(\phi\). For example, for admissible \(\phi\) it is \(\FSpace{L}{2}(G/H)\), however other cases are of interest as well~\cites{Kisil98a,Kisil13a,Kisil12d}. 

One can check that \eqref{eq:covariant-heisenberg-defn} is a special case of~\eqref{eq:covariant-transform-defn} for the Heisenberg group \(G=\Space{H}{n}\), its centre as the subgroup \(H\) and the Schr\"odinger representation. Then  \(\phi \in\FSpace{L}{2}(\Space{R}{n})\) is a common eigenvector  for all operators \(\uir{}{}(s,0,0)\),  \(s \in \Space{R}{}\).  

The adjoint notion to covariant transform is:
\begin{defn}
  \label{de:contravariant-transform}
  Let  \(\uir{}{}\) be a unitary representation of a group \(G\)  in a Hilbert space \(\FSpace{H}{}\). For some subgroup \(H\subset G\) and a fixed Borel section \(\map{s}: G/H \rightarrow G\) the \emph{contravariant transform} \(\contravar[\rho]{\psi}: \FSpace{L}{}(G/H) \rightarrow \FSpace{H}{}\) is
  \begin{equation}
    \label{eq:contravariant-transform-defn}
    \contravar[\rho]{\psi}: \ f(x) \mapsto \contravar[\rho]{\psi} (f) = \int_{G/H} f(x) \, \uir{}{}(\map{s}(x)) \psi\,\rmd x\,.
  \end{equation}
  Here \(\psi\in \FSpace{H}{}\) is known as \emph{reconstruction vector} and integration is performed over an invariant measure on \(G/H\).
\end{defn}
It is easy to check the relation between covariant and contravariant transforms:
\begin{equation}
  \label{eq:covariant-adjoint-contravariant}
  (\covar[\rho]{\phi})^* =\contravar[\rho]{\phi} \,.
\end{equation}

\begin{rem}
  To unload our notations we will continue to denote the co-/con\-tra\-vari\-ant transforms generated by the Schr\"odinger representation by \(\covar{\psi}\) and \(\contravar{\psi}\) respectively. On the contrast, co-/contravariant transforms for the left pulled action will be denoted as \(\covar[\Lambda ]{\psi}\) and \(\contravar[\Lambda ]{\psi}\) (omitting a reference to the fixed Planck constant \(\myhbar\)).
\end{rem}

Besides the generic relation \((\covar[\Lambda ]{\phi})^* =\contravar[\Lambda]{\phi} \)~\eqref{eq:covariant-adjoint-contravariant} there is a specific connection for the pulled action of the Heisenberg group:
\begin{equation}
  \label{eq:covariant-contravariant-left-right}
 \covar[\Lambda]{\Psi} =\contravar[\Lambda]{{\oper{R}\overline{\Psi}}} \,,
\end{equation}
where \( [\oper{R}f](x,y)=f(-x,-y)\) as before.
Indeed, using  identities~\eqref{eq:conjugation-left-right}--\eqref{eq:right-through-left} between the left and right pulled actions we find:
\begin{align*}
  [\covar[\Lambda]{\Psi} \Theta ](x,y)&=\scalar{\Theta}{\Lambda_\myhbar(0,x,y) \Psi}\\
  &=\int_{\Space{R}{2n}} \Theta(x',y')\, \overline{\Lambda_\myhbar(0,x,y) \Psi(x',y')}\,\rmd x' \, \rmd y'\\
  &=\int_{\Space{R}{2n}} \Theta(x',y')\, R_\myhbar(0,-x,-y) \overline {\Psi}(x',y')\,\rmd x' \, \rmd y'\\
  &=\int_{\Space{R}{2n}} \Theta(x',y')\, \Lambda_\myhbar(0,x',y') [\oper{R}\overline{\Psi}](x,y)\,\rmd x' \, \rmd y'\\
  &=[\contravar[\Lambda]{{\oper{R}{\overline{\Psi}}}} \Theta](x,y)\,.
\end{align*}
\begin{cor}
  For \(f\), \(\Psi \in \FSpace{L}{2}(\Space{R}{2n})\) we have \(\covar[\Lambda]{\Psi}  f \in \FSpace[\oper{R}\overline{\Psi }]{F}{}\)\,.
\end{cor}
In particular, since \(\oper{R}\overline{\Phi }=\Phi \), we recover already seen result: \(\covar[\Lambda]{\Phi_\tau}  f \in \FSpace[\tau]{F}{}\) for all \(f\in \FSpace{L}{2}(\Space{R}{2n})\).
\begin{cor}
  For \(f\in \FSpace[\theta ]{F}{}\), and \(\Psi \in \FSpace{L}{2}(\Space{R}{2n})\) we have
  \begin{equation}
    \covar[\Lambda]{\Psi}  f =\covar[\Lambda]{\Psi_1} f\, , \qquad \text{ where } \Psi_1= \oper{P}_\theta  \Psi\,.
  \end{equation}
\end{cor}
\begin{proof}
  Let \(\Psi = \Psi_1+\Psi_0\) for \(\Psi_1\in\FSpace[\tau]{F}{}\) and \(\Psi_0 \in  \FSpace[\perp]{F}{\phi}\). Then from the \(\Lambda_\myhbar\)-invariance of \(\FSpace[\perp]{F}{\phi}\) follows that \(\scalar{f}{\Lambda_\myhbar(0,x,y) \Psi_0}=0\)
 for all \((x,y)\in\Space{R}{2n}\). Thus:
  \begin{displaymath}
    [\covar[\Lambda]{\Psi}  f](x,y) = \scalar{f}{\Lambda_\myhbar(0,x,y) \Psi}  = \scalar{f}{\Lambda_\myhbar(0,x,y) \Psi_1} = [\covar[\Lambda]{\Psi_1} f](x,y)\,. 
  \end{displaymath}
\end{proof}

A \(\Lambda_\myhbar\)-covariant unitary operator \(\oper{U}_{\theta \psi}: \FSpace[\theta]{F}{} \rightarrow \FSpace[\psi]{F}{}\) can be constructed as the composition \(U= \covar{\psi}\circ \contravar{\theta}\) of contra- and covariant transforms:
\begin{align}
  \nonumber 
  [\oper{U}_{\theta \psi} f] (x,y)&= [\covar{\psi}\circ \contravar{\theta} f] (x,y)\\
  \nonumber 
             &= \scalar{\contravar{\theta} f}  {\uir{}{\myhbar}(0,x,y) \psi} \\
  \nonumber 
             &= \scalar{ \int_{\Space{R}{2n}} f(x',y')\, \uir{}{\myhbar}(0,x',y') \theta \, \rmd x' \,\rmd y'} {\uir{}{\myhbar}(0,x,y) \psi} \\
  \label{eq:intertwining-components}
  &= \int_{\Space{R}{2n}} f(x',y') \scalar{ \uir{}{\myhbar}(0,x',y') \theta} {\uir{}{\myhbar}(0,x,y) \psi} \, \rmd x' \,\rmd y'\,,
\end{align}
assuming we can change the order of integration, e.g. due to the Fubini theorem.
The last formula includes the reproducing kernel~\eqref{eq:reproducing-kernel} as the special case for \(\psi=\theta\).
 For \(\psi\in\FSpace[\theta]{F}{} \), operator \(\oper{U}_{\theta \psi}\) is a multiple of identity by the Schur lemma and irreducibility of the pulled action \(\Lambda_\myhbar\) on \(\FSpace[\theta]{F}{}\).
\begin{example}
  \label{ex:intertwining-kernel}
  We may consider the map \(\FSpace[\tau]{F}{} \rightarrow \FSpace[\varsigma]{F}{}\) as a composition \(\covar{\varsigma} \circ \contravar{\tau}\) of the contravariant and covariant transform for the respective reconstructing \(\phi_\tau\) and analysing \(\phi_\varsigma\) vectors. It turn out to be an integral operator:
  \begin{align}
    \nonumber 
    f_\varsigma(x,y) &=\sqrt[4]{\frac{\tau \varsigma }{(\tau+\varsigma)^2}}\int_{\Space{R}{2n}} f_\tau(x_1,y_1)\, \\
    \nonumber 
    & \qquad \qquad \rme^{-\frac{\pi \myhbar}{\tau + \varsigma} ((x-x_1)^2+ \tau \varsigma (y-y_1)^2 + 2\rmi (\varsigma x y_1- \tau x_1 y)+\rmi(xy+x_1 y_1)(\varsigma -\tau))} \, \rmd x_1\,\rmd y_1\\
    \nonumber 
           &=\int_{\Space{R}{2n}} f_\tau(x_1,y_1)\, \Lambda_h(x_1,y_1)\Phi_{\tau\varsigma}(x,y) \, \rmd x_1\,\rmd y_1
    =[\contravar[\Lambda_h]{\Phi_{\tau\varsigma}}f_\tau](x,y)
    \\
    \label{eq:intertwining-as-covariant}
           &=\int_{\Space{R}{2n}} f_\tau(x_1,y_1)\, \overline{\Lambda_h(x,y)\Phi_{\tau\varsigma}(x_1,y_1)} \, \rmd x_1\,\rmd y_1
    =[\covar[\Lambda_h]{\Phi_{\tau\varsigma}}f_\tau](x,y)
    \,.
  \end{align}
  where \(f_\tau \in \FSpace[\tau]{F}{}\),  \(f_\varsigma \in \FSpace[\varsigma]{F}{}\) and
  \begin{equation}
    \label{eq:Phi-tau-sigma}
    \Phi_{\tau\varsigma}(x,y) = \scalar{\phi_\varsigma}{\uir{}{}(g)\phi_\tau}
    = \sqrt[4]{\frac{\tau \varsigma }{(\tau+\varsigma)^2}}
    \rme^{-\pi \myhbar / (\tau + \varsigma) (\rmi xy (\tau-\varsigma)+x^2+\tau\varsigma y^2)}\,.
  \end{equation}
\end{example}
Obviously, for \(\varsigma = \tau\) the identity~\eqref{eq:intertwining-as-covariant} coinsides with the reproducing formula~\eqref{eq:reproducing-projection} on \(\FSpace[\tau]{F}{}\). We will call
\begin{equation}
  \label{eq:intertwining-kernel}
  K_{\tau\varsigma}(x_1,y_1;x,y) = \overline{\Lambda_h(x_1,y_1)\Phi_{\tau\varsigma}(x,y)}
\end{equation}
the \emph{intertwining kernel}. Also, the function \(\Phi_{\tau\varsigma}\)~\eqref{eq:Phi-tau-sigma} will be named \emph{mixed Gaussian}.  We present its main properties in Lem.~\ref{le:intertwining-kernel-properties} once some additional notions will be introduced.

\subsection{Twisted convolution and symplectic Fourier transform}
\label{sec:sympl-four-transf}

Many operators in analysis (e.g.~\eqref{eq:reconstruction-formula} and
many further examples in~\cites{Kisil13a,Kisil94e}) appear as
\emph{integrated representations} or \emph{relative convolution}
\begin{equation}
  \label{eq:integrated-rep}
  \uir{}{}(k)=\int_{X} k(x)\,\uir{}{}(\map{s}(x))\,\rmd x
\end{equation}
for a suitable representation \(\uir{}{}\) of a group \(G\), its
homogeneous space \(X=G/H\), a Borel section \(\map{s}:
X\rightarrow G\) and a
function (distribution) \(k\) defined on
\(X\)~\cite{Kisil94e,Kisil13a}. As usual, the above integral is
understood in the weak sense. In particular, \eqref{eq:contravariant-transform-defn} can be restated as \(\contravar{\psi} (f)=\uir{}{}(f)\psi\).

For an irreducible representation \(\uir{}{}\), a composition of two
integrated representations is an integrated representation again
\begin{equation} 
  \label{eq:twisted-conv-defn}
  \uir{}{}(k_1)\circ \uir{}{}(k_2) =
  \uir{}{}(k_1\twist k_2),
\end{equation}
for a kernel \( k_1\twist k_2\).
\begin{defn}
  The function (distribution) \( k_1\twist k_2\) defined
  by~\eqref{eq:twisted-conv-defn} is called the \emph{twisted
    convolution} of \(k_1\) and \(k_2\) (produced by the
  representation \(\uir{}{}\)).
\end{defn}
Note, that the twisted convolution~\eqref{eq:twisted-conv-defn}
depends on the equivalence class of \(\uir{}{}\) but \emph{not} on its
specific realisation.  In particular, for the left
\(\Lambda_\myhbar\)~\eqref{eq:left-action-pulled}, right
\(R_\myhbar\)~\eqref{eq:right-action-pulled} and the Schr\"odinger
representation \(\uir{}{\myhbar}\)~\eqref{eq:schroedinger-rep} the
twisted convolution on the homogeneous space \(\Space{H}{n}/Z\) is
given by~\citelist{\cite{Howe80b} \cite{Folland89}*{\S~1.3}}
\begin{align}
  \nonumber
  \lefteqn{(f_1 \twist f_2)(x,y)=  [\Lambda_\myhbar(f_1) f_2](x,y)}\qquad&\\
  \label{eq:twisted-heisenberg}
  &=\int_{\Space{R}{2n}}f_1(x',y')\,f_2(x-x',y-y')\,\rme^{\pi\myhbar\rmi (x'y-y'x)}\,dx'\,dy'\,. 
\end{align}
More accurately, we shall account the dependence on the Planck
constant and denote the twisted convolution by
\(f_1 \twist_{\myhbar} f_2\). However, the Planck constant
\(\myhbar\) is systematically unloaded from notations throughout this paper. As can be expected for an important concept, \eqref{eq:twisted-heisenberg} is known under many different names, e.g. \emph{Moyal star product}~\cites{Moyal49,Zachos02a}.

Using unitarity of the left regular representation and
\eqref{eq:conjugation-left-right} (alternatively, make a change of
variables) we can express the twisted convolution through the
integrated right regular representation:
\begin{equation}
  \label{eq:twisted-as-right-integrated}
  \begin{split}
    f_1 \twist f_2&=  \Lambda_\myhbar(f_1) f_2\\
&=R_\myhbar(\breve{f}_2) f_1, \quad \text{where } \breve{f}_2(x,y) \coloneqq [\oper{R}{f}_2](x,y)=f_2(-x,-y)\,.
  \end{split}
\end{equation}
\begin{example}
  Note that the FSB projection \(\oper{P}_\tau:
  \FSpace{L}{2}(\Space{R}{2n})\rightarrow
  \FSpace[\tau]{F}{}\)~\eqref{eq:reproducing-projection} is
  the twisted convolution with the Gaussian:
  \begin{equation}
    \label{eq:FSB-projection-as-twisted-convolution}
    \oper{P}_\tau f= f\twist \Phi_\tau\,.
  \end{equation}
\end{example}
For reasons emerging soon the following transformation can be a suitable
modification of the Fourier transform.
\begin{defn}
  The \emph{symplectic Fourier transform} (again depending on
  \(\myhbar\))~\citelist{\cite{Folland89}*{\S~2.1}
    \cite{Howe80b}*{\S~2.1}} is defined by:
  \begin{align}
    \label{eq:symplect-Fourier-defn}
                             \wideparen{\psi}(x,y)&= (\myhbar/2)^{n}\int_{\Space{R}{2n}} {\psi}(x',y')\, \rme^{\pi \rmi \myhbar
                                                    (x'y-y'x)}\,\rmd x'\,\rmd y'  \,.
  \end{align}
\end{defn}
A systematic usage of the symplectic Fourier transform in the context of the FSB spaces can be seen in~\cite{Coburn99,Coburn01b}. 
\begin{example}[Symplectic Fourier transform of the Gaussian]
  For the Gaussian \(\Phi(x,y)=\rme^{-\pi \myhbar /(2\tau )(x^2+\tau ^2
    y^2)}\)~\eqref{eq:FSB-Gaussian} the standard calculation shows:
  \begin{equation}
    \label{eq:sympl-Fourier-Gaussian}
    \wideparen{\Phi}(x,y)=\Phi(x,y)\,.
  \end{equation}
  An alternative demonstration based on group representations is given in
  Example~\ref{ex:gaussian-fsb}.
\end{example}
The symplectic Fourier transform is a reflection, that is of order two \((\wideparen{\psi})\wideparen{\phantom{a}}=\psi\) in contrast to
the usual Fourier transform, which is of order four:
\((\hat{\psi})\hat{\phantom{a}s}(t)=\psi(-t)\).  An ``intriguing
fact''~\cite{Howe80b}*{\S~2.1} is that the symplectic Fourier
transform can be represented (up to a scaling) as integrated
representations of \(\Lambda_\myhbar(0,x',y')\) and
\(R_\myhbar(0,x',y')\) applied to the constant function
\(\mathbf{1}(x,y)\equiv 1\):
\begin{equation}
  \label{eq:symplect-Fourier-left}
  \begin{split}
  \wideparen{a} = a\twist \mathbf{1}
  &= \Lambda_\myhbar(a)\,  \mathbf{1}
    = R_\myhbar(a)\, \mathbf{1}\\
   &= R_\myhbar( \mathbf{1})\, a
  = \Lambda_\myhbar(  \mathbf{1})\, a
                                      \,,
  \end{split}
\end{equation}
where the last line is based on~\eqref{eq:twisted-as-right-integrated}.
Some immediate consequences are:
\begin{enumerate}
\item Commutation of left regular representation and the symplectic
  Fourier transform:
\begin{equation}
  \label{eq:sympl-Fourier-intertwin-left}
  \wideparen{\ } \circ \Lambda_\myhbar (0,x,y)
  = \Lambda_\myhbar (0,x,y) \circ  \wideparen{\ }\,.
\end{equation}
\item The intertwining property for the right regular representation:
\begin{equation}
  \label{eq:sympl-Fourier-intertwin-right}
  \wideparen{\ } \circ R_\myhbar (0,x,y)
  = R_\myhbar (0,-x,-y) \circ  \wideparen{\ }\,.
\end{equation}

\item Combination of two previous intertwining properties
  with~\eqref{eq:shift-by-left-right}--\eqref{eq:mult-by-left-right} implies the
  fundamental intertwining properties of the (symplectic)
  Fourier transform:
  \begin{equation}
    \label{eq:fourier-intertw-shift-mult}
    \wideparen{\ } \circ \map{S} (x,y)
    = \map{E}_{\myhbar} (x,y) \circ  \wideparen{\ }\,,\qquad
    \wideparen{\ } \circ \map{E}_{\myhbar} (x,y)
    = \map{S} (x,y) \circ  \wideparen{\ }\,,
  \end{equation}
  where are operators of the Euclidean shift
  \(\map{S}\)~\eqref{eq:shift-defn} and multiplication
  \(\map{E}_{\myhbar}\)~\eqref{eq:exp-mult-defn}.
\item The symplectic Fourier transform of a twisted convolution is:
\begin{equation}
  \label{eq:sympl-Fourier-of-twisted}
  (f_1 \twist f_2)\wideparen{\ }= (f_1 \twist f_2) \twist \mathbf{1}
  = f_1 \twist (f_2 \twist \mathbf{1}) = f_1 \twist \wideparen{f}_2 \,.
\end{equation}

\end{enumerate}
\begin{rem}
  There is a dilemma discussed   in~\cite{Folland89}*{Prologue}: which Fourier transform to use---the
  ordinary or symplectic---in the context of \(\Space{H}{n}\)? In our
  opinion, these maps are clearly distinguished by their ranges. The
  symplectic transform~\eqref{eq:symplect-Fourier-defn} maps function
  on \(\Space{H}{n}/Z\) to functions on the same set. The ordinary
  Fourier transform
  \begin{equation}
    \label{eq:Fourier-Kirillov}
    \hat{\psi}(q,p)=\int_{\Space{R}{2n}} \psi(x,y)\, \rme^{-2\pi\rmi(qx
      +py)}\,\rmd x\,\rmd y
  \end{equation}
  sends a function \(\psi\) on \(\Space{H}{n}/Z\) to a function
  \(\hat{\psi}\) on the coadjoint orbit in the dual space
  \(\algebra{h}^*_{n}\) of the Lie \(\algebra{h}_{n}\) of
  \(\Space{H}{n}\). Here coordinates \((\myhbar,q,p)\) on
  \(\algebra{h}^*_{n}\) have the physical meaning of the Planck
  constant, coordinate and momentum on the phase space, respectively. More
  accurately, \((x,y)\) in~\eqref{eq:Fourier-Kirillov} are coordinates
  on the Lie algebra \(\algebra{h}_{n}\) transferred to
  \(\Space{H}{n}\) by means of the exponential map
  \(\algebra{h}_n\rightarrow \Space{H}{n}\).
  Then, \eqref{eq:Fourier-Kirillov} written in full as
  \begin{equation}
    \label{eq:Fourier-Kirillov-full}
    \hat{\psi}(\myhbar,q,p)=\int_{\Space{H}{n}} \psi(s,x,y)\,
    \rme^{-2\pi\rmi(\myhbar s+qx
      +py)}\,\rmd s\,\rmd x\,\rmd y
  \end{equation}
  is the basic example of the
  \emph{Fourier--Kirillov
    transform}~\cite{Kirillov04a}*{\S~4.1.4}.
\end{rem}

\begin{rem}[Physical units and their mathematical usage]
  \label{re:physical-units}
  The described difference between two transforms is also dictated by
  the physical dimensionality~\citelist{\cite{Kisil17a}*{\S~1.2}
  \cite{Kisil02e}*{\S~2.1}}.
  Let \(M\) be a unit of mass, \(L\)---of length, \(T\)---of
  time. Then coordinate \(q\) is measured in \(L\), momentum \(p\) in
  \(LM/T\) and \(\myhbar\) in their product
  \(L\times LM/T=L^2M/T\)---a unit of the action. Then dual variables
  \(s\), \(x\) and \(y\) are measured in the respective reciprocal
  units \(T/(L^2M)\), \(1/L\) and \(T/(LM)\) respectively. It is easy
  to see that formulae~\eqref{eq:symplect-Fourier-defn},
  \eqref{eq:Fourier-Kirillov}, \eqref{eq:Fourier-Kirillov-full} (as
  well as all other formulae in this paper) follow the following
  rules:
  \begin{enumerate}
  \item \label{it:addition} Only physical quantities of the \emph{same
      dimension} can be added or subtracted. However, there is no
    restrictions on multiplication/division.
  \item \label{it:functions} Therefore, mathematical functions, e.g.
    \(\exp(u)=1+u+u^2/2!+\ldots\) or \(\sin(u)\), can be naturally
    constructed out of a dimensionless number \(u\) only. For example,
    Fourier dual variables, say \(x\) and \(q\), should posses
    reciprocal dimensions because they have to form the expression
    like \(e^{2\pi \rmi xq}\).
  \end{enumerate}
  These rules have a natural physical origin and are
  mathematically valuable as well: validation of formulae is one
  example, the above discussion of ordinary and symplectic Fourier
  transforms---another.

  In the case of the time-frequency analysis units are simpler. A signal is described by a function \(f(t)\), where a variable the ``coordinate'' \(q\) is time in units \(T\), ``momentum'' \(p\) is frequency measured in \(1/T\). The ``Plank constant'' is dimensionless. The dual variables \(s\), \(x\) and \(y\) again have reciprocal units \(1\), \(1/T\) and \(T\) of \(\myhbar\), \(q\) and \(p\).

  The squeeze parameter \(\tau\) shall have units of \(x/y\), that is  \(M/T\) in quantum mechanics and \(T^2\) in time-frequency analysis.
\end{rem}

Finally we establish the connection between the twisted convolution and the
wavelet transform, cf. \cite{Folland89}*{(1.47)}:
\begin{lem} Let \(\phi_1\) and \(\phi_2\) be admissible analysing vectors and \(f_2\) be an admissible reconstructing vector. Then for any \(f_1\in \FSpace{H}{}\):
  \begin{equation}
    \label{eq:twisted-conv-of-wavelet-trans}
  \covar{\phi_1} (f_1) \twist
  \covar{\phi_2} (f_2)=\scalar{f_2}{\phi_1}  \covar{\phi_2} (f_1) \,.
\end{equation}
\end{lem}
\begin{proof}
  Note, that the reconstruction property~\eqref{eq:wave-trans-inverse}
  can be written through the integrated representation as:
  \begin{equation}
    \label{eq:wave-trans-inverse-integrated}
    \uir{}{\myhbar}(\covar{\phi} f) \psi = \scalar{\psi}{\phi} f\,,
    \qquad
    f,\phi,\psi\in \FSpace{L}{2}(\Space{R}{n})\,.
  \end{equation}
  Then, under the lemma's assumptions and for an admissible
  reconstructing vector \(\psi\in \FSpace{L}{2}(\Space{R}{n})\):
  \begin{align*}
    \uir{}{\myhbar}\left(\covar{\phi_1} ( f_1) \twist
    \covar{\phi_2} (f_2)\right)\psi
    &=\uir{}{\myhbar}(\covar{\phi_1} (f_1))
    \uir{}{\myhbar}(\covar{\phi_2} (f_2))\psi\\
    &= \uir{}{\myhbar}(\covar{\phi_1} (f_1)) \scalar{\psi}{\phi_2} f_2\\
    &= \scalar{\psi}{\phi_2} \uir{}{\myhbar}(\covar{\phi_1} (f_1)) f_2\\
    &= \scalar{\psi}{\phi_2} \scalar{f_2}{\phi_1} f_1\\
    &= \scalar{f_2}{\phi_1} \scalar{\psi}{\phi_2}  f_1\\
    &= \scalar{f_2}{\phi_1}  \uir{}{\myhbar}(\covar{\phi_2} (f_1)) \psi\,.
  \end{align*}
  Since \(\uir{}{\myhbar}\) is faithful and \(\psi\) is arbitrary we
  obtain~\eqref{eq:twisted-conv-of-wavelet-trans}.
\end{proof}

\section{Complex variables, analyticity and right shifts}
\label{sec:compl-vari-anal}

The reasons which make the Gaussian a preferred vacuum vector are linked to the complex structure and analyticity. 

\subsection{Right shifts and analyticity}
\label{sec:right-shifts-analyt}

Here we use the following
general result from~\cite{Kisil11c}*{\S~5}, see also further
developments in
\cites{Kisil12d,Kisil13a,Kisil13c,Kisil17a,AlmalkiKisil18a,AlameerKisil21a}.

Let \(G\) be a locally compact group and \(\uir{}{}\) be its
representation in a Hilbert space \(\FSpace{H}{}\). Let
\([\covar{\theta} v](g)=\scalar{v}{\uir{}{}(g) \theta }\) be the wavelet
transform defined by a vacuum state \(\theta \in \FSpace{H}{}\).  Then, the
right shift \(R(g): [\covar{\theta} v](g')\mapsto [\covar{\theta} v](g'g)\)
for \(g\in G\) coincides with the wavelet transform
\([\covar{\theta _g} v](g')=\scalar{v}{\uir{}{}(g')\theta _g}\) defined by the
vacuum state \(\theta _g=\uir{}{}(g) \theta \).  In other words, the covariant
transform intertwines%
\index{intertwining operator}%
\index{operator!intertwining} right shifts on the group \(G\) with the
associated action \(\uir{}{}\) on vacuum states,
cf.~\eqref{eq:left-shift-itertwine}:
\begin{equation}
  \label{eq:wave-intertwines-right}
  R(g) \circ \covar{\theta} = \covar{{\uir{}{}(g)\theta }}.
\end{equation}
This elementary observation has many fundamental consequences.
\begin{cor}[Analyticity of the wavelet transform~\cite{Kisil11c}*{\S~5}] 
  \label{co:cauchy-riemann-integ}
  Let \(G\) be a group and \(dg\) be a measure on \(G\). Let
  \(\uir{}{}\) be a unitary representation of \(G\), which can be
  extended by integration to a vector space \(V\) of functions or
  distributions on \(G\).  Let a mother wavelet \(\theta \in \FSpace{H}{}\) satisfy the
  equation
  \begin{displaymath}
    \int_{G} a(g)\, \uir{}{}(g) \theta \,\rmd g=0,
  \end{displaymath}
  for a fixed distribution \(a(g) \in V\). Then  any wavelet transform
  \(\tilde{v}(g)=\scalar{v}{\uir{}{}(g)\theta }\) obeys the condition:
  \begin{equation}
     \label{eq:dirac-op}
     D\tilde{v}=0,\qquad \text{where} \quad D=\int_{G} \overline{a}(g)\, R(g) \,\rmd g\,,
  \end{equation}
  with \(R\) being the right regular representation of \(G\).
\end{cor}

Some applications (including discrete ones) produced by the \(ax+b\)
group can be found in~\cite{Kisil12d}*{\S~6}, usage in quantum
mechanics is demonstrated in~\cite{AlmalkiKisil18a,AlmalkiKisil19a}.
We turn to the particular case of the Heisenberg group now.
\begin{example}[Gaussian and FSB transform]
  \label{ex:gaussian-fsb}
  Let us consider the squeezed Gaussian\index{Gaussian}
  \(\varphi(t)= \rme^{- \pi\myhbar t^2/\tau}\).  The parameter
  \(\tau\) has the physical dimension of mass times frequency,
  cf. Rem.~\ref{re:physical-units}. In other words, the (reduced)
  Planck constant \(\myhbar\) is in units reciprocal to the product
  \(xy\) and the parameter \(\tau\) is in units reciprocal to the
  ratio \(y/x\).  The physical meaning of \(\tau\) is squeeze
  parameter for coherent states, see \cite{AlmalkiKisil18a,AlmalkiKisil19a} and
  references therein. It is common to put \(\tau=1\)
  in mathematical texts.

  The Gaussian is a null-solution of the operator. cf.~\eqref{eq:schrodinger-derived}. 
  \begin{displaymath}
    \rmd\uir{ -\tau X+\rmi Y}{\myhbar}
    \coloneqq  -\tau \rmd\uir{X}{\myhbar}+\rmi \,\rmd\uir{Y}{\myhbar}
    =    \tau \partial_t+2\pi
  \myhbar t
  \end{displaymath}
For the centre
  \(Z=\{(s,0,0):\ s\in\Space{R}{}\}\subset \Space{H}{}\), we define
  the section \(\map{s}:\Space{H}{}/Z\rightarrow \Space{H}{}\) by
  \(\map{s}(x,y)=(0,x,y)\). Then, the corresponding induced wavelet
  transform~\eqref{eq:covariant-heisenberg-defn} with the measure renormalised by
  the factor \((\myhbar/\tau)^{n/2}\) is:
  \begin{equation}
    \label{eq:coherent-transf-gauss}
    \begin{split}
      \tilde{f}(x,y)&=\scalar{f}{\uir{}{}(\map{s}(x,y))\varphi}\\
      &= \left(\frac{\myhbar}{\tau}\right)^{n/2}
      \int_{\Space{R}{n}} f(t)\, \rme^{\pi \rmi \myhbar (2yt-xy)} \,
      \rme^{-\pi\myhbar (t-x)^2/\tau}\,\rmd t.
    \end{split}
  \end{equation}
  Note, that the normalising factor makes integration dimensionless in
  agreement with the Rem.~\ref{re:physical-units} as well.
  The transformation intertwines the
  Schr\"odinger representation~\eqref{eq:schroedinger-rep} and the left
  pulled action~\eqref{eq:left-action-pulled}. Cor.~\ref{co:cauchy-riemann-integ}
  ensures that the operator, cf.~\eqref{eq:right-action-derived}
  \begin{equation}
    \label{eq:fsb-cauchy-riemann}
    \begin{split}
      \rmd R^{\tau  X+\rmi Y}&= \tau(\pi \rmi \myhbar yI +\partial_x)
      +\rmi (-\pi \rmi \myhbar xI +\partial_y)\\
      &= \pi \myhbar
      (x+\rmi \tau y)I +(\tau \partial_x+ \rmi\partial_y)
    \end{split}
  \end{equation}
  annihilates any \(\tilde{f}(x,y)\) from the image space \(\FSpace[\varphi]{F}{}\)
  of transformation~\eqref{eq:coherent-transf-gauss}.
  Following the idea from~\cite{Howe80b} we use the above intertwining
  properties to evaluate
  \(\Phi(g)=\scalar{\varphi}{\uir{}{}(g)\varphi}\). Indeed, \(\Phi(g)\)
  shall be annihilated by both  \(\rmd R^{\tau  X+\rmi Y}\) and
  \(\rmd \Lambda ^{\tau  X-\rmi Y}\) and therefore is the
  simultaneous null-solution of the operators
  cf.~\eqref{eq:left-action-derived}--\eqref{eq:right-action-derived}: 
  \begin{equation}
    \label{ex:gaussian-annihilators}
    \rmd \Lambda ^{\tau  X-\rmi Y} + \rmd R^{\tau  X+\rmi Y} =
    \tau \pi \rmi \myhbar y + \rmi \partial_y
    \quad\text{and}\quad
    \rmd \Lambda ^{\tau  X-\rmi Y} - \rmd R^{\tau  X+\rmi Y} =
    -\tau \partial_x-\pi  \myhbar x.
  \end{equation}
  This determines
  \(\Phi(g)=\rme^{-\pi \myhbar /(2\tau )(x^2+\tau ^2
    y^2)}\)~\eqref{eq:FSB-Gaussian} up to a constant
  factor. Furthermore, \(\wideparen{\Phi}(x,y)\) shall be the null
  solution of the same operators~\eqref{ex:gaussian-annihilators} as
  \(\Phi(x,y)\) due to the intertwining
  properties~\eqref{eq:sympl-Fourier-intertwin-left}
  and~\eqref{eq:sympl-Fourier-intertwin-right} of the symplectic
  Fourier transform. Thus, \(\wideparen{\Phi}(x,y)=\Phi(x,y)\) once the
  constant factor is confirmed.
\end{example}

\begin{example}[Gaussian and a peeling map]{}
  \label{ex:gaussina-peeling}
  As described in Rem.~\ref{re:peeling-map}, we may look for a peeling
  \(\oper{E}_d: f(x,y)\mapsto
  \rme^{d(x,y)}f(x,y)\)~\eqref{eq:peeling-map}, which shall intertwine
  operator~\eqref{eq:fsb-cauchy-riemann} with the Cauchy--Riemann
  operator.  Then, a simple differential equation implies
  \(d(x,y)=\pi \myhbar/(2\tau) (x^2+\tau^2 y^2)= \frac{\pi}{2}
  \modulus{z}^2\) for
  \(z=\sqrt{\frac{\myh}{2\tau }}(x+\rmi \tau y)\) with \(\myh>0\)
  and \(\tau>0\) (see~\eqref{eq:right-shift-ladder-form} for
  motivation of this normalisation).  Thereafter, the \emph{peeling}
  map~\cite{Kisil13a}
  \begin{equation}
    \label{eq:peeling}
    \tilde{f}(z) \mapsto    F(z)= \rme^{\modulus{z}^2/2}\,\tilde{f}(z)= \rme^{\pi
      \myhbar/(2\tau)(x^2+\tau^2y^2)}\, \tilde{f}(x,y)
  \end{equation}
  produces the function \(F\) satisfying the Cauchy--Riemann equation
  \begin{displaymath}
    \overline{\partial}_zF(z)=(\tau \partial_x+ \rmi\partial_y)
    F(z)=0\,.
  \end{displaymath}
  The composition of the coherent state
  transform~\eqref{eq:coherent-transf-gauss} and the
  peeling~\eqref{eq:peeling} is:
  \begin{align}
    \nonumber 
      \tilde{f}(x,y)&= \left(\frac{\myhbar}{\tau}\right)^{n/2} \rme^{\pi
        \myhbar/(2\tau)(x^2+\tau^2y^2)}\, \int_{\Space{R}{n}}
      f(t)\, \rme^{\pi \rmi \myhbar (2yt-xy)} \,
                      \rme^{-\pi\myhbar (t-x)^2/\tau}\,\rmd t\\
    \nonumber 
      &= \left(\frac{\myhbar}{\tau}\right)^{n/2} \int_{\Space{R}{n}} f(t)\, \rme^{\pi
        \myhbar (- t^2 +2t (x+\rmi \tau y) -(x+ \rmi\tau
        y)^2/2)/\tau } \rmd t
      \\
    \label{eq:fsb-transform}
      &= \left(\frac{\myhbar}{\tau}\right)^{n/2}  \int_{\Space{R}{n}} f(t)\, \rme^{-
        \myh  t^2/(2\tau) +\sqrt{2\myh/ \tau }\cdot  t  z-z^2/2 } \rmd
      t\,,\quad
  \end{align}
  where  \(z=\sqrt{\frac{\myh}{2\tau }}(x+\rmi \tau y)\). Further
  discussion of peeling can be found in~\cites{Alamer19a,AlameerKisil21a}.
\end{example}
The integral~\eqref{eq:fsb-transform} is known as
Fock--Segal--Bargmann (FSB) transform%
\index{Fock--Segal--Bargmann!transform}%
\index{transform!Fock--Segal--Bargmann}. Many sources use this formula
for particular values \(\myhbar=1\) and \(\tau =1\) only.  With
variable value of \(\tau\) \eqref{eq:fsb-transform} becomes
the Fourier--Bros--Iagolnitzer (FBI) transform%
\index{Fourier--Bros--Iagolnitzer transform}%
\index{transform!Fourier--Bros--Iagolnitzer}%
\index{FBI transform!see{Fourier--Bros--Iagolnitzer transform}}%
\index{transform!FBI!see{Fourier--Bros--Iagolnitzer transform}}, see
\cite{Folland89}*{\S~3.3} for introduction. The image
\(\FSpace[\tau]{F}{}\) of \(\FSpace{L}{2}(\Space{R}{n})\) under FSB transform is called the
\emph{Fock--Segal--Bargmann (FSB) space}%
\index{Fock--Segal--Bargmann!space}%
\index{space!Fock--Segal--Bargmann}. More general,
\(\FSpace[\tau]{F}{p}\) is the closed subspace of
\(\FSpace{L}{p}(\Space{C}{n},\Phi_{\myhbar\tau}^{-2}(z)\,\rmd z)\)
consisting of analytic functions. The Stone--von Neumann theorem
implies the following result.
\begin{cor}
  \label{co:peeling-Gaussian}
  The action
  \(\widetilde{\Lambda}_\myhbar (g)=\oper{E}_d \circ \Lambda_\myhbar (g)
  \circ \oper{E}_d^{-1}\) with \(\rme^{-d(z)}=\Phi_\tau(z,\overline{z})\) is a
  unitary irreducible representations of \(\Space{H}{n}\) in the FSB
  space \(\FSpace[\tau]{F}{}\).  Two such actions
  \(\widetilde{\Lambda}_\myhbar\) and \(\widetilde{\Lambda}_{\myhbar'}\) are
  unitary equivalent if and only if \(\myhbar=\myhbar'\).

  Furthermore, \(\widetilde{\Lambda}_\myhbar (g)\) is an invertible isometric
  transformation of
  \(\FSpace[\tau]{F}{p}\rightarrow \FSpace[\tau]{F}{p}\),
cf.~\cite{Zhu11a} or any other space with shift-modulation invariant
norm, e.g~\eqref{eq:johansson-norm}.
\end{cor}

In this paper we do not require a peeling and will refer to the
coherent state transformation~\eqref{eq:coherent-transf-gauss} and its
image \(\FSpace[\varphi]{F}{}\) as \emph{pre-FSB transform} and \emph{pre-FSB space}
respectively. This prefix ``pre-'' is used to distinguish versions of FSB transform with and without peeling, see~\cites{AlameerKisil21a,Alamer19a} for further discussions.

\begin{rem}
  The Gaussian \(\varphi\) is the preferred vacuum state because
  \begin{itemize}
  \item it produces analytic functions through FSB transform;
  \item it is  the minimal uncertainty state.
  \end{itemize}
  Interestingly, the both properties are derived from the identity
  \(\rmd\uir{ -\tau X+ \rmi Y}{\myhbar} \varphi=0\),
  see~\cite{Kisil13c} for further discussion.
\end{rem}

\subsection{Ladder operators and complex variables}
\label{sec:ladd-oper-compl}

The analyticity considered in \S~\ref{sec:right-shifts-analyt}
suggests that a complexification of the derived representation can be useful. To avoid discussion of the complex Lie algebras we postpone a complexification till a particular representation in a complex Hilbert is already selected and complex scalars are consequently  present.

\begin{defn}
  \label{de:ladder-operators}
  Let \(\uir{}{}\) be a unitary irreducible representation of the
  Heisenberg group \(\Space{H}{n}\) with a non-zero Planck constant
  \(\myh=2\pi\myhbar = -\rmi \rmd\uir{S}{}\).  The \emph{ladder operators}, that is the
  pair of the \emph{creation} \(\ladder[\uir{}{},j]{+}\) and
  \emph{annihilation} \(\ladder[\uir{}{},j]{-}\) operators, are the
  following complex linear combination in the derived representation
  of the Weyl algebra:
  \begin{equation}
    \label{eq:ladder-heisenberg}
    \begin{split}
      \ladder[\uir{}{},j]{\pm}
      &= \frac{1}{\sqrt{2\modulus{\myh}\tau}} \left(\pm \tau\,
        \rmd\uir{X_j}{} + {\rmi}\,\rmd\uir{Y_j}{}\right),\qquad \text{ for }  j=1,\ldots,n \,,
    \end{split} \end{equation}
  where \(\tau>0\). We will incorporate the parameter \(\tau\) into the notation \(\ladder[\uir{}{},\tau]{\pm}\) due to course once the subscript \(j\) will be discharged in the one-dimensional situation (or implicitly used in multiindex setup).
\end{defn}
For a unitary representation \(\uir{}{}\) the derived representation
is skew-symmetric: \((\rmd\uir{U}{})^*= -\rmd\uir{U}{}\) for all \(U\)
in the Weyl algebra \(\algebra{h}_n\). Thus the ladder operators are
adjoint of each other:
\begin{equation}
  \label{eq:ladder-are-adjoint}
  (\ladder[\uir{}{},j]{+})^*=    \ladder[\uir{}{},j]{-}\,.
\end{equation}
The main motivation for the above definition is that the Heisenberg
commutator relations
\([\rmd \uir{X}{\myhbar} , \rmd \uir{Y}{\myhbar}]=\rmi \myh I\) imply
a simple commutator of ladder operators:
\begin{equation}
  \label{eq:ladder-commutator}
  [\ladder[\uir{}{},l]{-},\ladder[\uir{}{},m]{+}]=
\delta_{lm}   \sign(\myh) I \qquad \text{ for any } \tau>0 \text{ and } 1\leq
  l,m\leq n\,
\end{equation}
with the dimensionless right hand side. For \(\tau\) of the dimensionality
\(\text{(mass)}\times\text{(frequency)}\) the agreement from
Rem.~\ref{re:physical-units} is satisfied
in~\eqref{eq:ladder-heisenberg} and we do not meet meaningless
fractional powers of physical units.
Furthermore, the ladder
operators are dimensionless and it is also useful to introduce the
respective dimensionless complex variable
\begin{equation}
  \label{eq:z_tau-complex-defn}
  z_\tau=\sqrt{\frac{\modulus{\myh}}{2\tau }}(x+\rmi \tau y) \in
\Space{C}{n}.
\end{equation}
Accordingly our standard convention we will simply denote \(z_\tau\) by \(z\) unless its dependence on \(tau\) needs to be indicated. 

As usual with the Heisenberg group, everything
essential already happens for \(n=1\) and cases \(n>1\) have only minor technical distinctions. To minimise those we will use the
convenient  notation:
\begin{align*}
  z\cdot w &= z_1w_1+z_2w_2+\ldots+z_n w_n\,,\\
  z\cdot \ladder[\uir{}{}]{\pm} &=   z_1 \ladder[\uir{}{},1]{\pm}
                                  +z_2 \ladder[\uir{}{},2]{\pm}+\ldots
                                  +z_n \ladder[\uir{}{},n]{\pm}\,.
\end{align*}

\begin{example}
  Without loss of generality let us assume that \(n=1\) and \(\myh>0\)
  for expressions in \S~\ref{sec:induc-repr}. Recall
  from~\eqref{eq:left-action-derived}--\eqref{eq:right-action-derived},
  that the left \(\Lambda_\myhbar\) and the right \(R_\myhbar\)
  pulled actions have the opposite signs of the derived
  representation for \(S\in\algebra{h}_1\), thus creation/annihilation r\^oles of ladder operators for the left and right pulled actions  are opposite.  In terms
  of the above complex variable \(z\) and the respective derivatives:\footnote{  Among two possible notations \(\overline{\partial}_z\) and \(\partial_{\overline{z}}\) we prefer the former for a better visibility of complex conjugation.}
  \begin{equation}
    \label{eq:complex-derivatives}
    \partial_z=\frac{1}{\sqrt{2\modulus{\myh}\tau
      }}(\tau \partial_x-\rmi \partial_y)\,, \qquad
        \overline{\partial}_{z}=\frac{1}{\sqrt{2\modulus{\myh}\tau
      }}(\tau \partial_x+\rmi \partial_y)\,.\\
  \end{equation}
 We rewrite
  \eqref{eq:left-action-derived}--\eqref{eq:schrodinger-derived} as:
  \begin{align}
    \label{eq:ladder-left-explicit}
    \ladder[\Lambda]{+}&= zI-\overline{\partial}_{z}
                         \,,
                       &
    \ladder[\Lambda]{-}&=  \overline{z}I+\partial_{z}\,; \\
    \label{eq:ladder-right-explicit}
    \ladder[R]{+}&= zI+\overline{\partial}_{z}\,,
                   &\ladder[R]{-}&=  \overline{z}I-{\partial}_{z}\,.
  \end{align}
  Also for the Schr\"odinger representation ladder operators are:
  \begin{equation}
    \label{eq:ladder-schrodinger-explicit}
    \ladder[\uir{}{}]{\pm}= \frac{1}{\sqrt{2\modulus{\myh}\tau}}
    ( \myh \tau  tI \mp \partial_t).
  \end{equation}
  These expressions for ladder operators look
  rather similar. Interestingly, the peeled actions
  \(\widetilde{\Lambda}_\myhbar =\Phi_\tau^{-1} \circ \Lambda_\myhbar \circ
  \Phi_\tau\) and
  \(\widetilde{R}_\myhbar =\Phi_\tau^{-1} \circ R_\myhbar \circ \Phi_\tau\) from
  Cor.~\ref{co:peeling-Gaussian} highlight their different structures:
  \begin{align}
    \label{eq:ladder-left-right-peeled}
    \ladder[\widetilde{\Lambda}]{+}&=2zI-\overline{\partial}_{z}
                         \,,
                       &
    \ladder[\widetilde{\Lambda}]{-}&= \partial_{z}\,; &
    \ladder[\widetilde{R}]{+}&= \overline{\partial}_{z}\,,
                   &\ladder[\widetilde{R}]{-}&= 2\overline{z}I-\partial_{z}\,;
  \end{align}  
  Furthermore, the restriction of \(\ladder[\widetilde{\Lambda}]{+}\) to
  the irreducible subspace of functions annihilated by
  \(\ladder[\widetilde{R}]{+}=\overline{\partial}_{z}\) is 
  the operator of multiplication \(2zI\).
\end{example}
Consequently, for
\((0,x,y)\in\Space{H}{n}\) and the respective
\(z=\sqrt{\frac{\modulus{\myh}}{2\tau }}(x+\rmi \tau y) \in
\Space{C}{n}\) we can express the representation \(\uir{}{}\) as
exponentiation of ladder operators:
\begin{equation}
  \label{eq:exponentiation-ladder-operators}
  \begin{split}
    \uir{}{\myhbar}(0,x,y)&=\exp(x\cdot \rmd \uir{X}{}+ y\cdot
    \rmd \uir{Y}{})\\
    &=\exp(\overline{z}\cdot \ladder[\uir{}{}]{+}- {z}\cdot\ladder[\uir{}{}]{-})
   \eqqcolon \uir{}{\myhbar}(z,\overline{z})\,.
  \end{split}
\end{equation}
Here \(\uir{}{\myhbar}(z,\overline{z})\) is known\footnote{We shall point out that our formula~\eqref{eq:exponentiation-ladder-operators} swaps positions of \(z\) and \(\overline{z}\) in comparison to the established physics sources.} as the \emph{displacement
  operator} in quantum optics~\citelist{\cite{Glauber63a} \cite{GerryKnight05a}*{\S~3.2}} and the
\emph{Weyl operator} in mathematical physics. Its explicit expression in
complex coordinates is:
\begin{align}
  \label{eq:left-action-pulled-complex}
  [\Lambda_\myhbar(s,z)f] (z',\overline{z}') &= \rme^{
                                       2\pi\rmi \myhbar s-(z\cdot \overline{z}'-\overline{z}\cdot z')/2}\, f(z'-z,\overline{z}'-\overline{z})
                                        , \\
  \label{eq:right-action-pulled-complex}
  [R_\myhbar(s,z)f] (z',\overline{z}') &= \rme^{
                                      -2\pi \rmi \myhbar s-(z\cdot \overline{z}'-\overline{z}\cdot z')/2}\, f(z'+z,\overline{z}'+\overline{z})
                                       .
\end{align}

The composition formula for displacement operators is
fully determined by the commutator~\eqref{eq:ladder-commutator} and is representation-in\-depen\-dent, i.e. does not contain \(\myhbar\):
\begin{align}
  \label{eq:ladder-exp-composition}
  \lefteqn{\exp(\overline{z}_1\cdot \ladder[\uir{}{}]{+}-{z}_1\cdot\ladder[\uir{}{}]{-})
  \exp(\overline{z}_2\cdot
  \ladder[\uir{}{}]{+}-{z}_2\cdot\ladder[\uir{}{}]{-})}
  & \\
  \nonumber 
  &\qquad =
    \exp\left(\half \Im({z}_1\cdot \overline{z}_2)\right)
    \exp \left((\overline{z}_1+\overline{z}_2)\cdot \ladder[\uir{}{}]{+}-({z}_1+{z}_2)\cdot\ladder[\uir{}{}]{-}\right)\,.
\end{align}
A helpful technique~\cite{KermackMcCrea31a} is the separation of the ladder operators in the regular representation by the \emph{Kermack--McCrae identity}\footnote{See~\cite{Coutinho10a} for an interesting discussion of the Kermack--McCrae's papers as another example of lost opportunities.}:
\begin{align}
  \nonumber 
    R_{\myhbar}(0,z)
    &=\exp(\overline{z}\cdot \ladder[R]{+}-{z}\cdot\ladder[R]{-})\\
  \label{eq:right-shift-ladder-form}
    &=\Phi(z,\overline{z})^{-1}\exp(\overline{z} \cdot \ladder[R]{+}) \exp(-{z} \cdot
    \ladder[R]{-})\\
  \label{eq:right-shift-ladder-form1}
    &=\Phi(z,\overline{z})\exp(-{z} \cdot
    \ladder[R]{-})\exp(\overline{z} \cdot \ladder[R]{+}) \,.
\end{align}
Here
\begin{equation}
  \label{eq:complex-gaussian}
  \Phi(z,\overline{z})  \coloneqq  \rme^{-\modulus{z}^2/2} =
  \rme^{-\pi\myhbar/(2\tau)(x^2+\tau^2y^2)}
\end{equation}
is a complexified form of \(\Phi(x,y)\)~\eqref{eq:FSB-Gaussian}. In
the same fashion we have:
\begin{align}
  \label{eq:left-shift-ladder-form}
  \Lambda _{\myhbar}(z,\overline{z})&=\Phi(z,\overline{z})
                                 \exp(-\overline{z} \cdot \ladder[\Lambda ]{+})
                                    \exp({z} \cdot \ladder[\Lambda ]{-})\\
  \label{eq:left-shift-ladder-form1}
                               &=\Phi(z,\overline{z})^{-1}
                                 \exp({z} \cdot \ladder[\Lambda ]{-})
                                 \exp(-\overline{z} \cdot \ladder[\Lambda ]{+})\,.
\end{align}
Recall, that the difference between \eqref{eq:right-shift-ladder-form}--\eqref{eq:right-shift-ladder-form1} and \eqref{eq:left-shift-ladder-form}--\eqref{eq:left-shift-ladder-form1} is due to opposite signs of \(\rmd R_\myhbar^S \) and \(\rmd \Lambda_\myhbar^S \) echoing in the commutator~\eqref{eq:ladder-commutator}.

  An analytic meaning of operators \(\ladder[\Lambda]{\pm}\)  is the
  action by shifts
  \begin{equation}
    \label{eq:ortho-normal-basis}
    \ladder[\Lambda]{+}:\, \Phi_m
    \mapsto  \sqrt{m+1} \Phi_{m+ 1}\,, \qquad
    \ladder[\Lambda]{-}:\, \Phi_m
    \mapsto \sqrt{m} \Phi_{m- 1}\,.
  \end{equation}
  on the orthonormal basis
  \(\Phi_m=(\pi^m m!)^{-1/2} (\ladder[\Lambda]{+})^n \Phi\) within an
  irreducible component \(\FSpace[\varphi]{F}{}\) of the representation
  \(\Lambda_\myhbar\). On the other hand, operators
  \(\ladder[R]{\pm}\) acts in the similar fashion by ``shifting'' different irreducible
  components of \(\Lambda_\myhbar\) one into another. Namely:
  \begin{equation}
    \label{eq:irreducible-components-shift}
    \ladder[R]{\pm}:\, \FSpace[\varphi]{F}{m} \rightarrow
    \FSpace[\varphi]{F}{m\mp1}
    \quad \text{ where } \FSpace[\varphi]{F}{m}
    = (\ladder[R]{-})^m \FSpace[\varphi]{F}{}
    \ (m=0,1,\ldots)\,.
  \end{equation}
  For a  consistence we  set \(\FSpace[\varphi]{F}{-1}=\{0\}\)  in~\eqref{eq:irreducible-components-shift}. 
  From commutativity of the left and right pulled actions  \(\FSpace[\varphi]{F}{m}\) are \(\Lambda_\myhbar\)-invariant
  irreducible components of the orthogonal decomposition
  \(\FSpace{L}{2}(\Space{C}{n})= \operp_{m=0}^\infty
  \FSpace[\varphi]{F}{m} \), cf.~\cite{Vasilevski99b}. 
  Also these spaces
  \begin{displaymath}
    \FSpace[\varphi]{F}{m}
    \coloneqq  \FSpace[\varphi_m]{F}{}
  \end{displaymath}
  are image spaces~\eqref{eq:image-space-coav-tr} of the covariant transform for the Hermite functions \(\varphi_m=(\ladder[\myhbar]{+})^m \varphi\) as the respective ground vectors. Obviously, spaces \( \FSpace[\varphi]{F}{m}= (\ladder[R]{-})^m \FSpace[\varphi]{F}{}\) are annihilated by powers of the right ladder operator \(    (\ladder[R]{+})^n\) with \(n\geq m\). Spaces \(\FSpace[\varphi]{F}{m}\) were named true poly-analytic in~\cite{Vasilevski99b}. The concept was recently revised in~\cite{TurbinerVasilevski21a} from the representation theory viewpoint, where the authors employed the semidirect product of the Heisenberg group and \(\mathrm{SL}_2(\Space{R}{})\), known as the Schr\"odinger~\cite{Folland89}*{\S~1.2} or Jacobi~\cite{Berndt07a}*{\S~8.5} group. Although the latter and its various subgroups~\cite{AlmalkiKisil19a,AlmalkiKisil18a,Kisil21c} are very interesting and important objects, they may be excessive for the discussed poly-analytic function decomposition, which is completely manageable by the two-sided action of the Heisenberg group alone.
\begin{lem}
  The collection of functions \(\Phi_{jk}= (\ladder[\Lambda]{+})^j(\ladder[R]{-})^k \Phi\) forms an orthonormal basis of \(\FSpace{L}{2}(\Space{R}{2n})\) and, thereafter, an arbitrary \(v\in\FSpace{L}{2}(\Space{R}{2n})\) admits the presentation 
  \begin{equation}
    \label{eq:L2-decomposition-ladders}
    v = \sum_{j,k=0}^\infty v_{jk}(\ladder[\Lambda]{+})^j(\ladder[R]{-})^k \Phi, \quad \text{ where } v_{jk} = \scalar{v}{\Phi_{jk}}= \scalar{v}{(\ladder[\Lambda]{+})^j(\ladder[R]{-})^k \Phi}\in\Space{C}{}.
  \end{equation}
\end{lem}
See Fig.~\ref{fig:2D-lattice} for an illustration. Note, that the order of ladder operators in~\eqref{eq:L2-decomposition-ladders} is not important since  they commute.
This commutativity is also behind the following result.
\begin{figure}[htbp]
      \centering
      \(  \xymatrix@=4em@M=1em{
        & 
        \,0\,  &  
        \,0\,  & 
        \,0\,   & 
        \\
        0\, & 
        \,\Phi_{00}\,  \ar@<.4ex>[l]^-{\ladder[R]{+}}\ar@<.4ex>[r]^{\ladder[R]{-}}
        \ar@<.4ex>[u]^{\ladder[\Lambda]{-}} \ar@<.4ex>[d]^{\ladder[\Lambda]{+}} &  
        \,\Phi_{01}\, \ar@<.4ex>[l]^{\ladder[R]{+}} \ar@<.4ex>[r]^{\ladder[R]{-}}
        \ar@<.4ex>[u]^{\ladder[\Lambda]{-}} \ar@<.4ex>[d]^{\ladder[\Lambda]{+}} & 
        \,\Phi_{02}\,\ar@<.4ex>[l]^{\ladder[R]{+}}  \ar@<.4ex>[r]^-{\ladder[R]{-}}
        \ar@<.4ex>[u]^{\ladder[\Lambda]{-}} \ar@<.4ex>[d]^{\ladder[\Lambda]{+}}    & 
        \,\ldots\ar@<.4ex>[l]^-{\ladder[R]{+}}\\
        0\, & 
        \,\Phi_{10}\,  \ar@<.4ex>[l]^-{\ladder[R]{+}}\ar@<.4ex>[r]^{\ladder[R]{-}}
        \ar@<.4ex>[u]^{\ladder[\Lambda]{-}} \ar@<.4ex>[d]^{\ladder[\Lambda]{+}} &  
        \,\Phi_{11}\, \ar@<.4ex>[l]^{\ladder[R]{+}} \ar@<.4ex>[r]^{\ladder[R]{-}}
        \ar@<.4ex>[u]^{\ladder[\Lambda]{-}} \ar@<.4ex>[d]^{\ladder[\Lambda]{+}}& 
        \,\Phi_{12}\,\ar@<.4ex>[l]^{\ladder[R]{+}}  \ar@<.4ex>[r]^-{\ladder[R]{-}}
        \ar@<.4ex>[u]^{\ladder[\Lambda]{-}} \ar@<.4ex>[d]^{\ladder[\Lambda]{+}}    & 
        \,\ldots\ar@<.4ex>[l]^-{\ladder[R]{+}}\\
        0\, & 
        \,\Phi_{20}\,  \ar@<.4ex>[l]^-{\ladder[R]{+}}\ar@<.4ex>[r]^{\ladder[R]{-}}
        \ar@<.4ex>[u]^{\ladder[\Lambda]{-}} \ar@<.4ex>[d]^{\ladder[\Lambda]{+}} &  
        \,\Phi_{21}\, \ar@<.4ex>[l]^{\ladder[R]{+}} \ar@<.4ex>[r]^{\ladder[R]{-}} 
        \ar@<.4ex>[u]^{\ladder[\Lambda]{-}} \ar@<.4ex>[d]^{\ladder[\Lambda]{+}}& 
        \,\Phi_{22}\,\ar@<.4ex>[l]^{\ladder[R]{+}}  \ar@<.4ex>[r]^-{\ladder[R]{-}}
        \ar@<.4ex>[u]^{\ladder[\Lambda]{-}} \ar@<.4ex>[d]^{\ladder[\Lambda]{+}}    & 
        \,\ldots\ar@<.4ex>[l]^-{\ladder[R]{+}}\\
        & 
        \,\ldots\, \ar@<.4ex>[u]^{\ladder[\Lambda]{-}} &  
        \,\ldots\, \ar@<.4ex>[u]^{\ladder[\Lambda]{-}} & 
        \,\ldots\,  \ar@<.4ex>[u]^{\ladder[\Lambda]{-}}  &
        \save "2,2"."5,2"*[F--]\frm{}\restore
        \save "2,3"."5,3"*[F.]\frm{}\restore
       \save "2,4"."5,4"*[F.]\frm{}\restore}
      \)
      \caption[The action of hyperbolic ladder operators on a 2D
      lattice of eigenspaces]{The action of the ladder operators on
        the lattice of orthogonal vectors. Boxed columns form orthonormal bases of \(\Lambda_\myhbar\)-invariant subspaces \(\FSpace[\varphi]{F}{m}\). Left ladder operators shift elements within a basis. Right ladder operators map elements of one basis to the respective elements of an adjacent basis.}  
      \label{fig:2D-lattice}
    \end{figure}
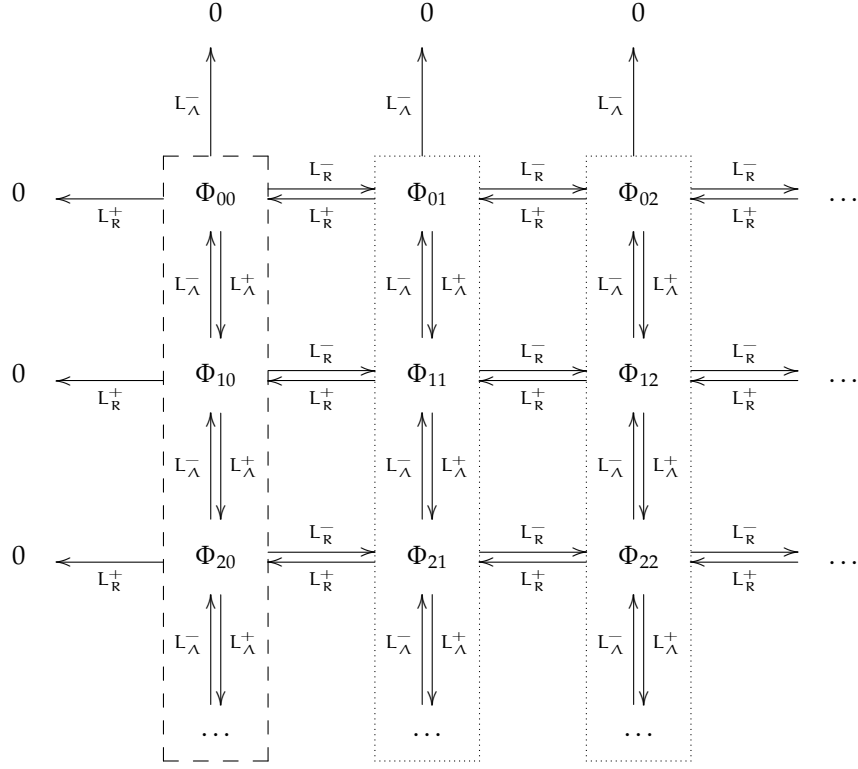
\begin{prop}
  Let a bounded operator \(A: \FSpace{L}{2}(\Space{C}{n}) \rightarrow \FSpace{L}{2}(\Space{C}{n})\) commute with all right ladder operators, i.e.  \([\ladder[R]{\pm},A]=0\). Then  all spaces \(\FSpace[\varphi]{F}{m}\) are \(A\)-invariant and \(A\) is a left relative convolution on \(\FSpace{L}{2}(\Space{C}{n}) \).
\end{prop}
\begin{proof}
  By induction we can show that \([\ladder[R]{\pm},A]=0\) is equivalent to \([(\ladder[R]{\pm})^m,A]=0\) for all natural numbers \(m\).  Assume towards a contradiction that for some \(m\) the space \(\FSpace[\varphi]{F}{m}\) is not invariant, that is there are exist \(v_m \in \FSpace[\varphi]{F}{m}\) and \(v_k \in \FSpace[\varphi]{F}{k}\) for some \(k\neq m \) such that \(\scalar{A v_m}{v_k}\neq 0\). Then:
  \begin{enumerate}
  \item If \(k < m\), let \(v \in \FSpace[\varphi]{F}{m-k}\) be such that \((\ladder[R]{-})^k v = v_m\), Then:
    \begin{align*}
      0 & \neq \scalar{A v_m}{v_k}\\
      & = \scalar{A (\ladder[R]{-})^k v}{v_k} \\
      & = \scalar{(\ladder[R]{-})^k A v}{v_k} \\
      & = \scalar{A v}{ (\ladder[R]{+})^k v_k} \\
      & = \scalar{A v}{ (\ladder[R]{+})^k v_k} \\
      & = \scalar{A v}{0}\\
      & =0\,.
    \end{align*}
  \item If \(k > m\), let \(v \in \FSpace[\varphi]{F}{k-m}\) be such that \((\ladder[R]{-})^m v = v_k\), Then:
    \begin{align*}
      0 & \neq \scalar{A v_m}{v_k}\\
      & = \scalar{A v_m}{ (\ladder[R]{-})^m v} \\
      & = \scalar{(\ladder[R]{+})^m A v_m}{v_k} \\
      & = \scalar{A (\ladder[R]{+})^m v_m}{ v_k} \\
      & = \scalar{0}{ v_k} \\
      & =0\,.
    \end{align*}
  \end{enumerate}
  We had obtained the contradiction.

  To show that \(A\) is a relative convolution, the Schur lemma implies that the restriction of \(A\) to the irreducible component \(\FSpace[\varphi]{F}{0}\) belongs to the \(C^*\)-algebra spanned by the left action of \(\Lambda_{\myhbar}\) of \(\Space{H}{n}\) restricted to \(\FSpace[\varphi]{F}{0}\). Thus the restriction of \(A\) is equal to a relative convolution on \(\FSpace[\varphi]{F}{0}\) with some kernel \(k(g)\). By commutativity with ladder operators, restrictions of \(A\) to any \(\FSpace[\varphi]{F}{m}\), \(m\in\Space{N}{}\) is again the relative convolution with the same kernel \(k\). The same is true for \(\FSpace{L}{2}(\Space{C}{n})\)---the direct sum of all \(\FSpace[\varphi]{F}{m}\), \(m\in\Space{N}{}\).
\end{proof}
The vanishing commutators \([\ladder[R]{\pm},A]=0\) is a sufficient but is not necessary condition for \(A\)-invariance of subspaces. A counterexample is a diagonal-type operator \(A: \FSpace{L}{2}(\Space{C}{n}) \rightarrow \FSpace{L}{2}(\Space{C}{n})\) such that restriction
of \(A\) to \(\FSpace[\varphi]{F}{m}\) is the identity operator times \(x_m\) for a bounded non-constant sequence \(x_m\).  We state the following simple result:
\begin{lem}
  The space \( \FSpace[\varphi]{F}{}\) is invariant under an operator \(A\) if and only if  \( \FSpace[\varphi]{F}{}\) is in the kernel of the commutator \([\ladder[R]{+}, A]\).
\end{lem}
\begin{proof}
  Let  \(Af\in \FSpace[\varphi]{F}{}\) for any \(f\in \FSpace[\varphi]{F}{}\). Then \(\ladder[R]{+} f =0\) implies \(\ladder[R]{+} Af =0\) and \([\ladder[R]{+}, A] f =0 \). Conversely, if \([\ladder[R]{+}, A] f =0 \) for some \(f\in \FSpace[\varphi]{F}{}\), then \(\ladder[R]{+} A f =0 \), that is  \(f\in \FSpace[\varphi]{F}{}\) is an invariant subspace of \(A\).
\end{proof}

To conclude this subsection we use complex variables to state the following straightforward properties of the mixed Gaussian~\(\Phi_{\tau\varsigma}\)~\eqref{eq:Phi-tau-sigma}:
\begin{lem}
\label{le:intertwining-kernel-properties}
\begin{enumerate}
  \item For \(\tau=\varsigma\) the mixed Gaussian coincides with the Gaussian: \(\Phi_{\tau\tau}=\Phi_\tau\).  
  \item Complex conjugations swaps parameters \(\tau\) and \(\varsigma\), that is \(\overline{\Phi}_{\tau\varsigma}=\Phi_{\varsigma\tau}\).
  \item \(\Phi_{\tau\varsigma}\) is fixed by the symplectic Fourier transform.\(\Phi_{\tau\varsigma} = \wideparen{\Phi}_{\tau\varsigma}\).
  \item Complex form of the mixed Gaussian in notations~\eqref{eq:z_tau-complex-defn} is:
    \begin{displaymath}
    \textstyle  \Phi_{\tau\varsigma}=\sqrt[4]{\frac{\tau \varsigma }{(\tau+\varsigma)^2}} \exp\left( -    \frac{\sqrt{\tau \varsigma}}{\tau+\sigma} 
    \overline{z}_\varsigma z_\tau\right).
    \end{displaymath}
  \item
    \label{it:mixed-vacuum-right}
    \(\Phi_{\tau\varsigma}\) is annihilated by the right ladder operator \(\ladder[R,\tau]{+}\) with squeeze \(\tau\).
  \item
    \label{it:mixed-vacuum-left}
    \(\Phi_{\tau\varsigma}\) is annihilated by the left ladder operator \(\ladder[\Lambda,\varsigma ]{-}\)  with squeeze \(\varsigma\).
  \end{enumerate}
\end{lem}

\subsection{Twisted convolution and symplectic Fourier transform}
\label{sec:twist-conv-sympl}

Exponentiation of ladder operators in the peeled
form~\eqref{eq:ladder-left-right-peeled} has an explicit meaning. Indeed, for a polynomial \(p(z,\overline{z})\) 
an algebraic manipulation with the Taylor expansion of the exponent
produces:
\begin{align}
  \label{eq:exp-ladder-left-right-operator}
                                                                            [\exp({w} \cdot \ladder[\widetilde{\Lambda} ]{-}) p](z,\overline{z})&=
                                                                                                                                             p(z+w,\overline{z})\,,
  &
  [\exp(\overline{w} \cdot \ladder[\widetilde{R}]{+}) p](z,\overline{z})&=
                                                                p(z,\overline{z}+\overline{w})\,.
\end{align}
These formulae can be extended by continuity to any space of
functions such that polynomials in \(z\) and \(\overline{z}\) form a
dense subspace.
Due to properties from Lem.~\ref{le:intertwining-kernel-properties}.\ref{it:mixed-vacuum-right}--\ref{it:mixed-vacuum-left} the respective exponents act on \(\Phi_{\tau\varsigma}\) trivially:
\begin{align}
  \label{eq:exp-ladder-on-gaussian-is-trivial}
  \exp(w \cdot \ladder[{\Lambda,\varsigma} ]{-}) \Phi_{\tau\varsigma}&= \Phi_{\tau\varsigma}\,,&
                                                                               \exp(\overline{w} \cdot \ladder[R,\tau]{+})\Phi_{\tau\varsigma}&= \Phi_{\tau\varsigma}\,.
\end{align}
The symplectic Fourier
transform~\eqref{eq:symplect-Fourier-defn} in complex variables is:
\begin{align}
  \label{eq:symplect-Fourier-complex}
  \wideparen{\psi}(z,\overline{z})&= (\myhbar/2)^{n}\int_{\Space{C}{n}}
                               {\psi}(z',\overline{z}')\,
                               \exp ({\pi \rmi (\overline{z}' z - z'\overline{z})})\,\rmd z'\wedge\rmd \overline{z}'  \,.
\end{align}
Of course, the intriguing relation~\eqref{eq:symplect-Fourier-left},
i.e. \(  \wideparen{a} = a \twist \mathbf{1}
= \Lambda_\myhbar(a)  \mathbf{1}
= R_\myhbar(a) \mathbf{1}\), is equally true in complex coordinates
as it was in real ones. Therefore, for ladder operators we have the following intertwining properties:
\begin{equation}
  \label{eq:sympl-Fourier-intertwin-left-complex}
  \wideparen{\ } \circ \ladder[{\Lambda} ]{\pm}
  =  \ladder[{\Lambda} ]{\pm}  \circ  \wideparen{\ }
  \quad \text{and} \quad
  \wideparen{\ } \circ \ladder[R ]{\pm}
  =  - \ladder[R]{\pm}  \circ  \wideparen{\ }\,.
\end{equation}

Furthermore,
the intertwining properties~\eqref{eq:sympl-Fourier-intertwin-left}
and~\eqref{eq:sympl-Fourier-intertwin-right} implying the following
complexified  versions:
\begin{equation}
  \label{eq:fourier-intertw-shift-mult-compl}
  \wideparen{\ } \circ \map{S} (z,\overline{z})
  = \map{E} (z,\overline{z}) \circ  \wideparen{\ }\qquad
  \wideparen{\ } \circ \map{E} (z,\overline{z})
  = \map{S} (z,\overline{z}) \circ  \wideparen{\ }\,,
\end{equation}
where  complexified versions of spatial~\eqref{eq:shift-defn} and
frequency~\eqref{eq:exp-mult-defn} shifts are:
\begin{align}
  \label{eq:shift-compl-defn}
  \map{S} (z,\overline{z})&=\Lambda_\myhbar (0,\half z)\circ R_\myhbar
           (0,-\half z):& f(z',\overline{z}') &\mapsto f(z'-z,\overline{z}'-\overline{z})\,,
  \\
 \label{eq:modulation-complex--defn}
  \map{E} (z,\overline{z})&=\Lambda_\myhbar (0,\half z)\circ R_\myhbar
           (0,\half z):& f(z',\overline{z}') &\mapsto \rme^{-
                                            (z\overline{z}'-\overline{z}z')/2}\,f(z',\overline{z}')\,.
\end{align}
Note, that in contrast to~\eqref{eq:exp-mult-defn} the complexified
frequency shift does not explicitly depend on the Planck constant, the
dependence is incorporated in the complex variable \(z\).  The
infinitesimal version of~\eqref{eq:fourier-intertw-shift-mult-compl}
is:
\begin{equation}
  \label{eq:fourier-intertw-shift-mult-compl-infinites}
  \wideparen{\ } \circ \partial_z
  = -\half\overline{z}I \circ  \wideparen{\ }\,,\quad
  \wideparen{\ } \circ \overline{\partial}_z
  =  \half zI \circ  \wideparen{\ }\,, \quad
    \wideparen{\ } \circ \overline{z}I
  = -2\partial_z \circ  \wideparen{\ }\,,\quad
  \wideparen{\ } \circ zI
  =  2\overline{\partial}_z  \circ  \wideparen{\ }\,.
\end{equation}

\section{Co- and contra-variant transforms of operators}
\label{sec:repr-oper-pdo}

Based on the previous consideration we are ready to describe relations within the family of  operators: localisation, Toeplitz, pseudodifferential and some their further generalisations.

\subsection{Localisation and Toeplitz operators}
\label{sec:local-toepl-oper}

The fundamental identity~\eqref{eq:wave-trans-inverse} implies that       \(\contravar{\theta } \circ \covar{\theta} = I\) for a vector \(\theta \) with unit norm. This identity is a source of a more sophisticated object, known as \emph{localisation operator}~\cites{BoggiattoCorderoGrochenig04a,Coburn01a,CorderoGrochenig,Coburn19a,AbreuFaustino15a}:
\begin{equation}
  \label{eq:localization-operator}
  \oper{L}_\psi = \contravar{\theta } \circ \psi I \circ \covar{\theta},
\end{equation}
where \(\psi\) is a function on the phase space called \emph{symbol} or \emph{weight function}. In a weak form the operator can be written as:
\begin{equation}
  \label{eq:localization-operator-weak}
  \scalar{\oper{L}_\psi u}{v} = \int_{\Space{R}{2n}} \psi(x,y) \scalar{u}{\uir{}{}(x,y)\theta}
  \scalar{\uir{}{}(x,y) \theta }{v}\rmd x\, \rmd y,
\end{equation}
where \(\uir{}{}\) is a unitary representation on a Hilbert space \(\FSpace{H}{}\) and \(u\), \(v\) are arbitrary vectors in \(\FSpace{H}{}\). Of course, identities similar to~\eqref{eq:localization-operator}--\eqref{eq:localization-operator-weak} can be used to define localisation operators for representations of other groups as well, but we remain with the Heisenberg group within this paper.

For the special choice of \(\uir{}{}=\Lambda_\myhbar\) and \(\theta =\Phi_\tau\), the operator in~\eqref{eq:localization-operator} is more known as the \emph{Toeplitz operators}
\(\oper{T}_\psi\) on (pre-)FSB space. Equivalently, Toeplitz operators \(\oper{T}_\psi\) and \(\tilde{\oper{T}}_\psi\)  on pre-FSB  space \( \FSpace[\tau]{F}{}\) and FSB spaces \(\FSpace{\tilde{F}}{\tau}\) respectively are defined by:
\begin{equation}
  \label{eq:two-toeplitz-definition}
  \oper{T}_\psi=\oper{P}_\tau \psi I
  \quad \text{ and } \quad
  \tilde{\oper{T}}_\psi=\tilde{\oper{P}}_\tau\psi I
  \,,
\end{equation}
where
\(\oper{P}_\tau :\; \FSpace{L}{2}(\Space{C}{n}) \rightarrow
\FSpace[\tau]{F}{}\) and
\(\tilde{\oper{P}}_\tau:\; \FSpace{L}{2}(\Space{C}{n},
\rme^{-\modulus{z}^2}\rmd z \wedge \rmd \overline{z}) \rightarrow
\FSpace{\tilde{F}}{\tau}\) are corresponding orthoprojections. 

Our study of  Toeplitz operators is preceded by a simple useful observation on their relations on FSB and pre-FSB spaces. Recall, that a unitary equivalence of \(\FSpace[\tau]{F}{} \) and
\(\FSpace{\tilde{F}}{\tau}\) is provided by the peeling map
\(\oper{E}_{\Phi}: f \mapsto \Phi f\)~\eqref{eq:peeling} and
\(\oper{P}_\tau = \oper{E}_\Phi \circ \tilde{\oper{P}} \circ
\oper{E}_\Phi^{-1}\).  The peeling operators \(\oper{E}_{\Phi}\) obviously commutes
with an operator of multiplication \(\psi: f \mapsto \psi
f\). Thereafter,
\begin{lem}
  For a bounded function \(\psi\), the Toeplitz
  operators~\eqref{eq:two-toeplitz-definition} are unitary equivalent
  through the peeling:
  \begin{displaymath}
    \oper{T}_\psi = \oper{E}_\Phi \circ \tilde{\oper{T}}_\psi \circ \oper{E}_\Phi^{-1}.
  \end{displaymath}
\end{lem}
In other words, the main properties of Toeplitz operators in terms of their symbol \(\psi\) are identical in both spaces.
Furthermore, consideration of pre-FSB spaces open the door for the following generalisation.
\begin{defn}
  For given positive reals \(\tau\), \(\varsigma\) and a function \(\psi(x,y)\) on \(\Space{R}{2n}\) we define \emph{cross-Toeplitz operator} \(\oper{T}_{\psi}: \FSpace[\tau]{F}{} \rightarrow \FSpace[\varsigma]{F}{}\) by:
  \begin{equation}
    \label{eq:toeplitz-defn}
    \oper{T}_\psi f =\oper{P}_\varsigma (\psi f), \qquad \text { where }
    f \in \FSpace[\tau]{F}{}\,.
  \end{equation}
  The function \(\psi\) is called the \emph{symbol} of the operator \(\oper{T}_\psi \). If parameters \(\tau\), \(\varsigma\) need to be explicitly indicated we will use the longer notation \(\oper{T}_\psi^{ \tau \varsigma }\) for the operator.
\end{defn}
Obviously, the traditional Toeplitz operators correspond to the case \(\tau=\varsigma\). We can also consider \(\oper{P}_\varsigma \psi I\) as a map \(\FSpace{L}{2}(\Space{R}{2n}) \rightarrow \FSpace[\varsigma]{F}{}\) which serves as an umbrella for all cross-Toeplitz operators  \(\oper{T}_\psi^{ \tau \varsigma }\) with any \(\tau>0\).

A natural family of symbols \(\psi\) for cross-Toeplitz operators \(\oper{T}_\psi\)~\eqref{eq:toeplitz-defn} requires that the set of functions \(f\in \FSpace[\tau]{F}{}\), such that \(\psi f \in \FSpace{L}{2}(\Space{R}{2n})\), is dense in  \(\FSpace[\tau]{F}{}\). The further discussion will give alternative expressions of cross-Toeplitz operators which is applicable for a wider set of symbols.

\begin{rem}
  Our cross-Toeplitz operators belongs to a more general class of operators prompted by the identity~\eqref{eq:wave-trans-inverse}. Let representations \(\uir{}{j}\) act in Hilbert spaces \(\FSpace[j]{H}{}\) with \(j=1,2\) and  \(\theta_j\in\FSpace[j]{H}{}\) be two fixed admissible vectors. Then for a bounded function \(\psi\) a generalised Anti-Wick operator \(\FSpace[1]{H}{}\rightarrow \FSpace[2]{H}{}\) is defined by~\cite{BoggiattoCorderoGrochenig04a}:
  \begin{equation}
    \label{eq:localization-operator-general}
    \oper{L}_\psi = \contravar{\theta_2 } \circ \psi I \circ \covar{\theta_1},
  \end{equation}
  That is, our cross-Toeplitz operator corresponds to the representations of the Heisenberg group and \(\theta_1\), \(\theta_2\) being the Gaussians (with possibly different squeeze parameters). Other groups and analysing/reconstructing vectors are out of our current scope. 
\end{rem}

\subsection{Toeplitz operators and PDO with heat flow}
\label{sec:toepl-oper-cald}

If a representation \(\uir{}{}\) is induced by a character of a
subgroup \(H\subset G\), then integration in~\eqref{eq:integrated-rep}
can be reduced to the homogeneous space \(X=G/H\) with the essentially
same set of resulting operators.  More specifically, for the natural
projection \(\map{p}: G \rightarrow X\) we fix a Borel section
\(\map{s}: X \rightarrow G\)~\cite[\S~13.2]{Kirillov76}, which is a
right inverse to \(\map{p}\).  Recall, we define an operator of
\emph{relative convolution}%
\index{relative convolution}%
\index{convolution!relative} on \(V\)~\cite{Kisil94e,Kisil13a},
cf.~\eqref{eq:integrated-rep}:
\begin{equation}
  \label{eq:relative-conv}
  \uir{}{}(k)=\int_{X} k(x)\,\uir{}{}(\map{s}(x))\,\rmd x 
\end{equation}
with a kernel \(k\) being a function on \(X=G/H\).  For example, if
\(G\) is the Heisenberg group and \(\uir{}{}\) is its Schr\"odinger
representation~\eqref{eq:schroedinger-rep}, then
\(\uir{}{}(\wideparen{a})\)~\eqref{eq:relative-conv}, where
\(\wideparen{a}\) is the symplectic Fourier transform of
\(a\)~\eqref{eq:symplect-Fourier-defn}, is a \emph{pseudodifferential operator} (PDO) \(a(D,X)\) with the
symbol \(a\)~\citelist{\cite{Howe80b} \cite{Folland89}*{\S~2.1}
  \cite{Kisil13a}}:
\begin{align}
  \nonumber
                   [a(D,X) f](t)& = \int_{\Space{R}{2n}} a(\xi, \textstyle \frac{1}{2}(t+r))\,
                                  \rme^{\pi\rmi \myhbar \xi (t-r)} f(r)\,\rmd r\,\rmd \xi
\end{align}
with the Weyl (symmetrized) symbol \(a(\xi,x)\).
For the sake of
completeness we provide a calculation parallel to
\cite{Folland89}*{\S~2.1} in our notations 
\begin{align*}
  \label{eq:pdo-def}
  [Af](t)&=[\uir{}{}(\wideparen{a})f](t)
           =\int_{\Space{R}{2n}} \wideparen{a}(x,y)\,\uir{}{\myhbar}(x,y)f(t)\,\rmd x\,\rmd y\\
         &= \int_{\Space{R}{2n}} \wideparen{a}(x,y)\,\rme^{\pi \rmi \myhbar x(2t-y)}\,
           f(t-y)\,\rmd x\,\rmd y\nonumber \\
    &= \myhbar^{n}\int_{\Space{R}{2n}}\int_{\Space{R}{2n}} a(x',y')\,
        \rme^{\pi\rmi\myhbar (x'y-y'x)} \,
      \rme^{\pi  \rmi \myhbar  x(2t-y)}\,
    f(t-y)\,\rmd x'\,\rmd y'\,\rmd x\,\rmd y\nonumber \\
    &= \int_{\Space{R}{2n}}\int_{\Space{R}{2n}} a(x',y')\,
        \rme^{\pi\rmi\myhbar x'y}
        \delta(\myhbar(\half (y+y')-t))\,
    f(t-y)\,\rmd x'\,\rmd y'\,\rmd y\nonumber \\
    &= \myhbar^{-n}\int_{\Space{R}{2n}} \textstyle  a(x',  t-\frac{1}{2}y)\,
    \rme^{\pi\myhbar\rmi x'y}
      \,    f(t-y)\,\rmd x'\,\rmd y\nonumber \\
    &= \myhbar^{-n} \int_{\Space{R}{2n}} \textstyle a(x', \frac{1}{2}(t+r))\,
    \rme^{\pi\myhbar\rmi x'(t-r)}
      \,    f(r)\,\rmd x'\,\rmd r.\nonumber 
\end{align*}

Now we revise the method used in~\cite[\S~3.1]{Howe80b} to
prove the Calder\'on--Vaillancourt estimations. It was described as ``rather
magical'' in~\cite[\S~2.5]{Folland89}. Use of the covariant
transform dispels the mystery without undermining the power of the
method. Relevantly for the present topic, the
demonstration~\cite[\S~3.1]{Howe80b} implicitly expresses a PDO as a
Toeplitz operator---the result which commonly attributed to a later
work of Guillemin~\cite{Guillemin84}*{(8.20)}, see also~\cite{Folland89}*{\S~2.7}.

We start from the following lemma, which has a transparent proof in
terms of covariant transform, cf. earlier presentations in~\cite[\S~3.1]{Howe80b}
and~\cite[(2.75)]{Folland89} in the case of \(\theta =\phi\).

\begin{lem}
  Let \(\uir{}{}\) be an irreducible square integrable representation
  a Lie group \(G\) in \(V\) and mother wavelets \(\phi\) and  \(\theta\) are
  admissible. Then
  \begin{equation}
    \label{eq:left-and-right-is-phi-rho}
    \Phi(g) \uir{}{}(g)
    =  \contravar{\theta} \circ (\Lambda \otimes R)(g,g)
    \circ\covar{\phi}
    \qquad \text{for all } g\in G
  \end{equation}
  and \(\Phi(g) = \scalar{\theta}{\uir{}{}(g)\phi}\).
\end{lem}
\begin{proof}
  We know from~\eqref{eq:wave-trans-inverse} that \(\contravar{\theta}
  \circ \covar{{\uir{}{}(g)\phi}} =
  \scalar{\theta}{\uir{}{}(g)\phi} I\) on \(V\), thus:
  \begin{displaymath}
    \contravar{\theta} \circ \covar{{\uir{}{}(g)\phi}} \circ \uir{}{}(g) =
    \scalar{\theta}{\uir{}{}(g)\phi} \uir{}{}(g) =\Phi (g) \uir{}{}(g)
    . 
  \end{displaymath}
  On the other hand, the intertwining
  properties~\eqref{eq:left-shift-itertwine}
  and~\eqref{eq:wave-intertwines-right} of the wavelet transform 
  imply:
  \begin{displaymath}
    \contravar{\theta} \circ \covar{{\uir{}{}(g)\phi}} \circ \uir{}{}(g) =
    \contravar{\theta} \circ (\Lambda \otimes R)(g,g)  \circ\covar{\phi}.
  \end{displaymath}
  A combination of the above two identities yields~\eqref{eq:left-and-right-is-phi-rho}.
\end{proof}
We will use a  specialisation of~\eqref{eq:left-and-right-is-phi-rho} to the
Schr\"odinger representation \(\uir{}{\myhbar}\) and the Gaussian
\(\phi\).
\begin{prop}
  \label{pr:toeplitz-as-PDO-Guillemin}
  The cross-Toeplitz operator
  \(\oper{P}_\varsigma \psi I : \FSpace[\tau]{F}{} \rightarrow \FSpace[\varsigma]{F}{}\) is unitary
  equivalent to PDO \(a(D,X): \Space{R}{n} \rightarrow \Space{R}{n}\)
  with the symbol:
  \begin{equation}
    \label{eq:Toeplitz-PDO-diffused}
    a_\psi(x,y)=
    {\myhbar^{2n}}
      \int_{\Space{R}{2n}} {\psi}(x',y')\,  {\Phi}_{\tau\varsigma}(x-2x',y-2y') \,\rmd x'\,\rmd y'\,,
    \end{equation}
    where \({\Phi}_{\tau\varsigma}\) is the intertwining kernel~\eqref{eq:Phi-tau-sigma}.
\end{prop}
\begin{proof}
  As we observed in~\eqref{eq:mult-by-left-right} that \((\Lambda \otimes R)(g,g)\)
  is symplectic phase shift \(\map{E}_\myhbar(x,y)\):
  \begin{equation}
    \label{eq:left-times-right-rep}
    [\Lambda_\myhbar(s,x,y)R_\myhbar(s,x,y)f] (x',y')
    = \map{E}_\myhbar(2x,2y) f(x', y') = \rme^{2\pi \rmi \myhbar (xy'-yx')} f(x', y').
  \end{equation}
  Thus, integrating the identity~\eqref{eq:left-and-right-is-phi-rho} with the
  function \((2\myhbar)^n \wideparen{\psi}(2x,2y)\) over \(G/H\) we obtain:
  \begin{equation}
    \label{eq:pdo-multiplication-intertwine}
    {a}_\psi(D,X)    = \contravar{\theta }  \circ \psi I  \circ\covar{\phi}
  \end{equation}
  where
  \(\wideparen{a}_\psi(x,y) = (2\myhbar)^n
  \wideparen{\psi}(2x,2y) \Phi_{\tau\varsigma}(x,y)\). 
  Therefore, the standard
  manipulations based on intertwining
  properties~\eqref{eq:fourier-intertw-shift-mult} of the symplectic
  Fourier transform present \(a_\psi\) as a sort of convolution:
  \begin{align}
    \nonumber 
    a_\psi(x,y)
    & =  (2\myhbar)^n \left(  \wideparen{\psi}(2x,2y) \Phi_{\tau\varsigma}(x,y) \right)
      \wideparen{\ }\\
     \label{eq:Fourier-product-conv}
    &=   {\myhbar^{2n}}
      \int_{\Space{R}{2n}} {\psi}(x',y')\,  \wideparen{\Phi}_{\tau\varsigma}(x-2x',y-2y') \,\rmd x'\,\rmd y'\\
        &=   {\myhbar^{2n}}
          \int_{\Space{R}{2n}} {\psi}(x',y')\,  {\Phi}_{\tau\varsigma}(x-2x',y-2y') \,\rmd x'\,\rmd y'\,,
\end{align}
  since \({\Phi}_{\tau\varsigma}\) is fixed by the symplectic transform: \(\wideparen{\Phi}_{\tau\varsigma} = {\Phi}_{\tau\varsigma}\).

  Finally, we rewrite the right-hand side 
  of~\eqref{eq:pdo-multiplication-intertwine}. Since \(\contravar{\varsigma}
  \circ \oper{P}_\varsigma= \contravar{\varsigma} \)  as operators  \(\FSpace{L}{2}(\Space{R}{2n})
  \rightarrow \FSpace{L}{2}(\Space{R}{n})\) we obtain:
  \begin{align*}
    \contravar{\varsigma}  \circ \psi I  \circ\covar{\phi}
    &= \contravar{\varsigma} \circ (\oper{P}_\varsigma \psi I ) \circ\covar{\tau}\,,
  \end{align*}
  where \(\oper{P}_\varsigma \psi I : \FSpace[\tau]{F}{} \rightarrow \FSpace[\varsigma]{F}{}\) is the   cross-Toeplitz operator. Since \(\contravar{\varsigma} \) and \(\covar{\tau} \)
  are unitary isomorphisms of \(\FSpace{L}{2}(\Space{R}{n})\) with \(\FSpace[\varsigma]{F}{}\) and  \(\FSpace[\tau]{F}{}\), respectively, we proved the statement.
\end{proof}

The convolution with a Gaussian is known as \emph{Weierstrass--Gauss
  transform} (also called Weierstrass or Gauss transforms
separately). Since the Gaussian provides the fundamental solution to
the \emph{heat equation},%
\index{heat equation}%
\index{equation!heat} the Weierstrass--Gauss transform can be
interpreted as a heat flow over a fixed period of time. Its appearance
in the context of FSB spaces can be traced back
to~\cite{Husimi40} at least. Being a smoothing operator the Weierstrass--Gauss
transform is not invertible on many function spaces. This puts
restrictions on applications of PDO calculus for Toeplitz operators as
was pointed in~\cite{Guillemin84} and extensively discussed
afterwords,
cf.~\cites{BergCob87,BergerCoburn94a,Coburn01b,BauerCoburnIsralowitz10a}.

The proof of Prop.~\ref{pr:toeplitz-as-PDO-Guillemin} is a rectification of ideas~\cites{Howe80b,Howe80a}
originally employed in the opposite direction: after multiplication of
both sides of~\eqref{eq:left-and-right-is-phi-rho} by \(\Phi^{-1}\), a
given PDO can be estimated through some Toeplitz operator.
Thereafter, the elementary bound of the norm of a Toeplitz
operator \(\oper{P}_\tau (\psi \Phi^{-1})I\) by
\(\norm[\infty]{\psi \Phi^{-1}}\) leads to the Calder\'on--Vaillancourt
theorem~\cite[Ch.~XIII]{MTaylor81}, which limits \(\norm{a(D,X)}\)
by \(\FSpace{L}{\infty}\)-norms of a finite number of partial
derivatives of \(a\). Prop.~\ref{pr:toeplitz-as-PDO-Guillemin} is also a specialisation of~\cite{BoggiattoCorderoGrochenig04a}*{Lem.~2.4}.

As another development of ideas from~\cite{Howe80b} we present a representation-theoretic derivation of the fundamental formula for composition of PDOs. Our treatment again relays on a simultaneous consideration of the left and right action of \(\Space{H}{n}\).

\begin{prop}
  The composition of two PDOs with symbols \(a_1\) and \(a_2\) has the symbol:
  \begin{align}
    \label{eq:PDO-composition-series}
    a &= \sum_{n,m=0}^\infty \frac{1}{n!\, m!\, (2\pi \rmi \myhbar)^{n+m}}
       (\rmd\Lambda_{\myhbar}^X - \rmd R_{\myhbar}^X)^m \,
      (\rmd\Lambda_{\myhbar}^Y - \rmd R_{\myhbar}^Y)^n \,
      a_1  \\
    \nonumber 
    &\qquad \qquad  \times
 (\rmd\Lambda_{\myhbar}^X - \rmd R_{\myhbar}^X)^n \,
      (\rmd\Lambda_{\myhbar}^Y - \rmd R_{\myhbar}^Y)^m \,
      a_2 \\
    \label{eq:PDO-composition-series-derivatives}
    &= \sum_{n,m=0}^\infty \frac{\rmi^{m-n}}{n!\, m!\, (2\pi \myhbar)^{n+m}}
       \partial_x^m \partial_y^n a_1 \, \partial_x^n \partial_y^m a_2 \,.
  \end{align}
\end{prop}
\begin{proof}
  The following computation is similar to one in~\cite{Howe80b}*{(2.2.5)}, but all elements are expressed in terms of representations of \(\Space{H}{n}\) now. For \(k_{1,2}=\wideparen{a}_{1,2}\), the symbol of  a composition of two PDOs is:
  \begin{align}
    \nonumber 
    (k_1 \twist k_2)\wideparen{\ }
    &= k_1 \twist \wideparen{k}_2 \\
    \nonumber 
    &= \Lambda_{\myhbar} (k_1) \, \wideparen{k}_2 \\
    \nonumber 
    &= \left(\int k_1(x,y) \,  \exp \left( x\,\rmd\Lambda_{\myhbar}^X + y\,\rmd\Lambda_{\myhbar}^Y\right) \rmd x \, \rmd y\right) \wideparen{k}_2\\
    \nonumber 
    &= \left(\int k_1(x,y) \, \exp \left(\half x(\rmd\Lambda_{\myhbar}^X - \rmd R_{\myhbar}^X) + \half y(\rmd\Lambda_{\myhbar}^Y - \rmd R_{\myhbar}^Y)\right) \right.\\
    \nonumber 
    &\qquad \qquad \left. \phantom{\int} \times
      \exp \left(\half x(\rmd\Lambda_{\myhbar}^X + \rmd R_{\myhbar}^X) + \half y(\rmd\Lambda_{\myhbar}^Y + \rmd R_{\myhbar}^Y)\right)  \rmd x \, \rmd y\right) \wideparen{k}_2\\
    \intertext{where we used \(\rme^{A/2+B/2}=\rme^{A/2}\rme^{B/2}\) for commuting operators \(A=x(\rmd\Lambda_{\myhbar}^X - \rmd R_{\myhbar}^X) + y(\rmd\Lambda_{\myhbar}^Y - \rmd R_{\myhbar}^Y)\) and \(B=x(\rmd\Lambda_{\myhbar}^X + \rmd R_{\myhbar}^X) + y(\rmd\Lambda_{\myhbar}^Y+ \rmd R_{\myhbar}^Y)\). Furthermore, commutativity of the left and right actions and vanishing commutator \([\rmd \Lambda_\myhbar^X - \rmd R_\myhbar^X, \rmd \Lambda_\myhbar^Y - \rmd R_\myhbar^Y]=0\) imply:}
    \nonumber 
    &= \left(\int k_1(x,y) \,  \exp \left(\half x(\rmd\Lambda_{\myhbar}^X - \rmd R_{\myhbar}^X)\right) \exp \left(\half y(\rmd\Lambda_{\myhbar}^Y - \rmd R_{\myhbar}^Y)\right) \right.\\
    \nonumber 
    &\qquad \qquad \left. \phantom{\int} \times
     \Lambda_{\myhbar}(\half x,\half y)\circ R_{\myhbar} (\half x,\half y)\,  \rmd x \, \rmd y\right) \wideparen{k}_2\\
    \nonumber 
    &= \left(\int k_1(x,y) \, \sum_{n,m=0}^\infty \frac{1}{n! \, m! \, 2^{n+m}}
      x^n(\rmd\Lambda_{\myhbar}^X - \rmd R_{\myhbar}^X)^n \, y^m(\rmd\Lambda_{\myhbar}^Y - \rmd R_{\myhbar}^Y)^m \right.\\
    \nonumber 
    &\qquad \qquad \left. \phantom{\int} \times
           \Lambda_{\myhbar}(\half x,\half y)\circ R_{\myhbar} (\half x,\half y)\,  \rmd x \, \rmd y\right) \wideparen{k}_2\\
    \nonumber 
    &= \sum_{n,m=0}^\infty \frac{1}{n! \, m! \, 2^{n+m}} \int k_1(x,y) \, 
      x^n  y^m    \Lambda_{\myhbar}(\half x,\half y)\circ R_{\myhbar} (\half x,\half y)\,   \rmd x \, \rmd y \\
    \nonumber 
    &\qquad \qquad  \times
 (\rmd\Lambda_{\myhbar}^X - \rmd R_{\myhbar}^X)^n
      (\rmd\Lambda_{\myhbar}^Y - \rmd R_{\myhbar}^Y)^m
      \wideparen{k}_2\\
    \intertext{The last integral is the symplectic Fourier transform of \(x^n  y^m k_1\) under the derived representation, thus the intertwining properties~\eqref{eq:fourier-intertw-shift-mult-compl} imply:}
    \nonumber 
    &= \sum_{n,m=0}^\infty \frac{\rmi^{m-n}}{n!\, m!\, (2\pi \myhbar)^{n+m}}
       (\rmd\Lambda_{\myhbar}^X - \rmd R_{\myhbar}^X)^m \,
      (\rmd\Lambda_{\myhbar}^Y - \rmd R_{\myhbar}^Y)^n \,
      \wideparen{k}_1  \\
    \nonumber 
    &\qquad \qquad  \times
 (\rmd\Lambda_{\myhbar}^X - \rmd R_{\myhbar}^X)^n \,
      (\rmd\Lambda_{\myhbar}^Y - \rmd R_{\myhbar}^Y)^m \,
      \wideparen{k}_2 \,,
  \end{align}
  that is the first form~\eqref{eq:PDO-composition-series} of the composition.  Using the expressions for the derived representations~\eqref{eq:left-action-derived}--\eqref{eq:right-action-derived}  we obtain~\eqref{eq:PDO-composition-series-derivatives}.
\end{proof}

A parallel computation with complex variables implies:
\begin{cor}
  The composition of two PDOs with symbols \(a_1\) and \(a_2\) has the symbol:
  \begin{align}
a      &= \sum_{n,m=0}^\infty \frac{(-1)^m}{n! \, m! \, 2^{n+m}}            \left(\ladder[R_\myhbar]{-}- \ladder[ \Lambda_\myhbar]{-} \right)^n  \left(  \ladder[R_\myhbar]{+}-\ladder[\Lambda_\myhbar]{+} \right)^m
        a_1\\
      \nonumber 
&\qquad \qquad \times 
    \left(\ladder[\Lambda_\myhbar]{+} - \ladder[R_\myhbar]{+}\right)^n  \left( \ladder[\Lambda_\myhbar]{-} -   \ladder[R_\myhbar]{-} \right)^m a_2\\
  &= \sum_{n,m=0}^\infty \frac{(-2)^{n+m}}{n! \, m! \, }       \partial_z^n \overline{\partial }_z^m 
    a_1 \cdot 
    \overline{\partial}_z^n  \partial_z^m a_2\,.
  \end{align}
\end{cor}

\subsection{Bounded operators and integrated representation}
\label{sec:pdo-representing-an}

It does not worth to be locked within the particular case of the Heisenberg
group while using group representation techniques. In fact, it is
often possible to get a more general statement at no extra cost,
Subsection~\ref{sec:right-shifts-analyt} provides an example. On the
other hand, the ultimate generality is not our purpose here either.
Thus, we are working in the comfortable assumption of a unimodular
group \(G\) and an \emph{irreducible square integrable} (possibly, \emph{modulo a
subgroup} \(H\)) representation \(\uir{}{}\) in a Hilbert space
\(\FSpace{H}{}\). Then, the orthogonality
relation~\eqref{eq:wigner-sesqui-linear} for
\(\oper{W}(f,\phi)(g)=\scalar{f}{\uir{}{}(g)\phi}\)
holds~\citelist{\cite{AliAntGaz14a}*{Thm.~8.2.1}
  \cite{FeichtingerPap14a}*{Thm.~1}}.

For an irreducible representation \(\uir{}{}\) of a group \(G\) on a
Hilbert space \(\FSpace{H}{}\), the Schur lemma implies that the
von~Neumann algebra generated by \(\uir{}{}(g)\), \(g\in G\) contains
all bounded linear operators on \(\FSpace{H}{}\). Thus, any such
operator can be considered as an integrated representation
\(\uir{}{}(k)\) with a suitable kernel \(k\) (possibly
being a distribution).
To find the kernel explicitly consider the induced covariant transform of the representation \(\uir{}{}\):
\begin{equation}
  \label{eq:wavelet-transform-product}
  b(g)= \oper{W}(u,v)(g)=\scalar[\FSpace{H}{}]{u}{\uir{}{}(g)v}\,,
  \quad \text{where } u, v \in \FSpace{H}{}
\end{equation}
as a map from \(\FSpace{H}{}\times \FSpace[*]{H}{}\) to a space of function \(\FSpace{L}{}(X)\) on the homogeneous space \(X=G/H\). This linear map \(\FSpace{H}{}\times \FSpace[*]{H}{} \rightarrow \FSpace{L}{}(X)\) can be extended to the tensor product \(\FSpace{H}{}\otimes \FSpace[*]{H}{} \rightarrow \FSpace{L}{}(X)\).

\begin{prop}
  \label{pr:integrated-representation-symbol}
  Let \(G\), \(H\), \(\FSpace{H}{}\) and \(\uir{}{}\) be as above.  For a
  Hilbert--Schmidt operator \(K\) considered as an element of
  \(\FSpace{H}{}\otimes \FSpace[*]{H}{}\), let us define a function
  \(b_K\) by the extension of~\eqref{eq:wavelet-transform-product} 
  from \(\FSpace{H}{}\times \FSpace[*]{H}{}\) to
  \(\FSpace{H}{}\otimes \FSpace[*]{H}{}\). Then, \(K\) is equal to the
  integrated representation: \(K=\uir{}{}(b_K)\).
\end{prop}
\begin{proof}
  For \(u\), \(v\in \FSpace{H}{}\), let \(K: \FSpace{H}{} \rightarrow \FSpace{H}{}\) be
  the rank one operator \(K: f \mapsto \scalar[\FSpace{H}{}]{f}{v} u\). 
  For the corresponding integrated representation \(\uir{}{}(b)\) and
  arbitrary \(f_1\), \(f_2\in \FSpace{H}{}\) with square-integrable
  matrix coefficients:
  \begin{align*}
    \scalar[\FSpace{H}{}]{\uir{}{}(b) f_1}{f_2}
      &= \scalar[\FSpace{H}{}]{\int_X b(g)\,\uir{}{}(g) f_1\,\rmd g}{f_2} \\ 
      &= \int_X b(g) \scalar[\FSpace{H}{}]{\uir{}{}(g) f_1}{f_2}\,\rmd g \\ 
      &= \int_X b(g)\,\overline{\scalar[\FSpace{H}{}]{f_2}{\uir{}{}(g) f_1}}\,\rmd g\\
      &= \scalar[\FSpace{L}{2}(X)] {\oper{W}(u,v)}{\oper{W}(f_2,f_1)}\\
        \intertext{
  then using the orthogonality relation~\eqref{eq:wigner-sesqui-linear}
  for wavelet transform we continue:}
    &=\scalar[\FSpace{H}{}]{u}{f_2} \overline{\scalar[\FSpace{H}{}]{v}{f_1}}\\
    &=\scalar[\FSpace{H}{}]{u}{f_2} \scalar[\FSpace{H}{}]{f_1}{v}\\
    &=\scalar[\FSpace{H}{}]{\scalar[\FSpace{H}{}]{f_1}{v}u}{f_2} \\
    &=\scalar{Kf_1}{f_2}.
  \end{align*}
  Since \(f_1\) and \(f_2\) with the square-integrable property are
  dense in \(\FSpace{H}{}\), the identity \(K=\uir{}{}(b)\) is shown
  in the simplest case \(K\in \FSpace{H}{}\times
  \FSpace[*]{H}{}\). Then, this identity generalises to any Hilbert--Schmidt
  operator  by the respective linear
  extension~\eqref{eq:wavelet-transform-product}.
\end{proof}

The above correspondence between operators and integrated
representations can be pushed beyond Hilbert--Schmidt and bounded operators on Hilbert
spaces through the \emph{coorbit theory}
(see~\cites{FeichGroech89a,FeichGroech89b} and the recent
survey~\cite{FeichtingerPap14a}) or wavelets in Banach spaces
technique~\cites{Kisil98a,Kisil13a,Kisil94e}. A motivating  example
is the case of \(G=\Space{H}{n}\), \(H\)---its centre \(Z\) and
\(\uir{}{}\)---the Schr\"odinger
representation~\eqref{eq:schroedinger-rep} in
\(\FSpace{H}{}=\FSpace{L}{2}(\Space{R}{n})\), which leads to the
following statement.

\begin{cor}[\citelist{\cite{Howe80b}*{\S~2.3} \cite{Folland89}*{\S~2.1}}]
  \label{pr:kernal-as-PDO-Wigner}
  Let \(K: \FSpace{S}{}(\Space{R}{n}) \rightarrow
  \FSpace[\prime]{S}{}(\Space{R}{n})\) be a linear operator with
  Schwartz kernel \(k\):
  \begin{equation}
    \label{eq:Schwartz-kernel}
    K: f(x) \mapsto [Kf](x)=\int_{\Space{R}{n}} k(x,y) \, f(y)\,\rmd y.
  \end{equation}
  Then \(K\) is PDO \(a_K(D,X)\) with the symbol \(a_K\) such that: 
  \begin{equation}
    \label{eq:PDO-symbol-from-kernel}
    a_K(x,y)=  \int_{\Space{R}{n}} \rme^{-\pi \rmi \myhbar
      (2t-y)x}\, k(t,y-t)\,\rmd t\,.
  \end{equation}
\end{cor}
\begin{proof}
  The extension of~\eqref{eq:wavelet-transform-product} to \( \FSpace{S}{}(\Space{R}{n}) \otimes   \FSpace[\prime]{S}{}(\Space{R}{n})\) for the integral operator \(K\)~\eqref{eq:Schwartz-kernel} obviously is:
  \begin{displaymath}
    b(g)=\int_{\Space{R}{n}} \uir{}{}(g) k(x,x)\,\rmd x,
  \end{displaymath}
  where the representation \(\uir{}{}(g)\) acts on the second variable of the kernel \(k(x,y)\).
  Then, for the Schr\"odinger representation~\eqref{eq:schroedinger-rep} it becomes:
  \begin{displaymath}
    b_K(x',y')= \int_{\Space{R}{n}} \rme^{-\pi \rmi \myhbar
      (2t-y')x'}\, k(t,t-y')\,\rmd t\,.
  \end{displaymath}
   By the previous Proposition, \(K\) is the integrated representation    of the function \(b_K\). An application of the symplectic Fourier transform \( {a}_K = \wideparen{b}_K\)   produces~\eqref{eq:PDO-symbol-from-kernel}, with tiny differences   between two expressions.
\end{proof}

For applications to cross-Toeplitz operators we again use a more general
setup and notations (\(G\), \(H\), \(X\), \(\chi\), \(\map{s}\), \(\phi\), \ldots) of induced
representations from Rem.~\ref{re:induced-representation}. An
additional assumption is: there are \(\uir{}{\chi}\)-irreducible unitary equivalent
components  \(\FSpace[\tau]{F}{}\) and \(\FSpace[\varsigma]{F}{}\) of \(\FSpace{L}{2}(X)\) with an intertwining integral kernel \(K_{\tau\varsigma}(w,z)\), which also provides an orthogonal
projection
\(\oper{P}_\varsigma: \FSpace{L}{2}(X)\rightarrow \FSpace[\varsigma]{F}{}\). A
reader, who is not interested in group representations, may safely
remain within the FSB space framework described in
subsection~\ref{sec:four-wign-transf}. However, there are other
interesting implementations of the scheme, which include the Hardy and
weighted Bergman spaces on the unit disk and upper half-plane
\citelist{\cite{Kisil97c} \cite{Kisil09e} \cite{Kisil13a}
  \cite{Kisil11c} \cite{FeichtingerPap14a}
  \cite{ChristensenOlafsson09a}}.
\begin{cor}
  \label{eq:toeplitz-as-convolution}  
  In the above assumptions and using notations of
  Rem.~\ref{re:induced-representation}, for a bounded function \(\psi\) on
  \(X\), the Toeplitz operator
  \(\oper{P}_\tau \psi I: \FSpace[\tau]{F}{} \rightarrow
  \FSpace[\tau]{F}{}\) is an integrated representation
  \(\uir{}{\chi}(b_\psi)\) for the function:
  \begin{align}
    \label{eq:Toeplitz-as-integrated-kernel}
    b_\psi(g)&=\int_X {\chi}\left(\map{r}(g^{-1}*\map{s}(z))\right)
               K_{\tau}(z,g^{-1}\cdot z ) \, {\psi}(g^{-1}\cdot z) \,\rmd z 
    \,,
  \end{align}
  where \(K_{\tau}(w,z)\) is the
  reproducing kernel~\eqref{eq:FSB-repro-kernel}. In particular, for FSB spaces:  
  \begin{equation}
    \label{eq:Toeplitz-as-integrated-kernel-FSB}
    \begin{split}
    b_\psi(g)&=\int_X
               {\chi}\left(\map{r}\left(g^{-1}*\map{s}(z)\right)\right)
               \overline{\chi}\left(\map{r}\left(\map{s}(g^{-1}\cdot z)^{-1}
               *\map{s}(z)\right)\right)\\
    &\qquad \times 
      \Phi_{\tau}((\map{s}(z)^{-1}*g)\cdot z)\,
      {\psi}(g^{-1}\cdot z) \,\rmd z
    \,,
  \end{split}
  \end{equation}
\end{cor}
\begin{proof}
  Obviously, the Toeplitz operator \(\oper{P}_\tau \psi I\) is an
  integral operator on \(\Space{C}{n}\) with the kernel
  \(k_\psi(w,z)=K_{\tau}(w,z){\psi}(z)\). Then, for the
  induced representation
  \(\uir{}{\chi}\)~\eqref{eq:induced-representation-def} the
  Wigner-type transform~\eqref{eq:wavelet-transform-product} produces
  the first form~\eqref{eq:Toeplitz-as-integrated-kernel}. The
  second form~\eqref{eq:Toeplitz-as-integrated-kernel-FSB} follows from the expression of the reproducing kernel
  \begin{displaymath}
    K_{\tau}(w,z)=\scalar{\uir{}{}(\map{s}(z)) \phi_\tau}{\uir{}{}(\map{s}(w)) \phi_\tau} =
    \overline{\chi}\left(\map{r}\left((\map{s}(z))^{-1}*\map{s}(w)\right)\right)
      \Phi_{\tau}((s(z))^{-1}\cdot w)
  \end{displaymath}
  through the wavelet transform of vacuum vectors (e.g. Gaussians, cf.~\eqref{eq:reproducing-kernel}).
\end{proof}
The above formula~\eqref{eq:Toeplitz-as-integrated-kernel-FSB} greatly
simplifies for the Heisenberg group and  the (pre-)FSB space. The
respective substitution produces the next result:
\begin{cor}
  The  Toeplitz operator  \(\oper{P}_\tau \psi I: \FSpace[\tau]{F}{} \rightarrow \FSpace[\tau]{F}{}\) with a symbol   \(\psi\)  is a relative convolution with the kernel:
\begin{align}
  \label{eq:Toeplitz-integrated-kernel}
  k_\psi(x,y)
   &=\int_{\Space{R}{2n}}
   \rme^{2\pi \rmi \myhbar (xy'-yx')}\,
     \Phi_{\tau}(x,y)\,{\psi}(x'-x,y'-y)\,\rmd x' \, \rmd y'\\
  \label{eq:toeplitz-kernel-fourier}
    &=\Phi_{\tau}(x,y)\, \wideparen{\psi}(2x,2y)\,.
\end{align}
\end{cor}
The result is a particular form for \(\tau=\varsigma\) of the identity~\eqref{eq:Toeplitz-PDO-diffused} since
\( \wideparen{k}_\psi(x,y)\) is equal to \(a(x,y)\)
from~\eqref{eq:Toeplitz-PDO-diffused}---the symbol of PDO representing
the cross-Toeplitz operator \(\oper{P}_{\tau} \psi I\) (with all its
deficiencies discussed above). A related expression can be found
in~\cite{BoggiattoCorderoGrochenig04a}*{\S~2.3}. We cannot get a cross-Toeplitz operator as an integrated representation of the left action since it preserves spaces \(\FSpace[\tau]{F}{} \).
Thus, for a better description of cross-Toeplitz operators we need
an expanded group action discussed in the next section.
\begin{rem}
  This section's methods and results are related to co- and contra-variant Berezin calculus, also known as Wick and anti-Wick symbolic calculus~\cite{Berezin86}, cf.~\citelist{\cite{Coburn01b} \cite{BoggiattoCorderoGrochenig04a} \cite{Folland89}*{\S~2.7}}. In the approach advocated here it falls into the general framework of operator covariant transform~\cites{Kisil98a,Kisil94e,Kisil10c,Kisil13a}, yet we keep this technique in low profile to avoid unnecessary abstraction. 
\end{rem}

\section{Two-sided convolutions and cross-Toeplitz operators}
\label{sec:two-sided-conv-toeplltz}
The previous consideration repeatedly used combinations of the left and right pulled actions~\eqref{eq:left-right-action-pulled-def}--\eqref{eq:left-action-pulled}. It is time to employ their union in a systematic way.

\subsection{Two-sided relative convolutions from the Heisenberg group}
\label{sec:two-sided-conv}

The relation~\eqref{eq:mult-by-left-right} (repeated in~\eqref{eq:left-times-right-rep}) suggests to present a general operator of multiplication
\(\psi I: f\mapsto \psi f\) as an integrated representation, cf.~\eqref{eq:relative-conv}:
\begin{align}
  \label{eq:multiplication-is-convolution}
  \psi I
  &=(\Lambda_\myhbar \otimes R_\myhbar)(\wideparen{\psi})\\
  \nonumber
        &=(2/\myhbar)^{n}\int_{\Space{R}{2n}}
          \wideparen{\psi}(2x,2y)\,
          \Lambda_\myhbar(0,x,y)\, R_\myhbar(0,x,y)\,dx\,dy\,,
\end{align}
where \(\wideparen{\psi}\) is the symplectic Fourier
transform~\eqref{eq:symplect-Fourier-defn} of \(\psi\).
Motivated by~\eqref{eq:multiplication-is-convolution} we introduce the
main tool of our investigation.
\begin{defn}
  For a function \(k\) on \(\Space{R}{4n}\), a \emph{two-sided
    relative convolution} \(\oper{D}(k)\) is defined by:
  \begin{equation}
    \label{eq:two-sided-convolution}
    \oper{D}(k)=
    \int_{\Space{R}{4n}}
    k(x_1,y_1,x_2,y_2)\,\Lambda_\myhbar(0,x_1,y_1)\,
    R_\myhbar(0,x_2,y_2)\,
    \rmd x_1 \,\rmd y_1 \,\rmd x_2 \,\rmd y_2\, .
  \end{equation}
\end{defn}
Such operators also arise as
representations of two-sided convolutions on the Heisenberg group
studied in~\citelist{\cite{VasTru88} \cite{VasTru94}
  \cite{Vasilevski94a} \cite{Kisil92} \cite{Kisil93e} \cite{Kisil93b}
  \cite{Kisil94a} \cite{Kisil94f} \cite{Kisil96e}} as local
(irreducible) representations of certain \(C^*\)-algebras. The connection
between these two sources was studied in~\cite{Kisil98e,Kisil12b}.

Since operators \(\Lambda_\myhbar\) and \(R_\myhbar\) are isometries
in \(\FSpace{L}{p}(\Space{R}{2n})\), \(1\leq p \leq \infty\) the
integral~\eqref{eq:two-sided-convolution} defines a bounded operator
in \(\FSpace{L}{p}(\Space{R}{2n})\) for any
\(k\in\FSpace{L}{1}(\Space{R}{4n})\). Alternatively, one can consider any shift-invariant
norm from Lem.~\ref{le:shift-invariant-norm}. Furthermore, the
Schwartz space \(\FSpace{S}{}(\Space{R}{2n})\) of smooth rapidly
decreasing functions is invariant under \(\Lambda_\myhbar\) and
\(R_\myhbar\).  Thus, for any \(u\),
\(v\in\FSpace{S}{}(\Space{R}{2n})\), the function
\begin{displaymath}
   \tilde{v}_u(x_1,y_1,x_2,y_2)=\scalar{v}{\Lambda_\myhbar(0,x_1,y_1)
     R_\myhbar(0,x_2,y_2) u} 
\end{displaymath}
is in \(\FSpace{S}{}(\Space{R}{4n})\). Therefore, a distribution
\(k\in\FSpace[\prime]{S}{}(\Space{R}{4n})\) defines a bounded operator
\(\oper{D}(k): \FSpace{S}{}(\Space{R}{2n}) \rightarrow
\FSpace[\prime]{S}{}(\Space{R}{2n})\) in the weak sense:
\(\scalar{\oper{D}(k)u}{v}=\scalar{k}{\tilde{v}_u}\).

\begin{example}
  \label{ex:as-2-side-conv}
  Several operators relevant to our consideration can be treated as
  two-sided convolutions.
  \begin{enumerate}
  \item The operator of multiplication \( f\mapsto \psi f\) by a
    function \(\psi\) is a two-sided convolution \(\oper{D}(k)\) with
    the distribution, cf.~\eqref{eq:multiplication-is-convolution}
    and~\eqref{eq:two-sided-convolution} 
    \begin{equation}
      \label{eq:two-kernel-multipl}
      k(x_1,y_1,x_2,y_2)=(2/\myhbar)^{n}\,
      \wideparen{\psi} (2x_1,2y_1)\, \delta(x_1-x_2,y_1-y_2),
    \end{equation}
    where \(\wideparen{\psi}\) is the symplectic Fourier transform of \(\psi\) and
    \(\delta\) is the Dirac delta function. 
  \item Integrated representations
    \(\Lambda_\myhbar(k)\) or \(R_\myhbar(k)\) with a kernel \(k\) on
    \(\Space{R}{2n}\) are two-sided convolutions
    \(\oper{D}(k\otimes \delta )\) or \(\oper{D}(\delta \otimes k )\),
    respectively. In particular, the symplectic Fourier
    transform~\eqref{eq:symplect-Fourier-left} is \(\oper{D}(\mathbf{1} \otimes \delta ) = \oper{D}(\delta \otimes \mathbf{1})\). 

    
  \item The orthogonal projection
    \(\oper{P}_{\tau}: \FSpace{L}{2}(\Space{R}{2n}) \rightarrow
    \FSpace[\tau]{F}{}\)~\eqref{eq:reproducing-projection} can be rearranged as
    an integrated representation \(R_\myhbar\), cf.~\eqref{eq:right-through-left} :
    \begin{align}
       \nonumber 
      [\oper{P}_{\tau} f](x,y)&=\scalar{f}{\Lambda_{\myhbar}(0,x,y)\Phi_{\tau}}\\
      \nonumber 
                            &=\scalar{\Lambda_{\myhbar}(0,-x,-y)f}{\Phi_{\tau}}\\
      \nonumber 
                            &=[(R_\myhbar(\overline{\Phi}_{\tau}) f] (x,y) \\
      \label{eq:reproducing-projection-double}
                               &=[(R_\myhbar({\Phi}_{\tau}) f] (x,y)\,
                                  ,
    \end{align}
    because \(\overline{\Phi}_{\tau}=\Phi_{\tau}\), see~\eqref{eq:Phi-tau-sigma}. 
    Therefore, \(\oper{P}_{\tau}\) is a two-sided convolution
    \begin{equation}
      \label{eq:fock-proj-two-conv}
      \oper{P}_{\tau} =\oper{D}(\delta \otimes {\Phi}_{\tau})
      \, .
    \end{equation}
\item    Similarly the intertwining operator  \(\oper{P}_{\tau\varsigma}: 
  \FSpace[\tau]{F}{} \rightarrow  \FSpace[\varsigma]{F}{}\) is a two-sided convolution with the kernel
    \begin{equation}
      \label{eq:intertwining-two-conv}
      \oper{P}_{\tau\varsigma} =\oper{D}(\delta \otimes {\Phi}_{\varsigma\tau})\,.
    \end{equation}
  Note the swapped order of squeeze parameters due to \(\overline{\Phi}_{\tau\varsigma}=\Phi_{\varsigma\tau}\).
  \item A (cross-)Toeplitz operator
    \(\oper{T}_\psi: f \mapsto \oper{P}_{\varsigma} (\psi f)\) is the
    composition of two-sided convolutions with
    kernels~\eqref{eq:two-kernel-multipl}
    and~\eqref{eq:fock-proj-two-conv}, which will be explicitly evaluated in
    Sect.~\ref{sec:comp-oper}.
\end{enumerate}
\end{example}

A simple change of variables in the integral shows that
\begin{lem}
  A two-sided convolution \(\oper{D}(k)\) with the kernel \(k(x_1,y_1;x_2,y_2)\) is an integral operator
  \begin{displaymath}
    [\oper{D}(k)f] = \int_{\Space{R}{2n}} \tilde{k}(x,y;x',y')\, f(x',y')
    \,\rmd x' \, \rmd y'
  \end{displaymath}
  with the Schwartz kernel
  \begin{equation}
    \label{eq:Schwartz-kernel-from-two-sided}
    \tilde{k}(x,y;x',y') = \int_{\Space{R}{2n}} \rme^{\pi \rmi \myhbar
                                      ((2x_2-x')y-(2y_2-y')x)} k(x-x'+x_2,y-y'+y_2;x_2,y_2)\,\rmd x_2 \, \rmd y_2.
  \end{equation}
\end{lem}
We can transit between Schwartz and convolution kernels in the opposite direction using Prop.~\ref{pr:integrated-representation-symbol}, cf. Cor.~\ref{pr:kernal-as-PDO-Wigner}:
\begin{cor}
  \label{co:two-sided-conv-from-Schwartz}
  An integral operator \(K\) on \(\Space{R}{2n}\) with a Schwartz kernel \(\tilde{k}(x_1,y_1;x_2,y_2)\) is equal to two-sided convolution  \(\oper{D}(k)\) with the kernel:
  \begin{equation}
    \label{eq:two-sided-kernel-from-Schwartz}
    k(x_1,y_1;x_2,y_2)=\int_{\Space{R}{2n}} \Lambda_\myhbar(-x_2,-y_2) R_\myhbar(-x_1,-y_1)
    \tilde{k}(x,y;x,y)\,\rmd x\,\rmd y,
  \end{equation}
  there both representations act on the second half of variables in \(\tilde{k}\).
\end{cor}
Reciprocity of formulae~\eqref{eq:Schwartz-kernel-from-two-sided} and~\eqref{eq:two-sided-kernel-from-Schwartz} can be checked by a direct calculation.

\subsection{Reduction of two-sided convolutions}
\label{sec:two-sided-conv-PDO}

The above examples indicate that an effective calculus of two-sided
convolutions is useful for the theory of cross-Toeplitz operators.
To this end we will reduce two-sided convolutions from
\(\Space{H}{n}\) to one-sided convolution on the Heisenberg
group \(\Space{H}{2n}\) of the doubled size~\cite{Kisil94f}, cf.~\cite{Street08a}
for applications of this method.

To begin we note,  that composition of two-sided convolutions on the entire group \(G\) is related to the left convolutions on the Cartesian product group \(G\times G\), that is we have   \(\oper{D}(k_1)\oper{D}(k_2)=\oper{D}(k)\) for:
\begin{equation}
  \label{eq:composition-two-sided-group}
  k(g_1,g_2)=\int_G\int_G k_1(g_3,g_4)k_1(g_3^{-1}g_1,g_4^{-1}g_2)\,\rmd g_3\, \rmd g_4\,.  
\end{equation}
However twisted convolutions for homogeneous spaces require a more accurate treatment shown below.

Let us introduce the action of \(\Space{H}{2n}\) on
\(\FSpace{L}{2}(\Space{R}{2n})\) by:
\begin{equation}
  \label{eq:pre-Xi-B-composition}
  \widetilde{\Xi}_\myhbar=(\Lambda_\myhbar\otimes R_\myhbar)\circ Y.
\end{equation}
Here the section
\begin{equation}
  Y: \Space{H}{2n} \rightarrow  \Space{H}{n}\times \Space{H}{n} :
  (s,x_1,x_2,y_1,y_2) \mapsto (s,x_1,y_1)\times (0,\upsilon
  y_2,x_2/\upsilon )
\end{equation}
is a right inverse of the group homomorphism \( \Space{H}{n}\times \Space{H}{n} \rightarrow (\Space{H}{n}\times \Space{H}{n})/Z_d  \simeq\Space{H}{2n} \), where
\begin{displaymath}
  Z_d=\{(s,0,0)\times (s,0,0) \such s\in \Space{R}{}\}
\end{displaymath}
is the diagonal of \(Z\times Z \subset \Space{H}{n}\times \Space{H}{n}\). 
In other words:
\begin{align}
  \label{eq:H-2n-action-two-sided}
  \widetilde{\Xi}_\myhbar(s,x_1,x_2,y_1,y_2)
  &=\Lambda_\myhbar(s,x_1,y_1) \, R_\myhbar(0,\upsilon y_2,x_2/\upsilon )\\
  \nonumber 
  &= \Lambda_\myhbar(0,x_1,y_1) \, R_\myhbar(-s,\upsilon y_2,x_2/\upsilon ),
    \
\end{align}
where \( s\in\Space{R}{}\) and \(x_1\), \(x_2\), \(y_1\), \(y_2\in\Space{R}{n}\). Note,
that the above formulae contain an arbitrary constant squeeze parameter \(\upsilon \neq 0\). From the physical point of view it adjusts units of the respective components of the phase space, cf. Rem.~\ref{re:physical-units}. The role of \(\upsilon\) is rather technical and it will disappear at later stage, cf. Rem.~\ref{re:upsilon-vanish}.

It is a straightforward calculation that \(\widetilde{\Xi}_\myhbar\) is an irreducible
unitary representation of \(\Space{H}{2n}\) on
\(\FSpace{L}{2}(\Space{R}{2n})\). The Stone--von Neumann theorem guarantees equivalence of \(\widetilde{\Xi}_\myhbar\) to the Schr\"odinger representation  and this is completely transparent from  the unitary map of a partial dilation
\begin{equation}
  \label{eq:remapping-squeeze}
  S_\upsilon :\; \FSpace{L}{2}(\Space{R}{2n},\rmd x\,\rmd y) \rightarrow
  \FSpace{L}{2}(\Space{R}{2n},\upsilon ^n\rmd t_1\,\rmd t_2):\;
  f(x,y) \mapsto \tilde{f}(t_1,t_2)=f(\upsilon t_1,t_2)\,.
\end{equation}
We conjugate the action \(\widetilde{\Xi}_\myhbar\) of \(\Space{H}{2n}\) on
\(\FSpace{L}{2}(\Space{R}{2n})\) by \(S_\upsilon\)
\begin{equation}
  \label{eq:Xi-S-composition}
  \Xi_\myhbar(g)= S_\upsilon \,  \widetilde{\Xi}_\myhbar(g) \, S_\upsilon ^{-1}, \qquad
  \text{ for } g\in \Space{H}{2n}\,.
\end{equation}
We do not indicate dependence of \(\Xi_\myhbar\) on \(\upsilon\) to keep notation simpler.
Substituting the above expressions 
in~\eqref{eq:left-action-pulled}--\eqref{eq:right-action-pulled} we
rewrite \(\Xi_\myhbar(g)\)~\eqref{eq:Xi-S-composition} explicitly as:
\begin{align}
  \nonumber 
    [\Xi_\myhbar(s,x_1,x_2,y_1,y_2)f](t_1,t_2)&=
    \chi_ \myhbar{
      (2s-(\upsilon  y_1+x_2)t_1+(x_1+\upsilon y_2)t_2+x_1x_2/\upsilon
      -\upsilon y_1y_2)}\\
  \label{eq:H-2n-action-two-sided-explicit}
    &\qquad\times f(t_1-x_1/\upsilon +y_2,t_2-y_1+x_2/\upsilon )\,.
\end{align}
The essential properties of \(\Xi_\myhbar\) are summarised as follows.
\begin{lem}
  The action \(\Xi_\myhbar\)~\eqref{eq:H-2n-action-two-sided} is a
  linear unitary irreducible representation of \(\Space{H}{2n}\) in
  \(\FSpace{L}{2}(\Space{R}{2n})\). Moreover, the Schr\"odinger
  representation~\eqref{eq:schroedinger-rep} of \(\Space{H}{2n}\)
  in \(\FSpace{L}{2}(\Space{R}{2n})\) is the composition of
  \(\Xi_\myhbar\) with the symplectic automorphism \(A\) of
  \(\Space{H}{2n}\):
  \begin{equation}
    \label{eq:two-sided-2-schroed}
    A:\; (t,x_1,x_2,y_1,y_2) \mapsto \textstyle
    (t,-\frac{1}{2} x_2-\frac{1}{2}\upsilon y_1,\frac{1}{2} x_1+\frac{1}{2} \upsilon y_2,
    x_1/\upsilon-y_2, - x_2/\upsilon+y_1)\,.
  \end{equation}
  The above symplectic map  \(A: (x_1,x_2,y_1,y_2) \mapsto
  (\tilde{x}_1,\tilde{x}_2,\tilde{y}_1,\tilde{y}_2)
  \) in matrix form is:
   \begin{equation}
    \label{eq:two-sided-2-schroed-matr}
     \begin{pmatrix}
       \tilde{x}_1\\\tilde{x}_2\\ \tilde{y}_1\\ \tilde{y}_2
     \end{pmatrix}
     =\begin{pmatrix}
       0&-\frac{1}{2} &-\frac{\upsilon}{2} &0  \\[.3em]
       \half& 0&0 & \frac{\upsilon}{2} \\[.3em]
       \frac{1}{\upsilon}&0 &0 &-1  \\[.3em]
       0& -\frac{1 }{\upsilon}&1 & 0  
     \end{pmatrix}
     \begin{pmatrix}
       x_1\\x_2\\y_1\\y_2
     \end{pmatrix}\,.
   \end{equation}
\end{lem}
\begin{proof}
  It follows from the group law~\eqref{eq:H-n-group-law} that a symplectic
  transformation of \(\Space{R}{4n}\)---in
  particular~\eqref{eq:two-sided-2-schroed}---produces an automorphism
  of \(\Space{H}{2n}\)~\cite{Folland89}*{\S~1.2}. It is another
  straightforward check, that the composition
  of \(A\)~\eqref{eq:two-sided-2-schroed} and the
  action~\eqref{eq:H-2n-action-two-sided-explicit} produces the
  Schr\"odinger representation~\eqref{eq:schroedinger-rep} of
  \(\Space{H}{2n}\) in \(\FSpace{L}{2}(\Space{R}{2n})\). A composition
  of a group automorphism with a representation produces again a
  representation of the group, thus \(\Xi_\myhbar\) as a composition
  of the inverse of~\eqref{eq:two-sided-2-schroed} and the Schr\"odinger
  representation is again a representation of \(\Space{H}{2n}\).  The
  rest follows from the properties of the Schr\"odinger representation.
\end{proof}
A two-sided convolution can be written as the integrated
representation  \(\widetilde{\Xi}_\myhbar\)  using a change of variables:
\begin{align}
  \nonumber 
  \oper{D}(k)
  &=\int_{\Space{R}{4n}} k(x_1,y_1,x_2,y_2)\,
    \Lambda_\myhbar(0,x_1,y_1)\, R_\myhbar(0,x_2,y_2)
    \,\rmd x_1\,\rmd x_2\,\rmd y_1\,\rmd y_2\\
  \nonumber 
  &=\int_{\Space{R}{4n}}
    k(x_1,y_1,x_2,y_2)\, \widetilde{\Xi}_\myhbar(0,x_1,\upsilon y_2,y_1,x_2/\upsilon )\,\rmd x_1\,\rmd x_2\,\rmd y_1\,\rmd y_2\\
  \label{eq:D-k-is-Xi}
  &=\int_{\Space{R}{4n}}
    k(x_1,y_1,\upsilon y_2',x_2'/\upsilon )\,
    \widetilde{\Xi}_\myhbar(0,x_1,x_2',y_1,y_2')\,\rmd x_1\,\rmd x_2'\,\rmd
    y_1\,\rmd y_2' \,.
\end{align}
We denote by \(B\) the above change of variables in the integration kernel:\footnote{We do not overload the notation \(B\) by the parameter \(\upsilon\) since it will be redundant later, cf. Rem.~\ref{re:upsilon-vanish}.}
\begin{equation}
  \label{eq:B-definition-R4n}
  B: \Space{R}{4n} \rightarrow \Space{R}{4n}: \  (x_1,x_2,y_1,y_2)
  \mapsto (x_1,y_1,\upsilon y_2,x_2/\upsilon )\,.
\end{equation}
Its inverse map obviously is:
\begin{equation}
  \label{eq:B-definition-R4n-inverse}
  B^{-1}: \Space{R}{4n} \rightarrow \Space{R}{4n}: \  (x_1,x_2,y_1,y_2)
  \mapsto (x_1,\upsilon y_2,x_2,y_1/\upsilon )\,.
\end{equation}
Then, computation~\eqref{eq:D-k-is-Xi} implies:
\begin{lem}
  \label{le:2-conv-int-Xi}
  The two-sided convolution with a kernel \(k(x_1,y_1,x_2,y_2)\) is the
  integrated representation \(\widetilde{\Xi}\) preceded by \(B\)~\eqref{eq:B-definition-R4n}:
  \begin{equation}
    \label{eq:2-sided-int-Xi}
    \oper{D}(k)=\widetilde{\Xi}_\myhbar(k\circ B)\,.
  \end{equation}
\end{lem}

The uniqueness (up to unitary equivalence) of the irreducible unitary
representation of \(\Space{H}{2n}\) for a given Planck constant
\(\myhbar\) implies the following:
\begin{prop}
  The Schr\"odinger representation
  \(\uir{}{\myhbar}(g)\)~\eqref{eq:schroedinger-rep} on
  \(\FSpace{L}{2}(\Space{R}{2n})\) and the representation
  \(\widetilde{\Xi}_\myhbar\) are intertwined by a unitary operator \(U:
  \FSpace{L}{2}(\Space{R}{2n}) \rightarrow \FSpace{L}{2}(\Space{R}{2n})\):
  \begin{equation}
    \label{eq:intertwine-xi-schroed}
    \uir{}{\myhbar}(g)=U\circ \widetilde{\Xi}_\myhbar(g)\circ U^{-1}, \qquad
    \text{for all }g\in \Space{H}{2n}\,.
  \end{equation}
\end{prop}
\begin{proof}
  The composition of a symplectic automorphism \(A\) and the
  Schr\"odinger representation \(\uir{}{\myhbar}\) is intertwined
  with \(\uir{}{\myhbar}\) by the unitary operator \(U(A)\)---the
  metaplectic representation of the double cover of the symplectic
  group~\citelist{\cite{Neretin11a}*{Ch.~1}
    \cite{Folland89}*{Ch.~4}}. 
\end{proof}
It is known that operator \(U\) is the Gaussian integral
operator~\citelist{\cite{Neretin11a}*{Ch.~1}
  \cite{Folland89}*{\S~4.4}}, but its explicit form will not be used
in this paper, merely its existence is important
for us. The proposition implies the same relation for the integrated
representations.
\begin{thm}
  \label{th:two-sided-PDO-symbol}
  The unitary operator \(U\)~\eqref{eq:intertwine-xi-schroed}
  conjugates 
  the operators of two-sided convolution \(\oper{D}(k)\) with the
  kernel \(k\) and the pseudodifferential operator \(a(D,X)\)   \begin{displaymath}
    U\circ \oper{D}(k)\circ U^{-1}=a(D,X),
  \end{displaymath}
  where the   Weyl symbol of \(a(D,X)\) is:
  \begin{equation}
    \label{eq:two-conv-PDO-symb}
    a(x_1,x_2,y_1,y_2)= \wideparen{k}(x_1,y_1,\upsilon y_2,x_2/\upsilon)\,,
  \end{equation}
  and \(\wideparen{k}\) is the symplectic Fourier transform of \(k\).
\end{thm}
\begin{proof}
  A combination of~\eqref{eq:2-sided-int-Xi} and the intertwining
  property~\eqref{eq:intertwine-xi-schroed} implies:
  \begin{align}
    \nonumber 
    U\circ\oper{D}(k)\circ U^{-1}
    &=\int_{\Space{R}{4n}}
      k(x_1,y_1,\upsilon y_2,x_2/\upsilon )\, \uir{}{\myhbar}(0,x_1,x_2,y_1,y_2)\,\rmd
      x_1\,\rmd x_2\,\rmd y_1\,\rmd y_2  \\
    \label{eq:PDO-intgrated}
    &=\int_{\Space{R}{4n}}
      \wideparen{a}(x_1,x_2,y_1,y_2)\, \uir{}{\myhbar}(0,x_1,x_2,y_1,y_2)\,\rmd x_1\,\rmd x_2\,\rmd y_1\,\rmd y_2 \,.
  \end{align}
  where \(\wideparen{a}(x_1,x_2,y_1,y_2)=k(x_1,y_1,\upsilon
  y_2,x_2/\upsilon )\) and~\eqref{eq:two-conv-PDO-symb} follows.
\end{proof}
We can now review Example~\ref{ex:as-2-side-conv}.
\begin{example}
  \begin{enumerate}
  \item   For the kernel \(k\)~\eqref{eq:two-kernel-multipl} of operator of
    multiplication by \(\psi\), we obtain
    \(\wideparen{k}(x_1,y_1,x_2,y_2)=\psi(\half(y_1+y_2),\half(x_1+x_2))\). Thus
    \begin{align*}
      a(x_1,x_2,y_1,y_2)&= \wideparen{k}(x_1,y_1,\tau y_2, x_2/\tau)\\
                        &=\psi(\half(y_1+x_2/\tau),\half(x_1 +\tau y_2))\,.
    \end{align*}
  \item Integrated representations \(\Lambda_\myhbar(k)\) and
    \(R_\myhbar(k)\), which are two-sided convolutions
    \(\oper{D}(k\otimes \delta )\) and \(\oper{D}(\delta \otimes k )\),
    respectively, are unitary equivalent to PDO's with symbols
    \begin{displaymath}
      a(x_1,x_2,y_1,y_2)=\wideparen{k}(x_1, y_1)\quad \text{and} \quad
      a(x_1,x_2,y_1,y_2)=\wideparen{k}(\tau y_2,x_2/\tau ),
    \end{displaymath}
    respectively.
  \item For \(\Phi_{\tau\varsigma}(x,y)=\scalar{\phi_\varsigma}{\uir{}{\myhbar}(x,y)\phi_\tau}\)~\eqref{eq:Phi-tau-sigma}, Lem.~\ref{le:intertwining-kernel-properties} states
    \(\wideparen{{\Phi}}_{\tau\varsigma}={\Phi}_{\tau\varsigma}\). Also,
    \({\Phi}_{\varsigma\tau}(\upsilon y_2,x_2/\upsilon )=\Phi_{\varsigma\tau}(x_2,y_2)\) for \(\upsilon=\sqrt{\tau \varsigma}\) by inspection.  Then, the intertwining operator \(\FSpace[\tau]{F}{} \rightarrow \FSpace[\varsigma]{F}{}\) from Ex.~\ref{ex:intertwining-kernel} is unitary
    equivalent to PDO with the symbol
    \begin{equation}
      \label{eq:FSB-proj-as-PDO}
      a(x_1,x_2,y_1,y_2)=\Phi_{\tau\varsigma}(x_2,y_2)\,.
    \end{equation}
    In the special case \(\tau=\varsigma\) we note, that  in the decomposition
    \(\FSpace{L}{2}(\Space{R}{2n})= \FSpace{L}{2}(\Space{R}{n})\otimes
    \FSpace{L}{2}(\Space{R}{n})\), PDO with
    symbol~\eqref{eq:FSB-proj-as-PDO} acts as
    \(I\otimes P_\phi\)---the identity in the first component and
    one-dimensional projection on the subspace spanned by the Gaussian
    \(\phi_\tau\) on the second. A similar representation of the projection
    \(\FSpace{L}{2}(\Space{C}{n}) \rightarrow \FSpace{F}{\phi}\) was obtained in~\cite{Vasilevski99b}*{Thm.~2.1}
    by a chain of explicit integral transformations.
  \end{enumerate}
\end{example}

\subsection{Composition and twisted convolution}
\label{sec:comp-oper}

Now we can use the composition of integrated
representations~\eqref{eq:twisted-conv-defn} to establish a calculus
of two-sided representations. This replaces the composition formula~\eqref{eq:composition-two-sided-group} of two-sided convolutions on groups.
\begin{lem}
  \label{le:2-con-twist-comp}
  The composition \(\oper{D}(k_1) \oper{D}(k_2)\) of two-sided convolutions with kernels \(k_1\) and
  \(k_2\)  is a two-sided convolution  \(\oper{D}(k)\) with the kernel  \(k =k_1 \Btwist k_2\) obtained by a  twisted convolution type formula:
  \begin{align}
    \label{eq:composition-kernel-abstract}
    k_1 \Btwist k_2 \coloneqq \left( (k_1\circ B) \twist (k_2\circ B) \right)\circ B^{-1}. 
  \end{align}
  Explicitly it is given by:
  \begin{equation}
    \label{eq:composition-kernel}
    \begin{split}
    k (x_1',y_1',x_2',y_2') &=\int_{\Space{R}{4n}}k_1(x_1,y_1,x_2,y_2)\,k_2(x_1'-x_1,y_1'-y_1,x_2'-x_2,y_2'-y_2)\,\\
    &\qquad \qquad \times \rme^{\pi\myhbar\rmi (y_1x'_1+x_2y'_2-x_1y_1'-y_2x_2')}\,dx_1\,dx_2\,dy_1\,dy_2\,. 
  \end{split}
  \end{equation}
\end{lem}
\begin{rem}
  \label{re:upsilon-vanish}
  Note, that the explicit expression \eqref{eq:composition-kernel} is independent from the parameter \(\upsilon\), which formally participates in the abstract formula~\eqref{eq:composition-kernel-abstract} through the map \(B\)~\eqref{eq:B-definition-R4n}.
\end{rem}
\begin{proof}
  Using the relation~\eqref{eq:2-sided-int-Xi} in both directions we obtain:
  \begin{align*}
    \oper{D}(k_1) \, \oper{D}(k_2)&=\widetilde{\Xi}_\myhbar(k_1\circ B)\,\widetilde{\Xi}_\myhbar(k_2\circ B)\\
    &=\widetilde{\Xi}_\myhbar\left((k_1\circ B) \twist (k_2\circ B)\right)\\
    &=\oper{D}\left(((k_1\circ B) \twist (k_2\circ B))\circ B^{-1}\right)\,.
  \end{align*}
  Thus, for the composition kernel defined by the identity \(k\circ B=(k_1\circ B) \twist
  (k_2\circ B)\)  we compute: 
  \begin{align*}
    \lefteqn{(k\circ B)(x_1',y_1',x_2',y_2' )}\qquad \qquad&\\
    &=\int_{\Space{R}{4n}}k_1(x_1,y_1,\tau y_2,x_2/\tau
      )\,k_2(x_1'-x_1,y_1'-y_1,\tau(y_2'-y_2),(x_2'-x_2)/\tau )\,\\
  &\qquad \quad \times \rme^{\pi\myhbar\rmi (y_1x'_1+y_2x'_2-x_1y_1'-x_2y_2')}\,dx_1\,dx_2\,dy_1\,dy_2\,. 
  \end{align*}
  Applying transformation \(B^{-1}\)~\eqref{eq:B-definition-R4n-inverse} to the last expression we obtain~\eqref{eq:composition-kernel}.
\end{proof}
\begin{rem}
  \label{re-almost-twist-conv}
  Formula~\eqref{eq:composition-kernel} is different from the twisted
  convolution on the Heisenberg group \(\Space{H}{2n}\) by the signs
  in the exponent for terms with \(x_2\) or \(y_2\). Thus, the twisted convolution part for \((x_2,y_2)\) phase space is intertwined by complex conjugation. 
\end{rem}
\begin{prop}
  A Toeplitz operator  \(\oper{T}_\psi: f \mapsto \oper{P}_{\varsigma}
  (\psi f)\) considered as an operator \(\FSpace{L}{2}(\Space{C}{n}) \rightarrow \FSpace[\varsigma]{F}{}\) is a two-sided convolution with the kernel \(k^{\#}_\psi\) given by
  either of the following expressions:
  \begin{align}
    \nonumber 
    k^{\#}_\psi(x_1,y_1,x_2,y_2)
    &= (2/\myhbar)^{n}\, \wideparen{\psi}
      (2x_1,2y_1)\, \rme^{\pi\myhbar\rmi (y_1x_2-x_1y_2)}\,
      \Phi_\varsigma(x_2-x_1,y_2-y_1) \\
   \label{eq:Toeplitz-two-sided-conv-kernel}
    &= (2/\myhbar)^{n}\, \wideparen{\psi} (2x_1,2y_1)\,
      \overline{[\Lambda_\myhbar(x_1,y_1) \Phi_\varsigma] (x_2,y_2)}\\
   \label{eq:toeplitz-kernel-shifted}
    &= (2/\myhbar)^{n}\, \wideparen{\psi} (2x_1,2y_1)\,
      [{\Lambda_\myhbar(x_2,y_2) \Phi_\varsigma] (x_1,y_1)}\\
    \label{eq:toeplitz-kernel-integral}
    &= (2/\myhbar)^{n}\, \wideparen{\psi} (2x_1,2y_1)\,
      K_{\varsigma}(x_1,y_1;x_2,y_2)\,.
  \end{align}
  Also, a cross-Toeplitz operator  \(\oper{T}_\psi: f \mapsto \oper{P}_{\varsigma}
  (\psi \oper{P}_\tau f)\) considered as an operator \(\FSpace[\tau]{F}{} \rightarrow \FSpace[\varsigma]{F}{}\) is a two-sided convolution with the kernel:
  \begin{equation}
    \label{eq:toeplitz-kernel-reduced}
    k^{\#}(x_1, y_1, x_2, y_2)
    =  (2/\myhbar)^{n}\, \wideparen{\psi} (2x_1,2y_1)\,
      \overline{[\Lambda_\myhbar(x_1,y_1) \Phi_{\tau \varsigma}] (x_2,y_2)}
\,. 
  \end{equation}
\end{prop}
\begin{proof}
  For the kernel~\eqref{eq:fock-proj-two-conv} of the Bargmann
  projection and the kernel~\eqref{eq:two-kernel-multipl} of
  the multiplication operator we use~\eqref{eq:composition-kernel} to
  calculate the kernel of the composition:
  \begin{align}
    \nonumber 
    k^{\#}_\psi(x'_1,y'_1,x'_2,y'_2)
    &= (2/\myhbar)^{n}
    \int_{\Space{R}{4n}}
      \delta(x_1,y_1)\,  \Phi_\varsigma(x_2,y_2) \, \wideparen{\psi} (2(x'_1-x_1),2(y'_1-y_1))\\
    \nonumber 
    & \quad \qquad \qquad \times
      \delta(x'_1-x'_2-(x_1-x_2),y'_1-y'_2-(y_1-y_2))\\
    \nonumber 
    & \quad \qquad \qquad \times \,\rme^{\pi\myhbar\rmi
      (y_1x_1'+x_2y_2'-x_1y_1'-y_2x'_2)}\,dx_1\,dx_2\,dy_1\,dy_2\\
    \nonumber 
    &= (2/\myhbar)^{n}
    \int_{\Space{R}{2n}}
      \Phi_\varsigma(x_2,y_2) \, \wideparen{\psi} (2x'_1,2y'_1)\\
    \nonumber 
    & \quad \qquad \qquad \times
      \delta(x'_1-x'_2+x_2,y'_1-y'_2+y_2) \,\rme^{\pi\myhbar\rmi
      (x_2y_2'-y_2x'_2)}\,dx_2\,dy_2\\
    \label{eq:Toeplitz-two-sided-kernel-first}
    &= (2/\myhbar)^{n}\,
       \wideparen{\psi} (2x'_1,2y'_1)\,\Phi_\varsigma(x'_2-x'_1,y'_2-y'_1) \, \rme^{\pi\myhbar\rmi
      (y'_1x'_2-x'_1y'_2)}
      \,.
  \end{align}

  That is the first form the kernel. Its comparison with~\eqref{eq:FSB-repro-kernel} gives
  formula~\eqref{eq:toeplitz-kernel-integral}. The
  expression~\eqref{eq:left-action-pulled} of \(\Lambda_\myhbar\)
  gives two other forms of the kernel.

  Composing~\eqref{eq:Toeplitz-two-sided-kernel-first} with the Bargmann projection \(\oper{P}_\tau\) on the right we obtain:
  \begin{align}
    \nonumber 
    k^{\#}_\psi(x'_1,y'_1,x'_2,y'_2)
    &= (2/\myhbar)^{n}
    \int_{\Space{R}{4n}}
       \wideparen{\psi} (2x_1,2y_1)\,\Phi_\varsigma(x_2-x_1,y_2-y_1) \, \rme^{\pi\myhbar\rmi
      (y_1x_2-x_1y_2)}\\
    \nonumber 
    & \quad \qquad \qquad \times  \delta(x'_1-x_1,y'_1-y_1)\,  \Phi_\tau(x'_2-x_2,y'_2-y_2)\\
    \nonumber 
    & \quad \qquad \qquad \times \rme^{\pi\myhbar\rmi
      (y_1x_1'+x_2y_2'-x_1y_1'-y_2x'_2)}
      \,dx_1\,dx_2\,dy_1\,dy_2\\
    \nonumber
        &= (2/\myhbar)^{n}
    \int_{\Space{R}{2n}}
       \wideparen{\psi} (2x'_1,2y'_1)\,\Phi_\varsigma(x'_1-x_2,y'_1-y_2)\, \Phi_\tau(x'_2-x_2,y'_2-y_2) \\
    \nonumber 
    & \quad \qquad \qquad \times 
     \rme^{\pi\myhbar\rmi
      (y'_1x_2-x'_1y_2+x_2y_2' -y_2x'_2)}
      \,dx_2\,dy_2\\
        \nonumber
    &= (2/\myhbar)^{n}
      \wideparen{\psi} (2x'_1,2y'_1)\,
      \rme^{\pi\myhbar\rmi (y_1x_2-x_1y_2)}\,
      \Phi_{\tau\varsigma} (x_1'-x_2',y_1'-y_2')\,,
  \end{align}
  which is equivalent to~\eqref{eq:toeplitz-kernel-reduced}.
\end{proof}
It is clear,  that for \(\tau=\varsigma\) relation~\eqref{eq:Toeplitz-two-sided-conv-kernel} is
an expanded version of~\eqref{eq:toeplitz-kernel-fourier}:
\begin{displaymath}
  k_\psi(x,y)=(\myhbar/2)^{n} k^{\#}_\psi(x,y,0,0)\,.
\end{displaymath}
On the other hand~\eqref{eq:toeplitz-kernel-reduced} extends~\eqref{eq:toeplitz-kernel-fourier} for the cross-Toeplitz operators.

\subsection{PDO from two-sided convolutions}
\label{sec:pdo-from-two}

The form of cross-Toeplitz operators as  two-sided convolution gives a new connection to PDO.
\begin{thm}
  A Toeplitz operator \(T_\psi: \FSpace{L}{2}(\Space{R}{2n}) \rightarrow \FSpace[\tau]{F}{}\) is unitary equivalent to PDO with the symbol:
  \begin{align}
    \label{eq:toeplitz-as-PDO}
    a_\psi^{\#}(x_1,x_2,y_1,y_2)
    &=
    \int_{\Space{R}{2n}} {\psi}(x,y)\,
      \Lambda_{\myhbar}(2x-x_1,2y-y_1)\, {\Phi}_\tau(x_2,y_2)
    \,\rmd x\,\rmd y \\
    \label{eq:toeplitz-as-PDO-transform}
    &=
     \int_{\Space{R}{2n}} {\psi}(x,y)\,
       \overline{\Lambda_{\myhbar}(x_2,y_2)\, {\Phi}_\tau(2x-x_1,2y-y_1)}
     \,\rmd x\,\rmd y 
      \,.
  \end{align}
\end{thm}
\begin{proof}
  We need to evaluate the symplectic Fourier transform
  of the kernel \(k=k^{\#}_\phi\)~\eqref{eq:toeplitz-kernel-integral} to be
  substituted in~\eqref{eq:two-conv-PDO-symb}. The transform can be
  explicitly calculated, this reduces to a chain of manipulations with Gauss-type
  integrals. However, it is more inspiring to split symplectic Fourier
  transform in \((x_1,y_1,x_2,y_2)\) into a composition of two
  commuting partial symplectic Fourier transforms in coordinates
  \((x_1,y_1)\) and \((x_2,y_2)\), we denote them by
  \(\wideparen{\ }^1\) and \(\wideparen{\ }^2\) respectively. For
  \(\wideparen{\ }^1\), a calculation identical
  to~\eqref{eq:Fourier-product-conv} produces the convolution type
  expression:
  \begin{align*}
        \wideparen{k}^1(x_1,y_1,x_2,y_2)
    &= \left(\wideparen{\psi} (2x_1,2y_1)\,
      K_{(x_2,y_2)}(x_1,y_1)\right)\wideparen{\ }^1 \\
    &= 
      \int_{\Space{R}{2n}} {\psi}(x,y)\,
      \wideparen{K_{(x_2,y_2)})}^1(x_1-2x,y_1-2y) \,\rmd
      x\,\rmd y\,.
  \end{align*}
  Now we shall use the following properties, which were already stated:
  \begin{enumerate}
  \item \(K_{(x_2,y_2)}=\overline{\Lambda_{\myhbar} (x_2,y_2) \Phi_\tau}\),
    cf.~\eqref{eq:reproducing-kernel};
  \item the symplectic Fourier transform propagates through the complex
    conjugation as follows:
    \((\overline{f})\wideparen{\ }(x,y)=\overline{\wideparen{f}(-x,-y)}\); 
  \item
    the symplectic Fourier transform commutes with left
    shifts~\eqref{eq:sympl-Fourier-intertwin-left}; and 
  \item the symplectic Fourier transform fixes the Gaussian:
    \(\wideparen{\Phi}_\tau(x,y)=\Phi_\tau (x,y)\),
    cf.~\eqref{eq:sympl-Fourier-Gaussian}.
  \end{enumerate}
  A combination of those produces:
  \begin{align*}
    \wideparen{K_{(x_2,y_2)})}^1(x_1,y_1)
    &=\overline{\Lambda_{\myhbar} (x_2,y_2) \wideparen{\Phi}_\tau(-x_1,-y_1)}\\
    &= \overline{\Lambda_{\myhbar} (x_2,y_2) {\Phi}_\tau(-x_1,-y_1)}\\
    &=\Lambda_{\myhbar}(-x_1,-y_1) {\Phi}_\tau(x_2,y_2)\,.
  \end{align*}
  Then, the symplectic Fourier transform of \(k\) as the composition
  of \(\ \wideparen{\ }^1\) and \(\ \wideparen{\ }^2\) is:
  \begin{align*}
    { \wideparen{k}(x_1,y_1,x_2,y_2)}
    & =  \left(\wideparen{k}^1(x_1,y_1,x_2,y_2)\right)\wideparen{\ }^2\\
     & =   
      \left(\int_{\Space{R}{2n}} {\psi}(x,y)\,
      \Lambda_{\myhbar}(2x-x_1,2y-y_1)\, {\Phi}_\tau(x_2,y_2)\,\rmd
      x\,\rmd y\right)\wideparen{\ }^2\\
      &=   
      \int_{\Space{R}{2n}} {\psi}(x,y)\,
      \Lambda_{\myhbar}(2x-x_1,2y-y_1)\, {\Phi}_\tau(x_2,y_2)\,\rmd
        x\,\rmd y
        \,,
  \end{align*}
  where we again use properties~\eqref{eq:sympl-Fourier-intertwin-left}
  and~\eqref{eq:sympl-Fourier-Gaussian} and linearity of the
  symplectic Fourier transform.  Thus, the PDO symbol in accordance
  with~\eqref{eq:two-conv-PDO-symb} is:
  \begin{align*}
    a_\psi^{\#}(x_1,x_2,y_1,y_2)
    &= \wideparen{k}(x_1,y_1, \tau y_2, x_2/\tau)\\
    &=   
      \int_{\Space{R}{2n}} {\psi}(x,y)\,
      \Lambda_{\myhbar}(2x-x_1,2y-y_1)\, {\Phi}_\tau(x_2,y_2 )\,\rmd
        x\,\rmd y\,,
  \end{align*}
  because \({\Phi}_\tau(\tau y_2, x_2/\tau )=\Phi_\tau(x_2,y_2)\). To obtain
  the second form we swap variables
  using~\eqref{eq:right-through-left} and  note that \(\Phi_\tau\) is an even
  function
  .
\end{proof}
In a similar way (or by CAS computation) we can demonstrate:
\begin{thm}
  The cross-Toeplitz operator \(T_\psi: \FSpace[\tau]{F}{} \rightarrow \FSpace[\varsigma]{F}{}\) is unitary equivalent to PDO with the symbol:
  \begin{align}
    \label{eq:cross-toeplitz-as-PDO}
    a_\psi^{\#}(x_1,x_2,y_1,y_2)
    &=
    \int_{\Space{R}{2n}} {\psi}(x,y)\,
      \Lambda_{\myhbar}(2x-x_1,2y-y_1)\, {\Phi}_{\tau\varsigma}(\varsigma y_2, x_2/\varsigma)
    \,\rmd x\,\rmd y \\
    \label{eq:cross-toeplitz-as-PDO-transform}
    &=
     \int_{\Space{R}{2n}} {\psi}(x,y)\,
       \overline{\Lambda_{\myhbar}(\varsigma y_2, x_2/\varsigma)\, {\Phi}_{\varsigma\tau}(2x-x_1,2y-y_1)}
     \,\rmd x\,\rmd y 
      \,.
  \end{align}
\end{thm}
\begin{rem}
  We note that \(a_\psi^{\#}\) contains all information from the
  Guillemin symbol \(a_\psi\)~\eqref{eq:Toeplitz-PDO-diffused} since
  \(a_\psi(x,y)=\myhbar^{2n} a_\psi^{\#}(x,0,y,0)\). However, the
  important difference of the obtained
  formula~\eqref{eq:toeplitz-as-PDO} and the Guillemin
  result~\eqref{eq:Toeplitz-PDO-diffused} is that the new symbol map
  \(\psi\mapsto a_\psi^{\#}\) is based on the invertible
  transformation. In fact, \eqref{eq:toeplitz-as-PDO} is a case of
  quadratic Fourier transform (aka short time Fourier transform) with
  the Gaussian window~\cite{Neretin11a}*{\S~4.3.5}.
\end{rem}
To confirm the invertibility of the symbol map
\(\psi\mapsto a_\psi^{\#}\) we
re-write~\eqref{eq:reproducing-projection} in order to express a
similarity with~\eqref{eq:toeplitz-as-PDO-transform} more explicitly:
\begin{align}
  \nonumber 
  a_\psi^{\#}(x_1,x_2,y_1,y_2)
  & =
     \int_{\Space{R}{2n}} {\psi}(x,y)\,
       \overline{\Lambda_{\myhbar}(\varsigma y_2, x_2/\varsigma)\, {\Phi}_{\varsigma\tau}(2x-x_1,2y-y_1)}
     \,\rmd x\,\rmd y  \\
  \nonumber 
  &=
    \int_{\Space{R}{2n}} {\psi}(x,y)\, \overline{[\Lambda _{\myhbar}(s,x_l,y_l)\otimes R(0,x_r,y_r)  \Phi_2](x,y)}
    \,\rmd x\,\rmd y\\
  \nonumber 
  &=
  \rme^{ \pi \rmi \myhbar({x_1  x_2 }/{\varsigma}- \varsigma y_1 y_2)/2} \\
    \label{eq:toeplitz-as-PDO-explicit}
  &\qquad \times
    \int_{\Space{R}{2n}} {\psi}(x,y)\, \overline{[\widetilde{\Xi} _{\myhbar}(0,x_l,\varsigma y_r,y_l, x_r/\varsigma)  \Phi_2](x,y)}
    \,\rmd x\,\rmd y\,,
\end{align}
where \(\Phi_2(x,y)=\Phi_{\varsigma\tau}(2x,2y)=\overline{\Phi_{\tau\varsigma}(2x,2y)}\) and:
\begin{align}
  \label{eq:Toeplitz-PDO-point-1}
  (s,x_l,x_r,y_l,y_r)
  &=
    \textstyle    \frac{1}{4}
({x_1  x_2 }/{\varsigma}- \varsigma y_1 y_2 , x_1+5  \varsigma y_2, -x_1+3  \varsigma y_2, 5 {x_2}/{\varsigma}+y_1, 3 {x_2}/{\varsigma}-y_1)\,,\\
  \label{eq:Toeplitz-PDO-point-2}
  (x_1,x_2,y_1,y_2)
  &=
    \textstyle
    \frac{1}{8} (3x_l-5x_r, \varsigma(y_l  + y_r), 3y_l - 5 y_r, (x_l + x_r)/\varsigma)\,.
\end{align}

From~\eqref{eq:toeplitz-as-PDO-explicit} we can directly deduce the following 
\begin{cor}
 PDO  symbol \(a_\psi^{\#}\)~\eqref{eq:toeplitz-as-PDO-transform} of the Toeplitz operator
  \(\oper{T}_\psi\) is expressed through the FSB transform of \(\psi\):
  \begin{equation}
    \begin{split}
      \label{eq:Toeplitz-PDO-as-wavelet-two-sided}
      a_\psi^{\#}(x_1,x_2,y_1,y_2)&= \rme^{\pi \rmi \myhbar (x_1 y_2 -y_1 x_2)} [\oper{W}^{\widetilde{\Xi}_\myhbar}_{\Phi_2}\psi](g) \\
      &= \rme^{\pi \rmi \myhbar (x_1 y_2 -y_1 x_2)} \scalar{\psi}{\widetilde{\Xi}_\myhbar(g)  \Phi_2}\,, 
    \end{split}
  \end{equation}
  where
  \begin{displaymath}
    g= \frac{1}{4}
    (0, x_1+5  \varsigma y_2, 3 {x_2}-{\varsigma}y_1, 5 {x_2}/{\varsigma}+y_1, -x_1/\varsigma+3   y_2)
    \,.
  \end{displaymath}
\end{cor}
This allows us to convert our knowledge of (pre-)FSB transform to
information on the symbol map \(\psi\mapsto a_\psi^{\#}\).

\subsection{Characteriseing the kernels of cross-Toeplitz operators}
\label{sec:char-kern-cross}

Using the annihilation properties for the mixed Gaussian from Lem.~\ref{le:intertwining-kernel-properties} we observe from~\eqref{eq:toeplitz-kernel-reduced} that kernels of a two-sided convolution generated by a cross-Toeplitz operator vanished under some derived representation.  Obviously,
\begin{displaymath}
  \ladder[R,\tau]{+}[2] k^{\#}_\psi(x_1,y_1,x_2,y_2) = 0\,,
\end{displaymath}
where the index \(2\) above the ladder operator indicates that it acts on the pair \((x_2,y_2)\). Also, the identity \(\ladder[\Lambda,\varsigma]{-} \Phi_{\tau\varsigma}=0\) implies that \(k^{\#}_\psi(x_1,y_1,x_2,y_2)\) is annihilated under the family of operators
\begin{align}
  \label{eq:kernel-second-annihilator}
  \overset{2}{\Lambda}_{\myhbar}(x_1,y_1) \circ \ladder[\Lambda,\varsigma]{-}[2] \circ \overset{2}{\Lambda}_{\myhbar}(-x_1,-y_1) &= \ladder[\Lambda,\varsigma]{-}[2] - \ladder[\Lambda,\varsigma]{-}[1] - \ladder[R,\varsigma]{-}[1]\\
  \nonumber 
  & = \ladder[\Lambda,\varsigma]{-}[2] + 2\overline{z}_1I \,,
\end{align}
for each fixed \((x_1,y_1)\in\Space{R}{2n}\) (with the above agreement on the numeric supersets). Obviously, an operator from the family~\eqref{eq:kernel-second-annihilator} is invariant under the right shift in \((x_2,y_2)\).
Therefore we are naturally arriving to the following definition.
\begin{defn}
  A kernel \(k(x_1,y_1,x_2,y_2)\) on \(\Space{R}{4n}\) is of \emph{cross-Toeplitz type} if it is annihilated by
  \begin{enumerate}
  \item the  left-invariant ladder operator \(\ladder[R,\tau]{+}[2]\); and
  \item the right-invariant operators \(\ladder[\Lambda,\varsigma]{-}[2] - \ladder[\Lambda,\varsigma]{-}[1] - \ladder[R,\varsigma]{-}[1]\), which is equal to:
  \begin{align*}
    \Lambda_{\myhbar,2}(x_1,y_1) \circ \ladder[\Lambda,\varsigma]{-}[2] \circ \Lambda_{\myhbar,2}(-x_1,-y_1) = \ladder[\Lambda,\varsigma]{-}[2] + 2\overline{z}_1I \,,
  \end{align*}
  where \(z_1=\sqrt{\frac{\modulus{\myh}}{2\tau }}(x_1+\rmi \varsigma y_1)\, \)
  .
  \end{enumerate}
\end{defn}
\begin{lem}
  Any differentiable cross-Toeplitz type kernel \(k\) is associated to a symbol \(\psi\) of Toeplitz operator for
  \begin{equation}
    \label{eq:Toeplitz-symbol-kernel}
    \wideparen{\psi}(2x_1,2y_1) = \frac {k(x_1,y_1,x_2,y_2) }{    [\Lambda_\myhbar(x_1,y_1) \Phi_{\tau\varsigma}] (x_2,y_2)}.
  \end{equation}
\end{lem}
\begin{proof}
  It is enough to show that the right-hand side of~\eqref{eq:Toeplitz-symbol-kernel} is independent of variables \((x_2,y_2)\).
  Let \(k_{(x_1,y_1)}(x,y)= k(x_1,y_1,x,y)\) and \( \Phi_{(x_1,y_1)} = \Lambda_\myhbar(x_1,y_1) \Phi_{\tau\varsigma}\). 
  Since \(\ladder[R,\tau]{+}k_{(x_1,y_1)} =0 \) and \(\ladder[R,\tau]{+}\Phi_{(x_1,y_1)} =0 \), we verify (omitting the parameter \(\tau\) of the complexification) that:
  \begin{align}
    \nonumber 
    0&=\frac{\ladder[R]{+}k_{(x_1,y_1)}}{\Phi_{(x_1,y_1)}} - \frac{k_{(x_1,y_1)} \cdot \ladder[R]{+}\Phi_{(x_1,y_1)}}{\Phi_{(x_1,y_1)}^2} \\
    \nonumber 
     &=\frac{({z}+\overline{\partial}_z)k_{(x_1,y_1)}}{\Phi_{(x_1,y_1)}} - \frac{k_{(x_1,y_1)} \cdot ({z}+\overline{\partial}_z)\Phi_{(x_1,y_1)}}{\Phi_{(x_1,y_1)}^2} \\
    \nonumber 
     &=\frac{\overline{\partial}_zk_{(x_1,y_1)}}{\Phi_{(x_1,y_1)}} - \frac{k_{(x_1,y_1)} \cdot \overline{\partial}_z\Phi_{(x_1,y_1)}}{\Phi_{(x_1,y_1)}^2} \\
    \label{eq:kernel-zero-derivative-bar}
     &=\overline{\partial}_z \left(\frac{k_{(x_1,y_1)}}{\Phi_{(x_1,y_1)}} \right).
  \end{align}
  Similarly, from the second condition on  right-invariant operators on cross-Toeplitz kernels we obtain that
  \begin{equation}
    \label{eq:kernel-zero-derivative}
    {\partial}_z \left(\frac{k_{(x_1,y_1)}}{\Phi_{(x_1,y_1)}} \right)=0,
  \end{equation}
  with the complexification for the parameter \(\varsigma\). Although, partial derivatives \(\overline{\partial}_z\) in~\eqref{eq:kernel-zero-derivative-bar} and \({\partial}_z\) in~\eqref{eq:kernel-zero-derivative} have different parameters of complexification we can conclude that \(k_{(x_1,y_1)}/\Phi_{(x_1,y_1)}\) is a constant for every fixed \((x_1,y_1)\).
\end{proof}

\begin{rem}
  Many results of this section remain true (with necessary adjustments) for the general localisation operators~\eqref{eq:localization-operator-general}. For example, the kernel of a two-sided convolution representing \(\oper{L}_\psi = \contravar{\Theta_2 } \circ \psi I \circ \covar{\Theta_1}\)~\eqref{eq:localization-operator-general} is cf.~\eqref{eq:toeplitz-kernel-reduced}:
  \begin{align*}
    k^{\#}(x_1. y_1, x_2, y_2)
    ,&=(2/\myhbar)^{n}\, \wideparen{\psi} (2x_1,2y_1)\,
       [\Lambda_\myhbar(x_1,y_1) \Theta] (x_2,y_2)
  \end{align*}
  where \(\Theta(x,y)=\scalar{\theta_1}{\uir{}{}(x,y) \theta_2}=\covar{\oper{R}{\Theta }_1}\Theta_2(x,y)\).
\end{rem}

\section*{Acknowledgments}
\label{sec:acknowledgments}
I am grateful to Prof.~N.L.~Vasilevski for the suggestion two-sided
convolutions as a research topic.  Prof.~Lewis Coburn provided interesting insights into the structure of FSB spaces. Author is also grateful to
Prof.~S.M.\; Sitnik for helpful discussions of Gauss-type integral
operators. Sitnik's visit to Leeds was supported by
London Mathematical  Society grant. An anonymous referee carefully read and thoroughly commented the initial version of the manuscript.

\small

\begin{bibdiv}
\begin{biblist}

\bib{AbreuFaustino15a}{article}{
      author={Abreu, Lu\'{\i}s~Daniel},
      author={Faustino, Nelson},
       title={On {T}oeplitz operators and localization operators},
        date={2015},
        ISSN={0002-9939},
     journal={Proc. Amer. Math. Soc.},
      volume={143},
      number={10},
       pages={4317\ndash 4323},
         url={https://doi.org/10.1090/proc/12211},
      review={\MR{3373930}},
}

\bib{AlameerKisil21a}{article}{
      author={Al~Ameer, Amerah~A.},
      author={Kisil, Vladimir~V.},
       title={Tuning co- and contra-variant transforms: the {Heisenberg} group
  illustration},
        date={2022},
     journal={SIGMA},
      volume={18},
      number={065},
       pages={21 p.},
        note={\arXiv{2105.13811}},
}

\bib{Alamer19a}{thesis}{
      author={Al~Ameer, Amerah~Ahmad},
       title={Singularities of analytic functions and group representations},
        type={Ph.D. Thesis},
     address={University of Leeds},
        date={2019},
        note={\url{https://etheses.whiterose.ac.uk/24776/}},
}

\bib{AliAntGaz14a}{book}{
      author={Ali, Syed~Twareque},
      author={Antoine, Jean-Pierre},
      author={Gazeau, Jean-Pierre},
       title={Coherent states, wavelets, and their generalizations},
     edition={Second},
      series={Theoretical and Mathematical Physics},
   publisher={Springer, New York},
        date={2014},
        ISBN={978-1-4614-8534-6; 978-1-4614-8535-3},
         url={http://dx.doi.org/10.1007/978-1-4614-8535-3},
      review={\MR{3154614}},
}

\bib{AlmalkiKisil18a}{article}{
      author={Almalki, Fadhel},
      author={Kisil, Vladimir~V.},
       title={Geometric dynamics of a harmonic oscillator, arbitrary minimal
  uncertainty states and the smallest step 3 nilpotent {Lie} group},
        date={2019},
     journal={J. Phys. A: Math. Theor},
      volume={52},
       pages={025301},
         url={https://doi.org/10.1088/1751-8121/aaed4d},
        note={\arXiv{1805.01399}},
}

\bib{AlmalkiKisil19a}{article}{
      author={Almalki, Fadhel},
      author={Kisil, Vladimir~V.},
       title={Solving the {S}chr\"odinger equation by reduction to a
  first-order differential operator through a coherent states transform},
        date={2020},
        ISSN={0375-9601},
     journal={Phys. Lett. A},
      volume={384},
      number={16},
       pages={126330},
  url={http://www.sciencedirect.com/science/article/pii/S0375960120301341},
        note={\arXiv{1903.03554}},
}

\bib{Bargmann61}{article}{
      author={Bargmann, V.},
       title={On a {Hilbert} space of analytic functions and an associated
  integral transform. {Part I}},
        date={1961},
     journal={Comm. Pure Appl. Math.},
      volume={3},
       pages={215\ndash 228},
}

\bib{BauerCoburnIsralowitz10a}{article}{
      author={Bauer, W.},
      author={Coburn, L.~A.},
      author={Isralowitz, J.},
       title={Heat flow, {BMO}, and the compactness of {Toeplitz} operators},
        date={2010},
        ISSN={0022-1236},
     journal={J. Funct. Anal.},
      volume={259},
      number={1},
       pages={57\ndash 78},
         url={http://dx.doi.org/10.1016/j.jfa.2010.03.016},
      review={\MR{2610379 (2011g:47056)}},
}

\bib{Berezin72}{article}{
      author={Berezin, F.~A.},
       title={Covariant and contravariant symbols of operators},
        date={1972},
     journal={Izv. Akad. Nauk SSSR Ser. Mat.},
      volume={36},
       pages={1134\ndash 1167},
        note={Reprinted in~\cite[pp.~228--261]{Berezin86}},
      review={\MR{50 \#2996}},
}

\bib{Berezin86}{book}{
      author={Berezin, F.~A.},
       title={Metod vtorichnogo kvantovaniya},
     edition={Second},
   publisher={``Nauka''},
     address={Moscow},
        date={1986},
        note={Edited and with a preface by M. K. Polivanov},
      review={\MR{89c:81001}},
}

\bib{BergerCoburn86b}{article}{
      author={Berger, C.~A.},
      author={Coburn, L.~A.},
       title={A symbol calculus for {T}oeplitz operators},
        date={1986},
        ISSN={0027-8424},
     journal={Proc. Nat. Acad. Sci. U.S.A.},
      volume={83},
      number={10},
       pages={3072\ndash 3073},
         url={https://doi.org/10.1073/pnas.83.10.3072},
      review={\MR{843308}},
}

\bib{BergerCoburn94a}{article}{
      author={Berger, C.~A.},
      author={Coburn, L.~A.},
       title={Heat flow and {Berezin}-{Toeplitz} estimates},
        date={1994},
        ISSN={0002-9327},
     journal={Amer. J. Math.},
      volume={116},
      number={3},
       pages={563\ndash 590},
         url={https://doi.org/10.2307/2374991},
      review={\MR{1277446}},
}

\bib{BergCob87}{article}{
      author={Berger, C.A.},
      author={Coburn, L.A.},
       title={Toeplitz operators on the {Segal-Bargmann} space},
        date={1987},
     journal={Trans. Amer. Math. Soc.},
      volume={301},
      number={2},
       pages={813\ndash 829},
}

\bib{Berndt07a}{book}{
      author={Berndt, Rolf},
       title={Representations of linear groups: {An} introduction based on
  examples from physics and number theory},
   publisher={Vieweg, Wiesbaden},
        date={2007},
        ISBN={978-3-8348-0319-1},
      review={\MR{2340988}},
}

\bib{BindenharnGustafsonLoheLouckMilne84a}{incollection}{
      author={{Bindenharn}, L.C.},
      author={{Gustafson}, R.A.},
      author={{Lohe}, M.A.},
      author={{Louck}, J.D.},
      author={{Milne}, S.C.},
       title={{Special functions and group theory in theoretical physics.}},
    language={English},
        date={1984},
   booktitle={Special functions: Group theoretical aspects and applications},
      editor={Askey, R.~A.},
      editor={Koornwinder, T.~H.},
      editor={Schempp, W.},
      series={Mathematics and its Applications},
      volume={18},
   publisher={D. Reidel Publ. Co.},
       pages={129\ndash 162},
}

\bib{BoggiattoCorderoGrochenig04a}{article}{
      author={Boggiatto, Paolo},
      author={Cordero, Elena},
      author={Gr\"ochenig, Karlheinz},
       title={Generalized anti-{Wick} operators with symbols in distributional
  {Sobolev} spaces},
        date={2004},
        ISSN={0378-620X},
     journal={Integral Equations Operator Theory},
      volume={48},
      number={4},
       pages={427\ndash 442},
         url={https://doi.org/10.1007/s00020-003-1244-x},
      review={\MR{2047590}},
}

\bib{ChristensenOlafsson09a}{article}{
      author={Christensen, Jens~Gerlach},
      author={{\'O}lafsson, Gestur},
       title={Examples of coorbit spaces for dual pairs},
        date={2009},
        ISSN={0167-8019},
     journal={Acta Appl. Math.},
      volume={107},
      number={1--3},
       pages={25\ndash 48},
      review={\MR{MR2520008}},
}

\bib{Coburn99}{article}{
      author={Coburn, L.~A.},
       title={The measure algebra of the {Heisenberg} group},
        date={1999},
        ISSN={0022-1236},
     journal={J. Funct. Anal.},
      volume={161},
      number={2},
       pages={509\ndash 525},
         url={https://doi.org/10.1006/jfan.1998.3354},
      review={\MR{1674651}},
}

\bib{Coburn01a}{incollection}{
      author={Coburn, L.~A.},
       title={The {Bargmann} isometry and {Gabor}-{Daubechies} wavelet
  localization operators},
        date={2001},
   booktitle={Systems, approximation, singular integral operators, and related
  topics ({Bordeaux}, 2000)},
      series={Oper. Theory Adv. Appl.},
      volume={129},
   publisher={Birkh\"auser, Basel},
       pages={169\ndash 178},
      review={\MR{1882695}},
}

\bib{Coburn01b}{article}{
      author={Coburn, L.~A.},
       title={On the {Berezin}-{Toeplitz} calculus},
        date={2001},
        ISSN={0002-9939},
     journal={Proc. Amer. Math. Soc.},
      volume={129},
      number={11},
       pages={3331\ndash 3338},
         url={https://doi.org/10.1090/S0002-9939-01-05917-2},
      review={\MR{1845010}},
}

\bib{Coburn12a}{article}{
      author={Coburn, L.~A.},
       title={Berezin transform and {Weyl}-type unitary operators on the
  {Bergman} space},
        date={2012},
        ISSN={0002-9939},
     journal={Proc. Amer. Math. Soc.},
      volume={140},
      number={10},
       pages={3445\ndash 3451},
         url={http://dx.doi.org/10.1090/S0002-9939-2012-11440-6},
      review={\MR{2929013}},
}

\bib{CoburnIsralowitzLi11a}{article}{
      author={Coburn, L.~A.},
      author={Isralowitz, J.},
      author={Li, Bo},
       title={Toeplitz operators with {BMO} symbols on the {Segal}-{Bargmann}
  space},
        date={2011},
        ISSN={0002-9947},
     journal={Trans. Amer. Math. Soc.},
      volume={363},
      number={6},
       pages={3015\ndash 3030},
         url={http://dx.doi.org/10.1090/S0002-9947-2011-05278-5},
      review={\MR{2775796 (2012h:47049)}},
}

\bib{Coburn19a}{incollection}{
      author={Coburn, Lewis~A.},
       title={{Fock} space, the {Heisenberg} group, heat flow, and {Toeplitz}
  operators},
        date={2019},
   booktitle={Handbook of analytic operator theory},
      series={CRC Press/Chapman Hall Handb. Math. Ser.},
   publisher={CRC Press, Boca Raton, FL},
       pages={1\ndash 15},
      review={\MR{3966624}},
}

\bib{CorderoGrochenig}{article}{
      author={Cordero, Elena},
      author={Gr\"{o}chenig, Karlheinz},
       title={Symbolic calculus and {F}redholm property for localization
  operators},
        date={2006},
        ISSN={1069-5869},
     journal={J. Fourier Anal. Appl.},
      volume={12},
      number={4},
       pages={371\ndash 392},
         url={https://doi.org/10.1007/s00041-005-5077-7},
      review={\MR{2256930}},
}

\bib{CordobaFefferman78a}{article}{
      author={C\'{o}rdoba, Antonio},
      author={Fefferman, Charles},
       title={Wave packets and {F}ourier integral operators},
        date={1978},
        ISSN={0360-5302},
     journal={Comm. Partial Differential Equations},
      volume={3},
      number={11},
       pages={979\ndash 1005},
         url={https://doi.org/10.1080/03605307808820083},
      review={\MR{507783}},
}

\bib{Coutinho10a}{article}{
      author={Coutinho, S.~C.},
       title={A lost chapter in the pre-history of algebraic analysis:
  {W}hittaker on contact transformations},
        date={2010},
        ISSN={0003-9519},
     journal={Arch. Hist. Exact Sci.},
      volume={64},
      number={6},
       pages={665\ndash 706},
         url={https://doi.org/10.1007/s00407-010-0063-0},
      review={\MR{2735747}},
}

\bib{deGosson11a}{book}{
      author={de~Gosson, Maurice~A.},
       title={Symplectic methods in harmonic analysis and in mathematical
  physics},
      series={Pseudo-Differential Operators. Theory and Applications},
   publisher={Birkh\"auser/Springer Basel AG, Basel},
        date={2011},
      volume={7},
        ISBN={978-3-7643-9991-7},
         url={http://dx.doi.org/10.1007/978-3-7643-9992-4},
      review={\MR{2827662 (2012m:53175)}},
}

\bib{deGosson13a}{article}{
      author={de~Gosson, Maurice~A.},
       title={Quantum blobs},
        date={2013},
        ISSN={0015-9018},
     journal={Found. Phys.},
      volume={43},
      number={4},
       pages={440\ndash 457},
         url={https://doi.org/10.1007/s10701-012-9636-x},
      review={\MR{3031620}},
}

\bib{FeichtingerPap14a}{incollection}{
      author={Feichtinger, H.~G.},
      author={Pap, M.},
       title={Coorbit theory and {Bergman} spaces},
        date={2014},
   booktitle={Harmonic and complex analysis and its applications},
      series={Trends Math.},
   publisher={Birkh\"auser/Springer, Cham},
       pages={231\ndash 259},
         url={https://doi.org/10.1007/978-3-319-01806-5_4},
      review={\MR{3203102}},
}

\bib{FeichGroech89a}{article}{
      author={Feichtinger, Hans~G.},
      author={Gr{\"o}chenig, K.~H.},
       title={Banach spaces related to integrable group representations and
  their atomic decompositions. {I}},
        date={1989},
        ISSN={0022-1236},
     journal={J. Funct. Anal.},
      volume={86},
      number={2},
       pages={307\ndash 340},
         url={http://dx.doi.org/10.1016/0022-1236(89)90055-4},
      review={\MR{1021139 (91g:43011)}},
}

\bib{FeichGroech89b}{article}{
      author={Feichtinger, Hans~G.},
      author={Gr{\"o}chenig, K.~H.},
       title={Banach spaces related to integrable group representations and
  their atomic decompositions. {II}},
        date={1989},
        ISSN={0026-9255},
     journal={Monatsh. Math.},
      volume={108},
      number={2-3},
       pages={129\ndash 148},
         url={http://dx.doi.org/10.1007/BF01308667},
      review={\MR{1026614 (91g:43012)}},
}

\bib{FeichtingerGrochenig88a}{incollection}{
      author={Feichtinger, Hans~G.},
      author={Gr\"{o}chenig, Karlheinz},
       title={A unified approach to atomic decompositions via integrable group
  representations},
        date={1988},
   booktitle={Function spaces and applications ({L}und, 1986)},
      series={Lecture Notes in Math.},
      volume={1302},
   publisher={Springer, Berlin},
       pages={52\ndash 73},
         url={https://doi.org/10.1007/BFb0078863},
      review={\MR{942257}},
}

\bib{Fock04a}{book}{
      author={Fock, V.~A.},
       title={Selected works. {Q}uantum mechanics and quantum field theory},
   publisher={Chapman \& Hall/CRC, Boca Raton, FL},
        date={2004},
        ISBN={0-415-30002-9},
        note={Edited by L. D. Faddeev, L. A. Khalfin and I. V. Komarov},
      review={\MR{2083062}},
}

\bib{Folland89}{book}{
      author={Folland, Gerald~B.},
       title={Harmonic analysis in phase space},
      series={Annals of Mathematics Studies},
   publisher={Princeton University Press},
     address={Princeton, NJ},
        date={1989},
      volume={122},
        ISBN={0-691-08527-7; 0-691-08528-5},
      review={\MR{92k:22017}},
}

\bib{Folland16a}{book}{
      author={Folland, Gerald~B.},
       title={A course in abstract harmonic analysis},
     edition={Second},
      series={Textbooks in Mathematics},
   publisher={CRC Press, Boca Raton, FL},
        date={2016},
        ISBN={978-1-4987-2713-6},
      review={\MR{3444405}},
}

\bib{GerryKnight05a}{book}{
      author={Gerry, Christopher},
      author={},
      author={Knight, Peter~L.},
       title={Introductory quantum optics},
   publisher={Cambridge University Press},
     address={Cambridge},
        date={2005},
        ISBN={9780521527354},
         url={https://books.google.co.uk/books?id=CgByyoBJJwgC},
}

\bib{Glauber63a}{article}{
      author={Glauber, Roy~J.},
       title={Coherent and incoherent states of the radiation field},
        date={1963Sep},
     journal={Phys. Rev.},
      volume={131},
       pages={2766\ndash 2788},
         url={https://link.aps.org/doi/10.1103/PhysRev.131.2766},
}

\bib{Grochenig01a}{book}{
      author={Gr\"{o}chenig, Karlheinz},
       title={Foundations of time-frequency analysis},
      series={Applied and Numerical Harmonic Analysis},
   publisher={Birkh\"{a}user Boston, Inc., Boston, MA},
        date={2001},
        ISBN={0-8176-4022-3},
         url={https://doi.org/10.1007/978-1-4612-0003-1},
      review={\MR{1843717}},
}

\bib{GrossmannLoupiasStein68}{article}{
      author={Grossmann, A.},
      author={Loupias, G.},
      author={Stein, E.~M.},
       title={An algebra of pseudodifferential operators and quantum mechanics
  in phase space},
        date={1968},
        ISSN={0373-0956},
     journal={Ann. Inst. Fourier (Grenoble)},
      volume={18},
      number={fasc. 2},
       pages={343\ndash 368, viii (1969)},
         url={http://www.numdam.org/item?id=AIF_1968__18_2_343_0},
      review={\MR{0267425}},
}

\bib{Guillemin84}{article}{
      author={Guillemin, Victor},
       title={Toeplitz operator in $n$-dimensions},
        date={1984},
     journal={Integral Equations Operator Theory},
      volume={7},
       pages={145\ndash 205},
}

\bib{Hormander83IV}{book}{
      author={H\"{o}rmander, Lars},
       title={The analysis of linear partial differential operators. {IV}},
      series={Classics in Mathematics},
   publisher={Springer-Verlag, Berlin},
        date={2009},
        ISBN={978-3-642-00117-8},
         url={https://doi.org/10.1007/978-3-642-00136-9},
        note={Fourier integral operators, Reprint of the 1994 edition},
      review={\MR{2512677}},
}

\bib{Howe80a}{article}{
      author={Howe, Roger},
       title={On the role of the {Heisenberg} group in harmonic analysis},
        date={1980},
        ISSN={0002-9904},
     journal={Bull. Amer. Math. Soc. (N.S.)},
      volume={3},
      number={2},
       pages={821\ndash 843},
      review={\MR{81h:22010}},
}

\bib{Howe80b}{article}{
      author={Howe, Roger},
       title={Quantum mechanics and partial differential equations},
        date={1980},
        ISSN={0022-1236},
     journal={J. Funct. Anal.},
      volume={38},
      number={2},
       pages={188\ndash 254},
      review={\MR{83b:35166}},
}

\bib{Husimi40}{article}{
      author={{Husimi}, K.},
       title={Some formal properties of the density matrix},
    language={English},
        date={1940},
        ISSN={0370-1239; 2185-2707/e},
     journal={Proc. Phys.-Math. Soc. Japan, III. Ser.},
      volume={22},
       pages={264\ndash 314},
}

\bib{Johansson08a}{article}{
      author={Johansson, Andreas},
       title={Shift-invariant signal norms for fault detection and control},
        date={2008},
        ISSN={0167-6911},
     journal={Systems Control Lett.},
      volume={57},
      number={2},
       pages={105\ndash 111},
      review={\MR{MR2378755 (2009d:93035)}},
}

\bib{KermackMcCrea31a}{article}{
      author={{Kermack}, W.~O.},
      author={{McCrea}, W.~H.},
       title={An operational method for the solution of linear partial
  differential equations},
    language={English},
        date={1931},
        ISSN={0370-1646},
     journal={{Proc. R. Soc. Edinburgh}},
      volume={51},
       pages={176\ndash 189},
}

\bib{Kirillov76}{book}{
      author={Kirillov, A.~A.},
       title={Elements of the theory of representations},
   publisher={Springer-Verlag},
     address={Berlin},
        date={1976},
        note={Translated from the Russian by Edwin Hewitt, Grundlehren der
  Mathematischen Wissenschaften, Band 220},
      review={\MR{54 \#447}},
}

\bib{Kirillov04a}{book}{
      author={Kirillov, A.~A.},
       title={Lectures on the orbit method},
      series={Graduate Studies in Mathematics},
   publisher={American Mathematical Society},
     address={Providence, RI},
        date={2004},
      volume={64},
        ISBN={0-8218-3530-0},
      review={\MR{2069175 (2005c:22001)}},
}

\bib{Kisil92}{article}{
      author={Kisil, Vladimir~V.},
       title={Algebra of two-sided convolutions on the {H}eisenberg group},
        date={1992},
        ISSN={0869-5652},
     journal={Dokl. Akad. Nauk},
      volume={325},
      number={1},
       pages={20\ndash 23},
        note={Translated in \emph{Russ.Acad. of Sci Doklady, Math}, v. {\bf
  46}(1994), pp. 12--16. \MR{93j:22018}},
}

\bib{Kisil93b}{article}{
      author={Kisil, Vladimir~V.},
       title={On the algebra of pseudodifferential operators that is generated
  by convolutions on the {H}eisenberg group},
        date={1993},
        ISSN={0037-4474},
     journal={Sibirsk. Mat. Zh.},
      volume={34},
      number={6},
       pages={75\ndash 85},
        note={(Russian) \MR{95a:47053}},
}

\bib{Kisil94a}{article}{
      author={Kisil, Vladimir~V.},
       title={Local behavior of two-sided convolution operators with singular
  kernels on the {H}eisenberg group},
        date={1994},
        ISSN={0025-567X},
     journal={Mat. Zametki},
      volume={56},
      number={2},
       pages={41\ndash 55, 158},
        note={(Russian) \MR{96a:22015}},
}

\bib{Kisil93e}{article}{
      author={Kisil, Vladimir~V.},
       title={The spectrum of the algebra generated by two-sided convolutions
  on the {H}eisenberg group and by operators of multiplication by continuous
  functions},
        date={1994},
        ISSN={0869-5652},
     journal={Dokl. Akad. Nauk},
      volume={337},
      number={4},
       pages={439\ndash 441},
        note={Translated in \emph{Russ. Acad. of Sci Doklady, Math}, v. {\bf
  50}(1995), No 1, pp. 92--97. \MR{95j:22019}},
}

\bib{Kisil94f}{article}{
      author={Kisil, Vladimir~V.},
       title={Connection between two-sided and one-sided convolution type
  operators on a non-commutative group},
        date={1995},
        ISSN={0378-620X},
     journal={Integral Equations Operator Theory},
      volume={22},
      number={3},
       pages={317\ndash 332},
        note={\MR{96d:44004}},
}

\bib{Kisil96e}{article}{
      author={Kisil, Vladimir~V.},
       title={Local algebras of two-sided convolutions on the {H}eisenberg
  group},
        date={1996},
        ISSN={0025-567X},
     journal={Mat. Zametki},
      volume={59},
      number={3},
       pages={370\ndash 381, 479},
      review={\MR{MR1399963 (97h:22006)}},
}

\bib{Kisil98e}{article}{
      author={{Kisil}, Vladimir~V.},
       title={Harmonic analysis and localization technique},
        date={1998},
     journal={Odessa University Herald},
      volume={3},
       pages={60\ndash 63},
        note={\arXiv{math/9902012}},
}

\bib{Kisil97c}{article}{
      author={Kisil, Vladimir~V.},
       title={Analysis in {$\mathbf{R}\sp {1,1}$} or the principal function
  theory},
        date={1999},
        ISSN={0278-1077},
     journal={Complex Variables Theory Appl.},
      volume={40},
      number={2},
       pages={93\ndash 118},
        note={\arXiv{funct-an/9712003}},
      review={\MR{MR1744876 (2000k:30078)}},
}

\bib{Kisil94e}{article}{
      author={Kisil, Vladimir~V.},
       title={Relative convolutions. {I}. {P}roperties and applications},
        date={1999},
        ISSN={0001-8708},
     journal={Adv. Math.},
      volume={147},
      number={1},
       pages={35\ndash 73},
        note={\arXiv{funct-an/9410001},
  \href{http://www.idealibrary.com/links/doi/10.1006/aima.1999.1833}{On-line}.
  \Zbl{933.43004}},
      review={\MR{MR1725814 (2001h:22012)}},
}

\bib{Kisil98a}{article}{
      author={Kisil, Vladimir~V.},
       title={Wavelets in {B}anach spaces},
        date={1999},
        ISSN={0167-8019},
     journal={Acta Appl. Math.},
      volume={59},
      number={1},
       pages={79\ndash 109},
        note={\arXiv{math/9807141},
  \href{http://dx.doi.org/10.1023/A:1006394832290}{On-line}},
      review={\MR{MR1740458 (2001c:43013)}},
}

\bib{Kisil02e}{article}{
      author={Kisil, Vladimir~V.},
       title={{$p$}-{M}echanics as a physical theory: an introduction},
        date={2004},
        ISSN={0305-4470},
     journal={J. Phys. A: Math. Theor},
      volume={37},
      number={1},
       pages={183\ndash 204},
        note={\arXiv{quant-ph/0212101},
  \href{http://stacks.iop.org/0305-4470/37/183}{On-line}. \Zbl{1045.81032}},
      review={\MR{MR2044764 (2005c:81078)}},
}

\bib{Kisil10c}{article}{
      author={Kisil, Vladimir~V.},
       title={Covariant transform},
        date={2011},
     journal={Journal of Physics: Conference Series},
      volume={284},
      number={1},
       pages={012038},
         url={http://stacks.iop.org/1742-6596/284/i=1/a=012038},
        note={\arXiv{1011.3947}},
}

\bib{Kisil11c}{incollection}{
      author={Kisil, Vladimir~V.},
       title={{E}rlangen programme at large: {An} overview},
        date={2012},
   booktitle={Advances in applied analysis},
      editor={Rogosin, S.V.},
      editor={Koroleva, A.A.},
   publisher={Birkh\"auser Verlag},
     address={Basel},
       pages={1\ndash 94},
        note={\arXiv{1106.1686}},
}

\bib{Kisil10a}{article}{
      author={Kisil, Vladimir~V.},
       title={Hypercomplex representations of the {H}eisenberg group and
  mechanics},
        date={2012},
        ISSN={0020-7748},
     journal={Internat. J. Theoret. Phys.},
      volume={51},
      number={3},
       pages={964\ndash 984},
         url={http://dx.doi.org/10.1007/s10773-011-0970-0},
        note={\arXiv{1005.5057}. \Zbl{1247.81232}},
      review={\MR{2892069}},
}

\bib{Kisil12b}{article}{
      author={Kisil, Vladimir~V.},
       title={Operator covariant transform and local principle},
        date={2012},
     journal={J. Phys. A: Math. Theor.},
      volume={45},
       pages={244022},
        note={\arXiv{1201.1749}.
  \href{http://stacks.iop.org/1751-8121/45/244022}{On-line}},
}

\bib{Kisil13b}{article}{
      author={Kisil, Vladimir~V.},
       title={Boundedness of relative convolutions on nilpotent {L}ie groups},
        date={2013},
     journal={Zb. Pr. Inst. Mat. NAN Ukr. (Proc. Math. Inst. Ukr. Ac. Sci.)},
      volume={10},
      number={4--5},
       pages={185\ndash 189},
        note={\arXiv{1307.3882}},
}

\bib{Kisil09e}{article}{
      author={Kisil, Vladimir~V.},
       title={Induced representations and hypercomplex numbers},
        date={2013},
     journal={Adv. Appl. Clifford Algebras},
      volume={23},
      number={2},
       pages={417\ndash 440},
         url={http://dx.doi.org/10.1007/s00006-012-0373-1},
        note={\arXiv{0909.4464}. \Zbl{1269.30052}},
}

\bib{Kisil13a}{article}{
      author={Kisil, Vladimir~V.},
       title={Calculus of operators: {C}ovariant transform and relative
  convolutions},
        date={2014},
     journal={Banach J. Math. Anal.},
      volume={8},
      number={2},
       pages={156\ndash 184},
         url={\url{http://www.emis.de/journals/BJMA/tex_v8_n2_a15.pdf}},
        note={\arXiv{1304.2792},
  \href{http://www.emis.de/journals/BJMA/tex_v8_n2_a15.pdf}{on-line}},
}

\bib{Kisil12d}{article}{
      author={Kisil, Vladimir~V.},
       title={The real and complex techniques in harmonic analysis from the
  point of view of covariant transform},
        date={2014},
     journal={Eurasian Math. J.},
      volume={5},
       pages={95\ndash 121},
        note={\arXiv{1209.5072}.
  \href{http://emj.enu.kz/images/pdf/2014/5-1-4.pdf}{On-line}},
}

\bib{Kisil13c}{incollection}{
      author={Kisil, Vladimir~V.},
       title={Uncertainty and analyticity},
    language={English},
        date={2015},
   booktitle={Current trends in analysis and its applications},
      editor={Mityushev, Vladimir~V.},
      editor={Ruzhansky, Michael~V.},
      series={Trends in Mathematics},
   publisher={Springer International Publishing},
       pages={583\ndash 590},
         url={http://dx.doi.org/10.1007/978-3-319-12577-0_64},
        note={\arXiv{1312.4583}},
}

\bib{Kisil17a}{incollection}{
      author={Kisil, Vladimir~V.},
       title={Symmetry, geometry and quantization with hypercomplex numbers},
        date={2017},
   booktitle={Geometry, integrability and quantization {XVIII}},
      editor={Mladenov, Iva\"{\i}lo~M.},
      editor={Meng, Guowu},
      editor={Yoshioka, Akira},
   publisher={Bulgar. Acad. Sci., Sofia},
       pages={11\ndash 76},
        note={\arXiv{1611.05650}},
      review={\MR{3616912}},
}

\bib{Kisil21c}{article}{
      author={Kisil, Vladimir~V.},
       title={Metamorphism---an integral transform reducing the order of a
  differential equation},
        date={2021},
        note={\arXiv{2105.12079}},
}

\bib{Kisil22a}{article}{
      author={Kisil, Vladimir~V.},
       title={Mixed coherent states decompositions on split domains},
        date={2022},
     journal={in prepapation},
}

\bib{Mackey70a}{incollection}{
      author={Mackey, George~W.},
       title={Induced representations of locally compact groups and
  applications},
        date={1970},
   booktitle={Functional {Analysis} and {Related} {Fields} ({Proc}. {Conf}. for
  {M}. {Stone}, {Univ}. {Chicago}, {Chicago}, {I}ll., 1968)},
   publisher={Springer},
     address={New York},
       pages={132\ndash 166},
      review={\MR{0425010 (54 \#12968)}},
}

\bib{MazyaSchmidt07a}{book}{
      author={Maz'ya, Vladimir},
      author={Schmidt, Gunther},
       title={Approximate approximations},
      series={Mathematical Surveys and Monographs},
   publisher={American Mathematical Society, Providence, RI},
        date={2007},
      volume={141},
        ISBN={978-0-8218-4203-4},
         url={https://doi.org/10.1090/surv/141},
      review={\MR{2331734}},
}

\bib{Miheisi10a}{article}{
      author={Miheisi, Nazar},
       title={Convolution operators on {B}anach lattices with shift-invariant
  norms},
        date={2010},
        ISSN={0378-620X},
     journal={Integral Equations Operator Theory},
      volume={68},
      number={2},
       pages={287\ndash 299},
         url={https://doi.org/10.1007/s00020-010-1817-4},
      review={\MR{2721088}},
}

\bib{Moyal49}{article}{
      author={Moyal, J.~E.},
       title={Quantum mechanics as a statistical theory},
        date={1949},
     journal={Mathematical Proceedings of the Cambridge Philosophical Society},
      volume={45},
      number={1},
       pages={99–124},
}

\bib{Neretin11a}{book}{
      author={Neretin, Yurii~A.},
       title={Lectures on {Gaussian} integral operators and classical groups},
      series={EMS Series of Lectures in Mathematics},
   publisher={European Mathematical Society (EMS)},
     address={Z\"urich},
        date={2011},
        ISBN={978-3-03719-080-7},
         url={http://dx.doi.org/10.4171/080},
      review={\MR{2790054}},
}

\bib{Perelomov86}{book}{
      author={Perelomov, A.},
       title={Generalized coherent states and their applications},
      series={Texts and Monographs in Physics},
   publisher={Springer-Verlag},
     address={Berlin},
        date={1986},
        ISBN={3-540-15912-6},
      review={\MR{87m:22035}},
}

\bib{RottensteinerRuzhansky20a}{article}{
      author={Rottensteiner, David},
      author={Ruzhansky, Michael},
       title={L'oscillateur harmonique sur le groupe de {H}eisenberg},
    language={English},
        date={2020},
        ISSN={1631-073X},
     journal={C. R., Math., Acad. Sci. Paris},
      volume={358},
      number={5},
       pages={609\ndash 614},
}

\bib{Segal60}{book}{
      author={Segal, Irving~E.},
       title={Mathematical problems of relativistic physics},
      series={Proceedings of the Summer Seminar (Boulder, Colorado, 1960)},
   publisher={American Mathematical Society},
     address={Providence, R.I.},
        date={1963},
      volume={II},
}

\bib{Street08a}{article}{
      author={Street, Brian},
       title={An algebra containing the two-sided convolution operators},
        date={2008},
        ISSN={0001-8708},
     journal={Adv. Math.},
      volume={219},
      number={1},
       pages={251\ndash 315},
         url={http://dx.doi.org/10.1016/j.aim.2008.04.014},
      review={\MR{2435424 (2009h:43005)}},
}

\bib{MTaylor81}{book}{
      author={Taylor, Michael~E.},
       title={Pseudodifferential operators},
      series={Princeton Mathematical Series},
   publisher={Princeton University Press},
     address={Princeton, N.J.},
        date={1981},
      volume={34},
        ISBN={0-691-08282-0},
      review={\MR{82i:35172}},
}

\bib{MTaylor86}{book}{
      author={Taylor, Michael~E.},
       title={Noncommutative harmonic analysis},
      series={Mathematical Surveys and Monographs},
   publisher={American Mathematical Society},
     address={Providence, RI},
        date={1986},
      volume={22},
        ISBN={0-8218-1523-7},
      review={\MR{88a:22021}},
}

\bib{TurbinerVasilevski21a}{article}{
      author={Turbiner, Alexander~V.},
      author={Vasilevski, Nikolai},
       title={Poly-analytic functions and representation theory},
        date={2021},
        ISSN={1661-8254},
     journal={Complex Anal. Oper. Theory},
      volume={15},
      number={7},
       pages={Paper No. 110, 24},
         url={https://doi.org/10.1007/s11785-021-01154-y},
        note={\arXiv{2103.12771}},
      review={\MR{4313586}},
}

\bib{Vasilevski99b}{incollection}{
      author={Vasilevski, N.~L.},
       title={Poly-{Fock} spaces},
        date={2000},
   booktitle={Differential operators and related topics, {V}ol. {I} ({Odessa},
  1997)},
      series={Oper. Theory Adv. Appl.},
      volume={117},
   publisher={Birkh\"auser, Basel},
       pages={371\ndash 386},
      review={\MR{1764974}},
}

\bib{Vasilevski94a}{article}{
      author={Vasilevski, Nikolai~L.},
       title={Convolution operators on standard {CR}-manifolds. {II}.
  {Algebras} of convolution operators on the {Heisenberg} group},
        date={1994},
        ISSN={0378-620X},
     journal={Integral Equations Operator Theory},
      volume={19},
      number={3},
       pages={327\ndash 348},
         url={https://doi.org/10.1007/BF01203669},
      review={\MR{1280127}},
}

\bib{VasTru88}{article}{
      author={Vasilevski, Nikolai~L.},
      author={Trujillo, Rafael},
       title={Group convolutions on standard {Cauchy--Riemann} manifolds},
         how={Preprint 17, Ser.~B},
        date={1988},
        note={(Russian)},
}

\bib{VasTru94}{article}{
      author={Vasilevski, Nikolai~L.},
      author={Trujillo, Rafael},
       title={Group convolution operators on standard {CR}-manifolds. {I}.
  {Structural} properties},
        date={1994},
        ISSN={0378-620X},
     journal={Integral Equations Operator Theory},
      volume={19},
      number={1},
       pages={65\ndash 104},
         url={https://doi.org/10.1007/BF01202291},
      review={\MR{1271239}},
}

\bib{Vilenkin68}{book}{
      author={Vilenkin, N.~Ja.},
       title={Special functions and the theory of group representations},
      series={Translations of Mathematical Monographs},
   publisher={American Mathematical Society},
     address={Providence, R. I.},
        date={1968},
      volume={22},
        note={Translated from the Russian by V. N. Singh},
      review={\MR{37 \#5429}},
}

\bib{VilenkinKlimyk95a}{book}{
      author={Vilenkin, N.~Ja.},
      author={Klimyk, A.~U.},
       title={Representation of {Lie} groups and special functions},
      series={Mathematics and its Applications},
   publisher={Kluwer Academic Publishers Group, Dordrecht},
        date={1995},
      volume={316},
        ISBN={0-7923-3210-5},
         url={https://doi.org/10.1007/978-94-017-2885-0},
        note={Recent advances, Translated from the Russian manuscript by V. A.
  Groza and A. A. Groza},
      review={\MR{1371383}},
}

\bib{Zachos02a}{article}{
      author={Zachos, Cosmas},
       title={Deformation quantization: Quantum mechanics lives and works in
  phase-space},
        date={2002},
        ISSN={0217-751X},
     journal={Internat. J. Modern Phys. A},
      volume={17},
      number={3},
       pages={297\ndash 316},
        note={\arXiv{hep-th/0110114}},
      review={\MR{1 888 937}},
}

\bib{Zhu11a}{article}{
      author={Zhu, Kehe},
       title={Invariance of {Fock} spaces under the action of the {Heisenberg}
  group},
        date={2011},
        ISSN={0007-4497},
     journal={Bull. Sci. Math.},
      volume={135},
      number={5},
       pages={467\ndash 474},
         url={https://doi.org/10.1016/j.bulsci.2011.04.002},
      review={\MR{2817458}},
}

\bib{Zhu12a}{book}{
      author={Zhu, Kehe},
       title={Analysis on {Fock} spaces},
      series={Graduate Texts in Mathematics},
   publisher={Springer, New York},
        date={2012},
      volume={263},
        ISBN={978-1-4419-8800-3},
         url={https://doi.org/10.1007/978-1-4419-8801-0},
      review={\MR{2934601}},
}

\end{biblist}
\end{bibdiv}
\providecommand{\noopsort}[1]{} \providecommand{\printfirst}[2]{#1}
  \providecommand{\singleletter}[1]{#1} \providecommand{\switchargs}[2]{#2#1}
  \providecommand{\irm}{\textup{I}} \providecommand{\iirm}{\textup{II}}
  \providecommand{\vrm}{\textup{V}} \providecommand{\cprime}{'}
  \providecommand{\eprint}[2]{\texttt{#2}}
  \providecommand{\myeprint}[2]{\texttt{#2}}
  \providecommand{\arXiv}[1]{\myeprint{http://arXiv.org/abs/#1}{arXiv:#1}}
  \providecommand{\doi}[1]{\href{http://dx.doi.org/#1}{doi:
  #1}}\providecommand{\CPP}{\texttt{C++}}
  \providecommand{\NoWEB}{\texttt{noweb}}
  \providecommand{\MetaPost}{\texttt{Meta}\-\texttt{Post}}
  \providecommand{\GiNaC}{\textsf{GiNaC}}
  \providecommand{\pyGiNaC}{\textsf{pyGiNaC}}
  \providecommand{\Asymptote}{\texttt{Asymptote}}

\end{document}